\documentclass{amsart}
\usepackage{amssymb,amsmath,amsthm,amsrefs}
\usepackage{bm}
\usepackage{fullpage}
\usepackage{enumerate}
\usepackage{verbatim}
\usepackage{graphicx,float}
\usepackage{comment}

\newtheorem{theorem}[equation]{Theorem}
\newtheorem{lemma}[equation]{Lemma}
\newtheorem{prop}[equation]{Proposition}

\newtheorem{defn}[equation]{Definition}
\newtheorem{remark}[equation]{Remark}

\numberwithin{equation}{section}

\newcommand{\Grp}{\mathcal{G}}
\newcommand{\GrpK}{\Grp_{cat}}
\newcommand{\GrpP}[1]{\Grp_#1}

\newcommand{\pre}[1]{\mathring{#1}}

\newcommand{\R}{\mathbb{R}}
\newcommand{\Sph}{\mathbb{S}}

\newcommand{\Z}{\mathbb{Z}}

\newcommand{\Bbar}{{\mathbb{B}}}
\newcommand{\cyl}{\mathcal{K}}

\newcommand{\K}{\mathbb{K}}

\newcommand{\rr}{\mathrm{r}}
\newcommand{\sss}{\mathrm{s}}
\newcommand{\ttt}{\mathrm{t}}
\newcommand{\xx}{\mathrm{x}}
\newcommand{\yy}{\mathrm{y}}
\newcommand{\zz}{\mathrm{z}}

\newcommand{\geuc}{g_{_E}}

\newcommand{\auxm}{{g_{_A}}}
\newcommand{\auxcf}{\Omega}

\newcommand{\chiK}{\widehat{\chi}_{_{cat}}}
\newcommand{\chiP}[1]{\widehat{\chi}_{_{#1}}}

\newcommand{\cate}{\widehat{\kappa}}
\newcommand{\rhohat}{\widehat{\rho}}
\newcommand{\pihat}{\widehat{\pi}}

\newcommand{\Hcal}{\mathcal{H}}
\newcommand{\Qcal}{\mathcal{Q}}
\newcommand{\Fcal}{\mathcal{F}}
\newcommand{\Jcal}{\mathcal{J}}

\newcommand{\cbar}{\underline{c}}

\newcommand{\Lcal}{\mathcal{L}}
\newcommand{\Lchi}{\mathcal{L}_\chi}
\newcommand{\Lcat}{\widehat{\Lcal}_{cat}}

\newcommand{\Rcal}{\mathcal{R}}
\newcommand{\Rcat}{\widehat{\Rcal}_{cat}}
\newcommand{\Rplanar}[1]{\widehat{\Rcal}_{#1}}

\newcommand{\Rcalt}{\widetilde{\Rcal}}

\newcommand{\Pcal}{\mathcal{P}}

\newcommand{\mubar}{\overline{\mu}}

\newcommand{\wbar}{\overline{w}}

\newcommand{\what}{\widehat{w}}
\newcommand{\wbarhat}{\widehat{\wbar}}
\newcommand{\vbarhat}{\widehat{\overline{v}}}

\newcommand{\taubar}{\underline{\tau}}
\newcommand{\abar}{\underline{a}}

\newcommand{\sech}{\operatorname{sech}}

\newcommand{\arcosh}{\operatorname{arcosh}}
\newcommand{\supp}{\operatorname{supp}}

\newcommand{\Sbar}{\underline{S}}

\newcommand{\sgn}{\operatorname{sgn}}

\newcommand{\abs}[1]{\left\lvert#1\right\rvert}
\newcommand{\norm}[1]{\left\|#1\right\|}
\newcommand{\cutoff}[2]{\psi\left[ #1,#2 \right]}

\title{Free boundary minimal surfaces with connected boundary in the $3$-ball
 by tripling the equatorial disc}
\author{Nikolaos Kapouleas}
\author{David Wiygul}
\address{Department of Mathematics, Brown University, Providence, RI 02912}
\email{nicos@math.brown.edu}
\address{Department of Mathematics, University of California, Irvine, CA 92697}
\email{dwiygul@uci.edu}
\date{\today}

\begin{document}

\begin{abstract}
In the Euclidean unit three-ball, 
we construct compact, embedded, two-sided free boundary minimal surfaces
with connected boundary and prescribed high genus,
by a gluing construction tripling the equatorial disc.
Aside from the equatorial disc itself, 
these are the first examples in the three-ball of compact free boundary minimal surfaces with connected boundary.
\end{abstract}

\maketitle

\section{Introduction}
\label{S:intro}
$\phantom{ab}$
\nopagebreak
\subsection*{The general framework}
$\phantom{ab}$
\nopagebreak

Free boundary minimal surfaces in a Riemannian manifold with boundary are defined to be critical 
for the area functional
(under compactly supported perturbations)
subject to the constraint
that their boundary is contained in the boundary of the 
manifold. 
They are minimal surfaces which meet (along their boundary)
the boundary of the manifold orthogonally. 
The simplest examples are the equatorial disc $\Bbar^2$
in the Euclidean closed three-ball $\Bbar^3$  
and the critical catenoid $\K$ \cite{fraser-schoen:1},
which is the portion in $\Bbar^3$ of a suitably scaled catenoid in $\R^3$.  
Early work on free boundary minimal surfaces was by Courant \cite{courant} 
and subsequently by 
Nitsche \cite{nitsche76}, 
Taylor \cite{taylor77}, Hildebrandt-Nitsche \cite{hildebrandt-nitsche79}, 
Gr\"{u}ter-Jost \cite{gruter-jost86} and Jost \cite{jost86}. 

Further progress has been made more recently:
General existence results were obtained by Fraser \cite{fraser00} for disk type solutions, 
and later by Chen-Fraser-Pang \cite{chen-fraser-pang15} for incompressible surfaces. 
For embedded solutions in compact $3$-manifolds
a general existence result using min-max constructions was obtained by Martin Li \cite{Li15}. 
The min-max theory for free boundary minimal hypersurfaces in the Almgren-Pitts setting was recently developed by Li-Zhou \cite{Li-Zhou}. 
Fraser \cite{fraser07,fraser02} used index estimates to study the topology of Euclidean domains with $k$-convex boundary and 
Fraser-Li \cite{Fraser-Li14} proved a smooth compactness result for embedded free boundary minimal surfaces 
when the ambient manifold has nonnegative Ricci curvature and convex boundary. 
Maximo-Nunes-Smith \cite{maximo-nunes-smith17} used this last result 
and degree theory to prove the existence of free boundary minimal annuli in such three-manifolds.

In a recent breakthrough, Fraser-Schoen \cite{fraser-schoen:1} discovered a deep connection between free boundary minimal surfaces 
in the Euclidean unit ball $\Bbar^n\subset\R^n$ and extremal metrics on compact surfaces with boundary associated with the Steklov eigenvalue problem. 
In a follow-up article \cite{fraser-schoen:2}
Fraser-Schoen constructed new examples of embedded free boundary minimal surfaces 
with genus zero and arbitrary number of boundary components.  
This motivated a doubling construction of the equatorial disc
in the spirit of \cite{kapouleas-yang10} by
Folha-Pacard-Zolotareva \cite{FPZ},
where examples of genus zero or one and a large number of 
connected components are constructed 
(plausibly the genus zero ones being the same as the Fraser-Schoen examples).  
Recently more examples desingularizing an equatorial disc intersecting a critical catenoid 
were constructed by Ketover \cite{ket16} using min-max methods  
and N.K.-Li \cite{kapouleas:li} by gluing PDE methods. 

There are similarities between closed minimal surfaces in the round three-sphere and the
compact properly immersed free boundary minimal surfaces in the unit three-ball. 
For example in some sense the equatorial disc $\Bbar^2$ and the critical catenoid $\K$ in $\Bbar^3$ 
are analogous to the equatorial sphere $\Sph^2$ and the Clifford torus in $\Sph^3(1)$. 
$\Bbar^2$ is the unique (immersed) free boundary minimal disk in $\Bbar^3$ by a result of Nitsche \cite{nitsche85}, 
and by a surprising recent result of Fraser-Schoen \cite{fraser-schoen:3} 
the unique free boundary minimal disk in $\Bbar^n$ for any $n\ge3$.
Although (as shown by Almgren \cite{almgren})
$\Sph^2$ is the unique minimally immersed topological sphere in $\Sph^3$,
the analogous result is known (see \cite{calabi}) to fail in higher dimensions.
This contrast suggests that the free boundary surfaces in the Euclidean balls are even more rigid 
than the closed minimal surfaces in the round spheres. 

The only known examples of compact properly
embedded free boundary minimal surfaces in $\Bbar^3$ are 
those already mentioned,
that is $\Bbar^2$, $\K$, and those in \cites{fraser-schoen:2,FPZ,ket16,kapouleas:li}. 
A natural question is which topological types can be realized as such surfaces. 
In particular the existence of such surfaces of nonzero genus 
with a given small number of boundary components 
is a very natural question. 
The uniqueness results mentioned earlier suggest that in some cases such surfaces may not exist. 
Gluing constructions provide a satisfactory answer 
for high genus examples with a small number of boundary components, as described now.  

Examples with connected boundary can be constructed
by both doubling and desingularizing gluing constructions. 
In this article, 
we construct examples with connected boundary (besides the disc) for the first time.  
They have arbitrarily prescribed high genus 
and they are constructed by tripling the equatorial disc.
Note that in ongoing work, 
N.K. and Martin Li construct more examples with connected boundary 
by desingularizing two orthogonal discs using methodology from \cite{kapouleas:compact}.  
(Although the latter construction was observed earlier,
the former is more symmetric and easier to implement.) 
The latter construction also provides examples with two boundary components. 
Examples with three boundary components are provided by desingularization methods in \cite{kapouleas:li} 
(and also by min-max methods in \cite{ket16}).   
Examples with four boundary components can be constructed by doubling the catenoid as in \cite{LDa}, which is in preparation. 

The general idea of doubling constructions
by gluing methods was proposed and discussed in 
\cites{kapouleas:survey,kapouleas:clifford,kapouleas:rs}.
Gluing methods have been applied extensively
and with great success in gauge theory by Donaldson, Taubes, and others.
The particular
gluing methods used in this article relate most
closely to the methods developed by Richard Schoen in \cite{schoen}
and N.K. in \cite{kapouleas:annals},
especially as they evolved and were systematized in
\cites{kapouleas:wente:announce,kapouleas:wente,kapouleas:imc}.
We refer to \cite{kapouleas:survey} for a general discussion of this gluing methodology 
and to \cite{kapouleas:rs} for a detailed general discussion of doubling by gluing methods.

Roughly speaking, in such doubling constructions,  
given a minimal surface $\Sigma$,
one constructs first a family of smooth, embedded, and approximately minimal initial surfaces. 
Each initial surface 
consists of two approximately parallel copies of $\Sigma$
with a number of discs removed and replaced by approximately catenoidal bridges.
Finally one of the initial surfaces in the family
is perturbed to minimality by partial differential equations methods.
Understanding such constructions in full generality
seems beyond the immediate horizon at the moment.
In the earliest such construction \cite{kapouleas:clifford},
where doublings of the Clifford torus are constructed, 
there is so much symmetry imposed 
that the position of the catenoidal bridges is completely fixed 
and all bridges are identical modulo the symmetries.
D.W. \cite{wiygul:thesis,wiygul:stacking} 
has extended that construction to ``multiple doublings'', or ``stackings'',  
with more than two copies of the Clifford torus involved
(and some less symmetric doublings also), 
where the symmetries do not determine the vertical
(that is perpendicular to $\Sigma$) position of the bridges.

The construction in this article is most closely related to this last construction. 
To keep this article as simple as possible, 
we limit our attention to ``triplings'', 
where we connect three copies of the equatorial disc with ``half-catenoidal'' bridges
(truncated by the boundary sphere of the ball approximately in half).  
This approach can be extended to apply to an arbitrary number of copies
as in \cite{wiygul:thesis,wiygul:stacking};
see Remark \ref{stacking} below.
The tripling of the equatorial disc has the surprising property
that the boundary of the surfaces obtained is connected. 
All our catenoidal bridges are equivalent modulo the symmetries
and their horizontal but not vertical
positions are
fixed by the symmetries 
and the boundary condition. 
The bridges are placed along the boundary circle of the equatorial disc
and alternate in connecting the top copy of the disc with the equatorial disc (middle copy) 
and the bottom copy with the middle copy. 
Note that ``Linearized Doubling'',
a methodology developed to deal with horizontal forces and other difficulties 
\cite{kapouleas:equator,kapouleas:ii,LDa}, 
is not required in the current construction. 

\subsection*{Brief discussion of the results}
$\phantom{ab}$
\nopagebreak

Our constructions depend on a large integer $m$,
which determines also the group $\Grp[m]$ of symmetries of the construction;
$\Grp[m]$ is an index-two subgroup of
the group of isometries of $\R^3$ 
preserving the union of $m$ lines on the equatorial plane through the center of the ball 
and arranged symmetrically. 
More precisely, 
$\Grp[m]$ is the subgroup preserving alternating sides of the wedges into which these $m$ lines subdivide 
the equatorial plane.  
We consider three parallel copies of the equatorial disc, one being the actual equatorial disc 
and the other two at equal distances above and below. 
The $m$ lines divide the equatorial circle (boundary of equatorial disc) into 
$2m$ equal arcs, the middle points of which form a collection $L_0$. 
We connect the three copies of the equatorial disc with catenoidal bridges (of appropriate size) 
which have vertical axes through the points of $L_0$ and alternate connecting the middle disc 
with the top or the bottom. 

This provides connected compact surfaces with connected boundary on the boundary of
$\Bbar^3$,
which we call ``pre-initial surfaces''. 
These surfaces are then perturbed using a bare-hands approach to the ``initial surfaces'' 
which satisfy the required boundary condition (orthogonal intersection with the boundary of the ball) 
but are only approximately minimal. 
Finally applying the Schauder fixed-point theorem, 
we prove that one of the initial surfaces can be perturbed 
to provide the desired surface (see
Theorem \ref{mainthm}): 

\begin{theorem}[Main Theorem]
\label{mainthmIntro}
If $m$ is large enough, 
one of the initial surfaces outlined above can be perturbed to a
compact, embedded, two-sided
free boundary minimal surface in the unit three-ball 
which is invariant under $\Grp[m]$ and has connected boundary and genus $m-1$. 
The minimal surfaces obtained tend as $m\to\infty$ to the equatorial disc with multiplicity three
and the length of their boundary tends to $6\pi$. 
\end{theorem}

\begin{remark}
  \label{stacking}
More generally, 
for any positive integer $N$ and any sufficiently large integer $m$, 
we can produce by gluing methods embedded free-boundary minimal surfaces
resembling $N$ parallel copies of the equatorial disc,
with each adjacent pair of copies joined by $m$ catenoidal strips,
in maximally symmetric fashion.
For odd $N$ the resulting surfaces have connected boundary
and genus $(m-1)\frac{N-1}{2}$,
while for even $N$
the resulting surfaces have $m$ boundary components
and genus $(m-1)\frac{N-2}{2}$.
To simplify the presentation in this article we carry out only the case $N=3$.
\end{remark}

\begin{remark}
  \label{onecomponent}
Tripling constructions
(or more generally stacking constructions with an odd number of copies) 
of any free boundary minimal surface
with half catenoidal bridges at the boundary in the manner of the current construction
produce examples with the same number of boundary components as the original surface.
\end{remark}

\subsection*{Outline of the approach}
$\phantom{ab}$
\nopagebreak

The families of the pre-initial surfaces we construct are parametrized by two continuous parameters: 
one which controls the size of the catenoidal bridges and one which controls the distance of the centers 
of the catenoidal bridges from the equatorial plane. 
The choice of these parameters is motivated by the Geometric Principle
(see \cite{kapouleas:survey,kapouleas:rs} for example)
because they control the ``dislocation'' of the initial surfaces,
that is the (vertical) repositioning of the 
copies of the disc and the catenoidal bridges. 
Motivated by this, 
we solve the linearized equation modulo a two-dimensional extended substitute kernel spanned 
by one function which allows us to solve orthogonally to the constants on the top
(equivalently by symmetry the bottom)
disc, 
and one which allows us to ensure decay away from the middle disc.
(The decay from the top disc is not obstructed.)

We use a conformal metric in the ball to describe the graphical perturbation of the initial surfaces. 
This has various advantages and simplifies the treatment of the boundary equation,  
provided that the initial surfaces intersect the boundary orthogonally. 
This is easy to arrange by appropriately modifying the pre-initial surfaces: 
the initial surfaces have the same boundary
as the pre-initial surfaces but the intersection angle with the boundary of the ball has been corrected to orthogonality.

\subsection*{Organization of the presentation}
$\phantom{ab}$
\nopagebreak

Besides the introduction, this article has six more sections. 
In Section \ref{nc}, we fix some useful notation. 
In Section \ref{S:pi}, we construct and provide precise estimates for the pre-initial surfaces. 
In Section \ref{S:i}, we modify the pre-initial surfaces to construct the initial surfaces which we also estimate carefully. 
In Section \ref{graphs}, we define the ambient conformal metric and study the new graphical perturbations of the initial surfaces. 
In Section \ref{S:lp}, we study the linearized equation on the initial surfaces. 
In Section \ref{end}, we state and prove the main theorem.

\subsection*{Acknowledgments}
$\phantom{abbbbbbbbbbb}$ 
\nopagebreak

The authors would like to thank Richard Schoen for his continuous support and interest in the results of this article. 
N.K. would like to thank the Mathematics Department at the University of California, Irvine 
and the Mathematics Department and the MRC at Stanford University, 
for providing a stimulating mathematical environment and generous partial financial support during the academic year 2015-2016. 
N.K. was also partially supported by NSF grant DMS-1405537.

\section{Notation and conventions}
\label{nc}
We will employ standard Cartesian coordinates $(\xx,\yy,\zz)$ on $\R^3$
as well as the radial spherical coordinate $\rr$, measuring distance from the origin,
and the angular cylindrical coordinate $\theta$, relative to the positive $\xx$-axis
as usual.
We will routinely identify $\R^2$ with $\{\zz=0\} \subset \R^3$
and restrict to it the above coordinate functions without relabelling.
We set $\Bbar^3:=\{\rr \leq 1\} \subset \R^3$,
the closed unit $3$-ball in $\R^3$ centered the origin,
and we set $\Bbar^2:=\Bbar^3 \cap \{\zz=0\}$,
the closed unit disc (or $2$-ball)
centered at the origin in $\R^2$.

We fix a smooth, nondecreasing $\Psi: \R \to [0,1]$ with $\Psi$  identically $0$
on $(-\infty,-1]$, identically $1$ on $[1,\infty)$,
and such that $\Psi - \frac{1}{2}$ is odd. We then define, for any $a,b \in \R$, the function $\cutoff{a}{b}: \R \to [0,1]$ by
  \begin{equation}
    \cutoff{a}{b} = \Psi \circ L_{a,b},
  \end{equation}
where $L_{a,b}: \R \to \R$ is the linear function satisfying $L(a)=-3$ and $L(b)=3$.

Given an open set $\Omega$ of a submanifold (possibly with boundary)
immersed in an ambient manifold (possibly with boundary) endowed with metric $g$,
an exponent $\alpha \in [0,1)$, 
and a tensor field $T$ on $\Omega$,
possibly taking values in the normal bundle,
we define the pointwise H\"{o}lder seminorm
\begin{equation}
  [T]_\alpha(x) := \sup_{y \in B_x} 
  \frac{\abs{T(x)-P_y^xT(y)}_g}{d(x,y)^\alpha},
\end{equation}
where $B_x$ denotes the open geodesic ball (or possibly half ball),
with respect to $g$, 
with center $x \in \Omega$ and radius the minimum of $1$ and the injectivity radius at $x$;
$\abs{\cdot}_g$ denotes the pointwise norm induced by $g$; 
$P_y^x$ denotes parallel transport, also induced by $g$, 
from $y$ to $x$ along the unique geodesic in $B_x$ 
joining $y$ and $x$;
and $d(x,y)$ denotes the distance between $x$ and $y$.

Given further a positive function
$f: \Omega \to \R$ and a nonnegative integer $k$,
assuming that all order-$k$ partial derivatives of the section $T$ 
(with respect to any coordinate system) exist and are continuous, we set
\begin{equation}
  \norm{T: C^{k,\alpha}(\Omega,g,f)}
  := \sum_{j=0}^k \sup_{x \in \Omega} \frac{\abs{D^jT(x)}}{f(x)}
  + \frac{\left[D^k T\right]_\alpha(x)}{f(x)},
\end{equation}
where $D$ denotes the Levi-Civita connection determined by $g$.
In this article the vector space $\R$ is always endowed with its standard norm $\abs{\cdot}$,
the absolute value,
and every product of normed vector spaces is always endowed with the norm
defined as the sum of the norms of its factors.

Now let $\Sigma$ be a two-sided hypersurface in a Riemannian manifold $(M,g)$,
with global unit normal $\nu$,
and let $\Grp$ be a subgroup of the group of isometries of $M$
that preserve the set $\Sigma$.
We say that a function $u: \Sigma \to \R$
is \emph{$\Grp$-odd} if for every $\mathfrak{g}: M \to M$ belonging to $\Grp$
  \begin{equation}
    u \circ \mathfrak{g}|_{_\Sigma}
    =
    g|_{_\Sigma}(\nu \circ \mathfrak{g}|_{_\Sigma}, \, \mathfrak{g}_*\nu)u
    =\pm u,
  \end{equation}
the last sign depending on whether
$\mathfrak{g}$ preserves or reverses the sides of $\Sigma$;
similarly we say that $u: \Sigma \to \R$ is \emph{$\Grp$-even}
(or simply $\Grp$-invariant) if for every $\mathfrak{g}: M \to M$ belonging to $\Grp$ instead
  \begin{equation}
    u \circ \mathfrak{g}|_{_\Sigma} = u.
  \end{equation}

Last, for any subset $S$ of a Riemannian manifold $(M,g)$
we write $d[S,M,g]: M \to \R$ for the function
whose value at a point in $M$ is that point's distance to $S$;
when context permits we will frequently abbreviate
$d[S,M,g]$ to $d[S,g]$ or even $d[S]$.
When $S$ is a submanifold of $(M,g)$
(not necessarily a hypersurface),
we will sometimes write $g|_{_S}$
for the Riemannian metric on $S$ induced by $g$.

\section{The pre-initial surfaces}
\label{S:pi} 

\subsection*{Construction}
$\phantom{ab}$
\nopagebreak

Given an integer $m>0$ let $\Grp[m]$
be the subgroup of $O(3)$ generated by reflection through
the plane $\{\yy=0\}$
and by reflection through the line
$\left\{\yy=\xx \tan \frac{\pi}{2m}\right\} \cap \{\zz=0\}$.
In particular $\Grp[m]$ includes
(i) the $m$ rotations about the $\zz$-axis through angles of the form
$j\frac{2\pi}{m}$ for each $j \in \Z$,
(ii) the $m$ reflections through vertical planes of the form 
$\left\{\yy=\xx \tan j\frac{\pi}{m}\right\}$
for each $j \in \Z$,
and (iii) the $m$ reflections through horizontal lines of the form
$\left\{\yy=\tan j\frac{\pi}{2m}\right\} \cap \{\zz=0\}$ for each $j \in 2\Z+1$.
It will often be useful to isolate the subgroup $\Grp^+[m]<\Grp[m]$
generated by just reflection through the plane $\{\yy=0\}$
and rotation about the $\zz$-axis through angle $\frac{2\pi}{m}$.
Thus $\Grp^+[m]$ is isomorphic to the dihedral group of order $2m$.
For convenience of notation we will sometimes
prefer to denote $\Grp[m]$ by $\Grp_0[m]$ and $\Grp^+[m]$
by $\Grp_1[m]$, so that
  \begin{equation}
  \label{Grp+def}
    \Grp_1[m]=\Grp^+[m]<\Grp[m]=\Grp_0[m].
  \end{equation}

We also define the sets
  \begin{equation}
  \begin{aligned}
  \label{Wmdef}
    &W[m]:=\left\{(\xx,\yy) \in \Bbar^2 \; | \; 
      \xx \geq 0 \mbox{ and } \abs{y} \leq \xx \tan \frac{\pi}{2m} \right\},
    \qquad
    \widehat{W}[m]:=W[m] \times \R \subset \R^3, \\
    &L_i[m]:=\Grp_i[m] \{(1,0,0)\} \subset \partial \Bbar^2,
    \quad \mbox{ and } \quad
    \widehat{L}_i[m]:=L_i[m] \times \R \qquad \mbox{for each $i \in \{0,1\}$},
  \end{aligned}
  \end{equation}
so that $L_0[m]$ consists of the $2m$\textsuperscript{th} roots of unity on the equator,
$L_1[m]$ consists of the $m$\textsuperscript{th} roots of unity on the equator,
$\Grp[m] W[m] = \Bbar^2$,
$\Grp[m]\widehat{W}[m]=\Bbar^2 \times \R$, and
$\Bbar^2$ ($\Bbar^2 \times \R$)
is the union of $2m$ congruent copies of $W[m]$ ($\widehat{W}[m]$ respectively).
The surfaces we build will depend on the positive integer $m$
and two additional parameters
$\zeta, \xi \in \R$.

\begin{remark}
\label{suppress}
Throughout the construction we will define many objects depending on
the three pieces of data $m \in \Z^+$ and $\zeta,\xi \in \R$,
but we will routinely suppress this dependence from our notation
when there is little danger of confusion.
\end{remark}
Given $m \in \Z^+$ and $\zeta,\xi \in \R$ we define the constants
$\tau$, $a$, and $\zz_0$,
as well as the functions $\phi_{cat}: [\tau,\infty) \to \R$
and $\phi_0, \phi_1: \{\rr \geq \tau\} \subset \R^2 \to \R$
by
  \begin{equation}
  \label{phi01}
    \begin{aligned}
    &\tau := \tau[m,\zeta] := \frac{1}{20m}e^{\zeta-\frac{1}{2}m}, \quad
    \phi_{cat}(\rr):=\phi_{cat}[m,\zeta](\rr):=\tau \arcosh \frac{\rr}{\tau}, \quad
    \zz_0 := \zz_0[m,\zeta,\xi] := \tau \xi + \tau a, \\
    &a := a[m,\zeta] := \tau^{-1}\phi_{cat}\left(\frac{1}{20m}\right)
      = \frac{1}{2}m-\zeta+\ln \left(1+\sqrt{1-e^{2\zeta-m}}\right)
      = \frac{1}{2}m-\zeta+\ln 2 + O\left(m^2\tau^2\right), \\
    &\phi_0(\xx,\yy):=\phi_0[m,\zeta,\xi](\xx,\yy)
    :=
    \left[\zz_0 - \phi_{cat}(\rr)\right] \cdot \cutoff{\frac{1}{10m}}{\frac{1}{20m}}(\rr),
      \mbox{ and} \\
    &\phi_1(\xx,\yy):=\phi_1[m,\zeta,\xi](\xx,\yy)
    :=
    \zz_0 + \phi_{cat}(\rr) \cdot \cutoff{\frac{1}{10m}}{\frac{1}{20m}}(\rr) 
      + \tau a \cdot \cutoff{\frac{1}{20m}}{\frac{1}{10m}}(\rr).
    \end{aligned}
  \end{equation}
In the estimate for $a$ the term $O\left(m^2\tau^2\right)$
satisfies the inequality $\abs{O\left(m^2\tau^2\right)} \leq Cm^2\tau^2$
for some constant $C>0$ independent of $\zeta$, $\xi$, and $m$.
Motivation for the choice of $\tau$ can be found in Section \ref{end}.

In a neighborhood of the $\zz$-axis the graphs of $\phi_0$ and $\phi_1$
coincide with the lower ($\phi_0$) and upper ($\phi_1$) 
halves of the catenoid of waist radius $\tau$
centered at $(0,0,\zz_0)$,
but at a distance of order $\frac{1}{m}$ from the $\zz$-axis
the lower half levels off to agree with the $\{\zz=0\}$ plane
outside a cylinder of radius $\frac{1}{10m}$ about the $\zz$-axis,
while the upper half instead levels off to agree, outside the same cylinder,
with the horizontal plane
which it intersects at a distance $\frac{1}{20m}$ from the $\zz$-axis.
We will use translated truncations of these graphs,
along with an additional, exactly planar region,
to glue together our pre-initial surface.

Specifically we introduce
  \begin{equation}
  \label{Gammadef}
    \begin{aligned}
      &\digamma_0:=\digamma_0[m,\zeta,\xi]:=\Bbar^3 \cap
        \left\{ (\xx,\yy,\phi_0(\xx-1,\yy))  \; | \; 
        (\xx,\yy) \in W[m] \mbox{ and } (\xx-1)^2+\yy^2 \geq \tau^2 \right\}, \\
      &\digamma_1:=\digamma_1[m,\zeta,\xi]:=\Bbar^3 \cap
        \left\{ (\xx,\yy,\phi_1(\xx-1,\yy)) \; | \;
        (\xx,\yy) \in W[m] \mbox{ and } (\xx-1)^2+\yy^2 \geq \tau^2 \right\}, \mbox{ and} \\
      &\digamma_{-1}:=\digamma_{-1}[m,\zeta,\xi]:=\Bbar^3 \cap
        \left\{ (\xx,\yy, -\zz_0-\tau a) \; (\xx,\yy) \in W[m] \right\}
    \end{aligned}
  \end{equation}
and then define the pre-initial surface
  \begin{equation}
  \label{preinitsdef}
    \pre{\Sigma}
    :=
    \pre{\Sigma}[m,\zeta,\xi]
    :=
    \Grp[m] \left(\digamma_{-1} \cup \digamma_0 \cup \digamma_1\right),
  \end{equation}
the union of $2m$ congruent copies, having pairwise disjoint interiors,
of $\pre{\Sigma} \cap \widehat{W}[m]=\digamma_{-1} \cup \digamma_0 \cup \digamma_1$.
Following Remark \ref{suppress} we will abbreviate $\pre{\Sigma}[m,\zeta,\xi]$
by $\pre{\Sigma}$ whenever context permits.
Thus, outside a small tubular neighborhood of $\widehat{L}_0[m]$ (with radius of order $1/m$)
the pre-initial surface $\pre{\Sigma}$ is a union of three horizontal discs
close to the equatorial one;
on the other hand
the intersection of $\pre{\Sigma}$ with a small tubular neighborhood
(with radius of order $1/m$) of a single line in $\widehat{L}_1$
is a truncated catenoid cut nearly in half along its axis by the sphere $\partial \Bbar^3$.

  \begin{figure}[H]
    \includegraphics[width=2in]{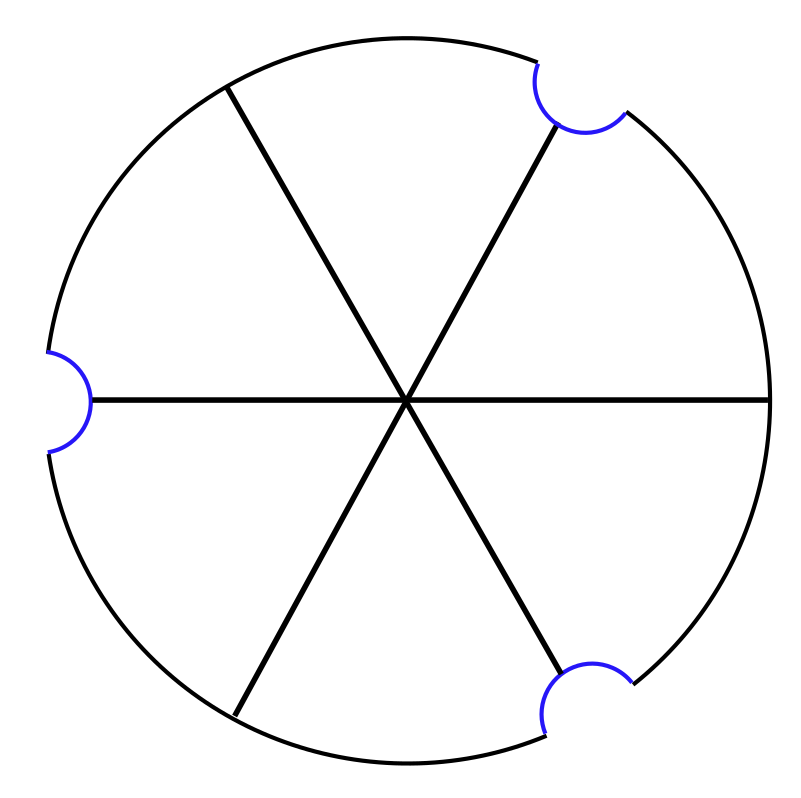}
    \includegraphics[width=2in]{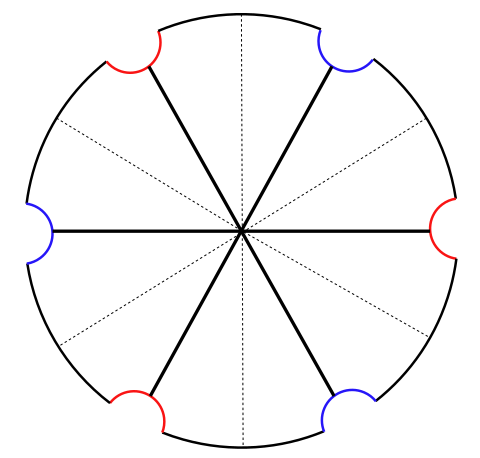}
    \includegraphics[width=2in]{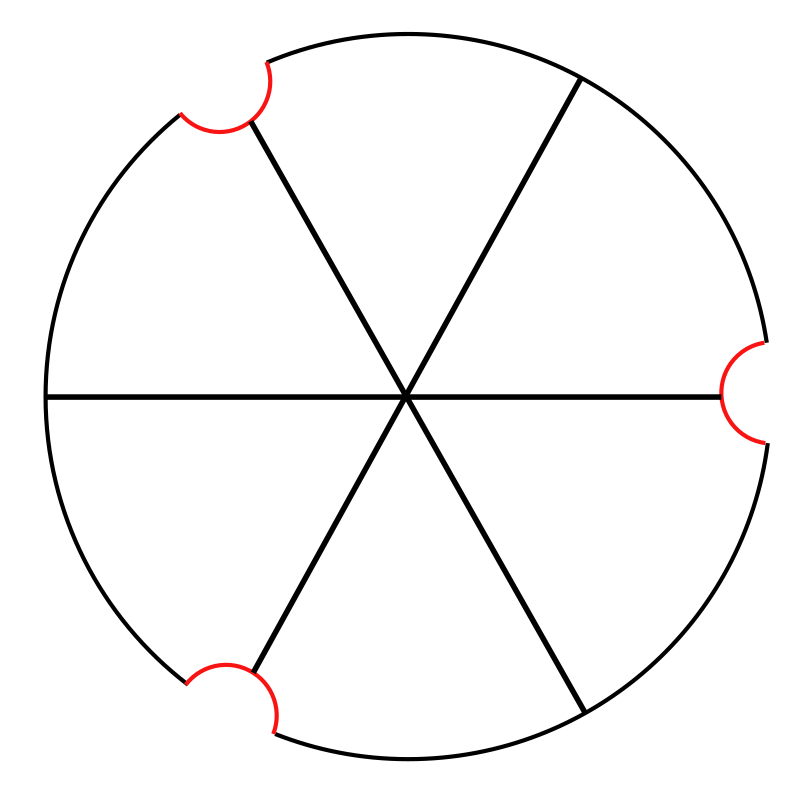}
    \caption{A schematic depiction of (from left to right)
      the bottom, middle, and top planar regions
      (defined in (\ref{planardef}))
      of a pre-initial surface with $m=3$.
      The three indented discs are connected through catenoidal strips
      attaching to the blue and red arcs in alternating fashion.
      The solid segments indicate vertical planes of symmetry,
      while the dotted segments indicate horizontal lines of symmetry.}
  \end{figure}

\begin{prop}[Topology of the pre-initial surfaces]
\label{preinittop}
Given $c>0$, there exists $m_0>0$ such that
for every integer $m>m_0$ and for every $\zeta,\xi \in [-c,c]$
the pre-initial surface $\pre{\Sigma}[m,\zeta,\xi]$
is a compact, orientable, smooth surface with boundary;
$\pre{\Sigma}$ is embedded in $\Bbar^3$, has genus $m-1$, and is preserved as a set by $\Grp[m]$;
and its boundary $\partial \pre{\Sigma}$ has a single connected component,
which is embedded in $\partial \Bbar^3$.
\end{prop}

\begin{proof}
It is clear from the definition that the set
$C:=\partial \pre{\Sigma} \cap \{0 \leq \theta \leq 2\pi/m\}$ is connected.
The connectedness of $\partial \pre{\Sigma}$
then follows from the observations that
$C$ intersects its image under rotation about the $\zz$-axis
through angle $2\pi/m$
and that $\partial \pre{\Sigma}$ is the union over $j \in \Z$
of the image of $C$ under rotation about the $\zz$-axis by $2j\pi/m$.
By topologically gluing $\pre{\Sigma}$ to a homeomorphic copy of itself
along $\partial \pre{\Sigma}$ (that is doubling $\pre{\Sigma}$ in the standard topological sense)
one obtains a closed surface of genus $2m-2$, since each of the three discs
in $\pre{\Sigma}$ extends to a sphere, $2$ of the $2m$ catenoids connect these three spheres,
and each remaining catenoid contributes genus to the surface.
On the other hand the genus of the topological double
of a surface with $k$ boundary components and genus $g$
is $2g+(k-1)$,
so $\pre{\Sigma}$ has genus $m-1$.

\end{proof}

\subsection*{Catenoidal and planar regions}
$\phantom{ab}$
\nopagebreak

By construction every pre-initial surface decomposes into overlapping regions,
each of which resembles either a catenoid or a plane,
truncated and subjected to small perturbations.
Modulo the symmetries we have one catenoidal region and two planar regions.

Given $s>0$ we write
  \begin{equation}
    \cyl_s := [-s,s] \times \Sph^1
  \end{equation}
for the standard cylinder of radius $1$ and length $2s$,
and, recalling the definition of $a$ in (\ref{phi01}),
we define the embedding $\cate: \cyl_a \to \R^3$ by
  \begin{equation}
  \label{cate}
    \cate(t,\vartheta)
    :=
    (1,0,\zz_0)+\tau(-\cosh t \cos \vartheta, \cosh t \sin \vartheta, t)
    =
    (1-\tau \cosh t \cos \vartheta, \tau \cosh t \sin \vartheta, \tau t + \zz_0),
  \end{equation}
where $(t,\vartheta)$ are the standard global coordinates
for the universal cover of $\cyl_a$.
We also define
  \begin{equation}
  \label{precatr}
    \pre{\Sigma}_{cat}
    :=
    \pre{\Sigma} \cap \cate \left( \cyl_a \right)
    \subset
    \cate\left(\cyl_a^+\right),
  \end{equation}
where $\cyl_a^+$ is the half cylinder
  \begin{equation}
    \cyl_a^+:=\left\{\abs{\vartheta} \leq \frac{\pi}{2} \right\} \cap \cyl_a;
  \end{equation}
for large $m$ we can regard $\pre{\Sigma}_{cat}$ as a small perturbation
of the corresponding half catenoid $\cate\left(\cyl_a^+\right)$.

We write
  \begin{equation}
    \partial_{\Sph^2}\pre{\Sigma}_{cat}
    :=
    \partial \pre{\Sigma}_{cat} \cap \partial \Bbar^3
    =
    \pre{\Sigma}_{cat} \cap \partial \pre{\Sigma}
  \end{equation}
for the spherical part of the boundary of $\pre{\Sigma}_{cat}$.
It is easy to check from \ref{cate} that
  \begin{equation}
  \label{catbdycurve}
    \partial_{\Sph^2}\pre{\Sigma}_{cat}
    =
    \cate
      \left(
        \left\{ \left. \left(t,\frac{\pi}{2}\lambda(t)\right) \; \right| \; 
          t \in [-a,a] \right\}
      \right)
    \cup
      \cate\left(
        \left\{ \left. \left(t,-\frac{\pi}{2}\lambda(t)\right) \; \right| \; 
          t \in [-a,a] \right\}
      \right),
  \end{equation}
where $\lambda: [-a,a] \to \R$ is a height-dependent scaling factor
(for the $\vartheta$ interval)
defined by
  \begin{equation}
  \label{lambdadef}
    \lambda(t):
    =1-\frac{2}{\pi}\arcsin \frac{\tau^2 \cosh^2 t + (\tau t + \zz_0)^2}{2\tau \cosh t}.
  \end{equation}
These curves are small perturbations of two catenaries in the plane $\xx=1$,
the images under $\cate$ of the vertical lines $\vartheta=\pm \frac{\pi}{2}$ on $\cyl_a$.
To identify $\pre{\Sigma}_{cat}$ with the exact half cylinder $\cyl_a^+$
we introduce
the map $\pre{\kappa}: \cyl_a^+ \to \pre{\Sigma}_{cat}$ given by
  \begin{equation}
  \label{prekappadef}
    \pre{\kappa}(t,\vartheta)=\cate(t, \lambda(t)\vartheta),
  \end{equation}
which for sufficiently large $m$ is clearly a diffeomorphism.
Moreover $\pre{\kappa}$ restricts to a diffeomorphism
of $\partial \cyl_a^+ \backslash \{\abs{t}=a\}$ onto
$\partial_{\Sph^2}\pre{\Sigma}_{cat}$.

The planar regions are small graphs over the equatorial disc $\Bbar^2$
indented near the catenoids.
Referring to (\ref{phi01})
we see that by taking $m$ large enough in terms of $\zeta$
we can ensure that $m/4<a$.
With this in mind we define the top, middle, and bottom planar regions---$\pre{\Sigma}_1$,
$\pre{\Sigma}_0$, and $\pre{\Sigma}_{-1}$ respectively---as the connected
components of $\pre{\Sigma} \backslash \Grp[m] \left(\cate\left(\cyl_{m/4}\right)\right)$,
so that
  \begin{equation}
  \label{planardef}
    \begin{aligned}
      &\pre{\Sigma}=\pre{\Sigma}_1 \cup \pre{\Sigma}_0 \cup \pre{\Sigma}_{-1}
        \cup \Grp[m] \left(\cate\left(\cyl_{m/4}\right)\right) \mbox{ and} \\
      &\pre{\Sigma}_i \supset \digamma_i \mbox{ for each } i \in \{-1,0,1\},
    \end{aligned}
  \end{equation}
recalling (\ref{Gammadef}).
(Each boundary circle of $\cate\left(\cyl_{m/4}\right)$ has radius
$\tau \cosh \frac{m}{4}
 =\frac{e^{\zeta/2}}{2\sqrt{20m}}\sqrt{\tau}
   +\frac{1}{2}e^{-\zeta/2}\sqrt{20m}\tau^{3/2}$.)
For each $i \in \{-1,0,1\}$
we write
  \begin{equation}
    \partial_{\Sph^2}\pre{\Sigma}_i
    :=
    \partial \pre{\Sigma}_i \cap \partial \Bbar^3
    =
    \pre{\Sigma}_i \cap \partial \pre{\Sigma}
  \end{equation}
for the spherical part of the boundary of $\pre{\Sigma}_i$.
See Figure 1 for an overhead view of the various planar regions.

For $m$ large enough $\abs{\zz}$ is small on $\pre{\Sigma}$
and $\pre{\pi}_i: \pre{\Sigma}_i \to \Bbar^2$ ($i \in \{0,1\}$)
defined by
  \begin{equation}
  \label{prepidef}
    \pre{\pi}_i(\xx,\yy,\zz)=\frac{(\xx,\yy)}{\sqrt{1-\zz^2}}
  \end{equation}
is a diffeomorphism onto a relatively open subset of $\Bbar^2$,
and its restriction to $\partial_{\Sph^2} \pre{\Sigma}_i$
is a diffeomorphism onto a subset of the equator $\partial \Bbar^2$.
Note from (\ref{phi01}) that the factor $\left(1-\zz^2\right)^{1/2}$
is almost constantly $1$
on each $\pre{\Sigma}_i$
and is needed to ensure that $\pre{\pi}_i$ maps
$\partial_{\Sph^2} \pre{\Sigma}_i$
to the boundary of the equatorial disc $\Bbar^2$.
Clearly $\pre{\Sigma}=\pre{\Sigma}_0 \cup \Grp[m] \pre{\Sigma}_1 \cup \Grp[m] \pre{\Sigma}_{cat}$.

\subsection*{Geometric estimates}
$\phantom{ab}$
\nopagebreak

We write $\geuc$ for the Euclidean metric on $\R^3$,
$\pre{X}: \pre{\Sigma}\to \R^3$ for the inclusion map for $\pre{\Sigma}$,
$\pre{g}=\pre{X}^*\geuc$ for the induced metric,
$\pre{\nu}: \pre{\Sigma} \to \R^3$ for the unit normal to
$\pre{\Sigma}$ which points upward on the top planar region
(so downward at the origin $(0,0,0) \in \pre{\Sigma}$),
$\pre{A}_{ab}:=-\pre{\nu}_{i|b}\pre{X}^i_{\;,a}$
for the corresponding second fundamental form
(in which expression the vertical bar indicates covariant differentiation as
induced by $\geuc$ on the bundle $\pre{X}^*T^*\R^3$ over $\pre{\Sigma}$),
$\pre{H}:=\pre{g}^{ab}\pre{A}_{ab}$ for the corresponding mean curvature,
and $\pre{\Theta}: \partial \pre{\Sigma} \to \R$ for the Euclidean inner product
  \begin{equation}
  \label{preTheta}
    \pre{\Theta}
    :=\langle \pre{X}, \pre{\nu} \rangle |_{\partial \pre{\Sigma}}
  \end{equation}
of $\pre{\nu}$
with $\pre{X}$---here regarded as the position vector field on $\pre{\Sigma}$---along
$\partial \pre{\Sigma}$.
Note that $\pre{\Theta}$ is $\Grp[m]$-odd, in the sense defined in Section \ref{nc}.

For future applications
we introduce coordinates $(\pre{s},\pre{\sigma})$
on a neighborhood in $\pre{\Sigma}$ of $\partial \pre{\Sigma}$
small enough so that within it
the map of nearest-point projection onto $\partial \pre{\Sigma}$
is well-defined and smooth;
for any point $p$ in this neighborhood
$\pre{\sigma}(p)$ is the distance in $\pre{\Sigma}$ of $p$ from $\partial \pre{\Sigma}$
and $\pre{s}(p)$ is the distance in $\partial \pre{\Sigma}$
of the nearest-point projection of $p$ onto $\partial \pre{\Sigma}$
from an arbitrarily fixed reference point on $\partial \pre{\Sigma}$.
In particular,
along $\partial \pre{\Sigma}$
the coordinate vector field $\partial_{\pre{\sigma}}$ is the inward unit conormal
for $\pre{\Sigma}$,
and $\{\pre{X}_*\partial_{\pre{s}},\pre{X}_*\partial_{\pre{\sigma}},\nu\}$
is an orthonormal basis for $\R^3$ at each point of $\partial \pre{\Sigma}$.

Obviously $\pre{\Theta}$ defined above equivalently encodes the angle
between the conormal to $\pre{\Sigma}$ and the normal to $\partial \Bbar^3$.
Note that $\partial_{\pre{s}}$ is tangential to $\partial \Bbar^3$
along $\partial \pre{\Sigma}$, so orthogonal to $\pre{X}$, and therefore
  \begin{equation}
    \pre{X}|_{\partial \pre{\Sigma}}
    =
    \langle \pre{X}, \pre{X}_*\partial_{\pre{\sigma}} \rangle
      \pre{X}_*\partial_{\pre{\sigma}}|_{\partial \pre{\Sigma}}
      + \langle \pre{X}, \pre{\nu} \rangle \pre{\nu}|_{\partial \pre{\Sigma}}.
  \end{equation}
Since the left-hand side has unit length,
it follows, using the Pythagorean theorem and (\ref{preTheta}), that
  \begin{equation}
  \label{preconormal}
    \langle \pre{X}, \pre{X}_*\partial_{\pre{\sigma}} \rangle|_{\partial \pre{\Sigma}}
    =
    -\sqrt{1-\pre{\Theta}^2},
  \end{equation}
keeping in mind that $\partial_{\pre{\sigma}}$ points toward the interior of $\Bbar^3$
while $\pre{X}$ points toward the exterior.
   
It will be useful to equip $\pre{\Sigma}$ not only with its natural metric $\pre{g}$
but also with the conformal metric
  \begin{equation}
    \pre{\chi} = \pre{\rho}^2 \pre{g},
  \end{equation}
with conformal factor
  \begin{equation}
    \pre{\rho}: \pre{\Sigma} \to \R
  \end{equation}
defined to be the unique $\Grp[m]$-even function
having the restrictions 
  \begin{equation}
  \label{prerho}
    \begin{aligned}
      &\pre{\rho}\left|_{\pre{\Sigma}_{cat}} \right.
      :=
      m \cutoff{\frac{1}{10m}}{\frac{1}{5m}}
        \circ \left.d\left[\widehat{L}_0,\R^3,\geuc\right]\right|_{ \pre{\Sigma}_{cat} }
        + \frac{1}{d\left[\widehat{L}_0,\R^3,\geuc\right]\left|_{ \pre{\Sigma}_{cat} } \right.}
        \cutoff{\frac{1}{5m}}{\frac{1}{10m}}
        \circ d\left[\widehat{L}_0,\R^3,\geuc\right]\left|_{ \pre{\Sigma}_{cat} } \right. 
        \mbox{ and} \\
      &\pre{\rho}\left|_{ \pre{\Sigma}_i } \right.
      := 
      m \cutoff{\frac{1}{10m}}{\frac{1}{5m}}
        \circ \left.d\left[\widehat{L}_i,\R^3,\geuc\right]\right|_{ \pre{\Sigma}_i }
        + \frac{1}{d\left[\widehat{L}_i,\R^3,\geuc\right]\left|_{ \pre{\Sigma}_i } \right.}
        \cutoff{\frac{1}{5m}}{\frac{1}{10m}}
        \circ d\left[\widehat{L}_i,\R^3,\geuc\right]\left|_{ \pre{\Sigma}_i } \right.
        \mbox{ for each } i \in \{0,1\},
    \end{aligned}
  \end{equation}
where we recall from Section \ref{nc}
that the function $d\left[\widehat{L}_i,\R^3,\geuc\right]$ 
measures the distance (relative to the Euclidean metric)
of its argument from the set $\widehat{L}_i$,
(defined in (\ref{Wmdef}))
of axes of half catenoids attaching to $\Sigma_i$.

In particular
  \begin{equation}
    \pre{\rho}^{-1}(\xx,\yy,\zz)
    =
    d\left[\widehat{L}_0\right](\xx,\yy,\zz)
    =
    d\left[\widehat{L}_1\right](\xx,\yy,\zz)
    =\sqrt{(\xx-1)^2+\yy^2}
    \mbox{ when } \sqrt{(\xx-1)^2+\yy^2} \leq \frac{1}{10m},
  \end{equation}
and $\pre{\rho}^{-1}$ provides a pointwise measure of the
natural local scale
of the pre-initial surface,
satisfying $\pre{\rho}(\cate(t,\vartheta))=\tau^{-1} \sech t$
on $\cate^{-1}\left(\pre{\Sigma}_{cat}\right) \subset \cyl_a^+$ and transitioning,
on small annuli, to take the constant value $m$ away from
$\Grp[m] \pre{\Sigma}_{cat}$.
Clearly every element of $\Grp[m]$ preserves $\pre{\chi}$.

For each $i \in \{0,1\}$ we also define the function
  \begin{equation}
  \label{rhohat}
    \begin{aligned}
      &\rhohat_i: \Bbar^2 \backslash L_i \to \R \mbox{ by} \\
      &\rhohat_i
      :=
      m \cutoff{\frac{1}{10m}}{\frac{1}{5m}} \circ d\left[L_i,\Bbar^2,\geuc\right]
        + \frac{1}{d\left[L_i,\Bbar^2,\geuc\right]}
        \cutoff{\frac{1}{5m}}{\frac{1}{10m}} \circ d\left[L_i,\Bbar^2,\geuc\right].
    \end{aligned}
  \end{equation}

Analogously to $(\pre{s},\pre{\sigma})$ we introduce the coordinate system
$\left(\pre{s}_{\pre{\chi}},\pre{\sigma}_{\pre{\chi}}\right)$,
well-defined on a sufficiently small neighborhood of $\partial \pre{\Sigma}$,
so that for any point $p$ in this neighborhood
$\pre{\sigma}_{\pre{\chi}}(p)$ is the $\pre{\chi}$
distance from $\partial \pre{\Sigma}$ to $p$ and
$\pre{s}_{\pre{\chi}}(p)$ is
the $\pre{\chi}$ distance along $\partial \pre{\Sigma}$
from an arbitrarily fixed reference point to
the $\pre{\chi}$ nearest-point projection of $p$ onto $\partial \pre{\Sigma}$.
Note that in particular
  \begin{equation}
  \label{prechiconormal}
    \left.\partial_{\pre{\sigma}}\right|_{\partial \pre{\Sigma}}
    =
    \pre{\rho}\left.\partial_{\pre{\sigma}_{\pre{\chi}}}\right|_{\partial \pre{\Sigma}}.
  \end{equation}

In Proposition \ref{preinitest} below
we estimate some of the above quantities.
To state the estimate of $\pre{H}$
we define
$\pre{\wbar} \in C^\infty(\pre{\Sigma})$---which will also play
an important, closely related role in the next section---to be the unique $\Grp[m]$-odd function
supported on $\pre{\Sigma}_0$
and having restriction to $\digamma_0$
  \begin{equation}
  \label{prewbar}
    \pre{\wbar} \left|_{\digamma_0} \right.
    :=
    \pre{\pi}_0^*\rhohat_0^{-2}\Delta_{\geuc}
      \left(\cutoff{\frac{1}{10m}}{\frac{1}{20m}} \circ \sqrt{(1-\rr)^2+\theta^2}\right).
  \end{equation}

\begin{prop}[Estimates on the pre-initial surfaces]
\label{preinitest}
There exists $\epsilon>0$
and, for each nonnegative integer $k$,
there exists a constant $C(k)>0$ such that given any $c>0$,
for every positive integer $m$ sufficiently large
in terms of $c$ and for every $\zeta,\xi \in [-c,c]$
the following estimates hold on the pre-initial surface $\pre{\Sigma}[m,\zeta,\xi]$:
  \begin{enumerate}[(i)]
    \item \label{precatchiest}
      $\norm{\pre{\kappa}^*\pre{\chi}-(dt^2+d\vartheta^2):
      C^k(\cyl_a^+, dt^2+d\vartheta^2)} \leq C(k)m^{-1}$,

    \item \label{precatAest}
      $\norm{\pre{\kappa}^*\pre{A}-\tau(dt^2-d\vartheta^2):
      C^k(\cyl_a^+, dt^2+d\vartheta^2)} \leq C(k)m^{-1}$,

    \item \label{precatHest} $\pre{H}|_{\pre{\Sigma}_{cat}}=0$,

    \item \label{preintplanar0est}
      on $\pre{\Sigma}_0 \cap \left\{\pre{\pi}_0^*\rr<1-\frac{1}{m}\right\}$ we have
      $\pre{\chi}=m^2\pre{\pi}_0^*(d\xx^2+d\yy^2)$,
      $\pre{A}=0$, and $\pre{H}=0$,

    \item \label{preintplanar1est}
      on $\pre{\Sigma}_1 \cap \left\{\pre{\pi}_0^*\rr<1-\frac{1}{m}\right\}$ we have
      $\pre{\chi}=\frac{m^2}{1-(\zz_0+a\tau)^2}\pre{\pi}_1^*(d\xx^2+d\yy^2)$,
      $\pre{A}=0$, and $\pre{H}=0$,

    \item \label{preplanarchiest}
      $\norm{\pre{\chi}-\pre{\pi}_i^*\rhohat_i^2(d\xx^2+d\yy^2)
      : C^k(\pre{\Sigma}_i, \pre{\chi})}
      \leq C(k) \left(\sech^2 \frac{m}{4}+m\tau\right)$ for each $i \in \{0,1\}$,

    \item \label{preAest}
      $\norm{\pre{A}:
      C^k(\pre{\Sigma},\pre{\chi})}
      \leq C(k)(1+c)\tau$,

    \item \label{preHest}
      $\norm{\pre{H}-\pre{\rho}^2\tau \xi \pre{\wbar}:
      C^k(\pre{\Sigma},\pre{\chi})}
      \leq C(k) m^2\tau$,

    \item \label{prerhoest}
      $\norm{\pre{\rho}: C^k(\pre{\Sigma},\pre{\chi},\pre{\rho})}
        +\norm{\pre{\rho}^{-1}: C^k(\pre{\Sigma},\pre{\chi},\pre{\rho})}
      \leq C(k)$,

    \item \label{prezest}
      $\norm{\zz|_{\pre{\Sigma}}:
      C^k(\pre{\Sigma},\pre{\chi})}
      \leq C(k)\left(\sech^2 \frac{m}{4}+m^{-1}+m\tau\right)$,

    \item \label{precatbdymet}
      $\norm{\pre{\kappa}^*\pre{\chi}|_{\partial \pre{\Sigma}}-dt^2:
        C^k(\{\abs{\theta}=\pi/2\}, dt^2)}
      \leq C(k)m^{-1}$,

    \item \label{preplanarbdymet}
      $\norm{\pre{\chi}|_{\partial \pre{\Sigma}}-\pre{\pi}_i^*\rhohat_i^2 d\theta^2:
        C^k\left(\partial_{\Sph^2}\pre{\Sigma}_i,\pre{\chi}|_{\partial \pre{\Sigma}}\right)}
      \leq C(k) \left(\sech^2 \frac{m}{4}+m\tau \right)$ for each $i \in \{0,1\}$,

    \item \label{precatconormal}
      $\norm{\pre{\rho}^{-1}\partial_{\pre{\sigma}}
        +(\sgn \vartheta)\pre{\kappa}_*\partial_{\vartheta}:
      C^k\left(\partial_{\Sph^2}\pre{\Sigma}_{cat},
        \pre{\chi}|_{\partial_{\Sph^2}\pre{\Sigma}_{cat}}\right)}
      \leq C(k)m^{-1}$,

    \item \label{preplanarconormal}
      $\norm{\pre{\pi}_{i*}\pre{\rho}^{-1}\partial_{\pre{\sigma}}+\rhohat_i^{-1}\partial_{\rr}:
      C^k\left(\partial_{\Sph^2}\pre{\Sigma}_i,\rhohat_i^2 d\theta^2\right)}
      \leq C(k) \left(\sech^2 \frac{m}{4}+m\tau\right)$ for each $i \in \{0,1\}$,
      
    \item \label{prebdy}
      $\left(\pre{s}_{\pre{\chi}},\pre{\sigma}_{\pre{\chi}}\right)$
      is a $C^k$ coordinate system on $\left\{\pre{\sigma}_{\pre{\chi}}<\epsilon\right\}$,

    \item \label{pressigmaest}
      $\norm{d\pre{s}_{\pre{\chi}}:
      C^k\left(\left\{\pre{\sigma}_{\pre{\chi}}<\epsilon\right\},\pre{\chi}\right)}
      +
      \norm{\pre{\sigma}_{\pre{\chi}}:
      C^k\left(\left\{\pre{\sigma}_{\pre{\chi}}<\epsilon\right\},\pre{\chi}\right)}
      \leq C(k)$, and

    \item \label{preThetaest}
      $\norm{\pre{\Theta}:
      C^k\left(\partial \pre{\Sigma},\pre{\chi}|_{\partial \pre{\Sigma}},
       m\tau\left(1+m\tau^2\pre{\rho}^2|_{\partial \pre{\Sigma}}\right)\right)}
      \leq C(k)$.
  \end{enumerate}
\end{prop}

\begin{proof}
Using (\ref{phi01}) and (\ref{lambdadef}) it is elementary to check that
for each nonnegative integer $k$ there is a constant $C(k)>0$ so that
  \begin{equation}
    \norm{\lambda-1:C^k([-a,a],dt^2)} \leq C(k)m^{-1}.
  \end{equation}
Consequently, from (\ref{cate}) and (\ref{prekappadef}), we obtain
  \begin{equation}
  \label{kappacatediff}
    \norm{\pre{\kappa}-\cate: C^k\left(\cyl_a^+,dt^2+d\vartheta^2,\tau \cosh t\right)}
    \leq
    C(k)m^{-1}.
  \end{equation}
Using also
  \begin{equation}
  \label{categeom}
    \begin{aligned}
      &\pre{\rho}(\cate(t,\vartheta))=\tau^{-1} \sech t
        \mbox{ for } (t,\vartheta) \in \cate^{-1}\left(\pre{\Sigma}_{cat}\right)
        \subset \cyl_a^+, \\
      &\pre{\rho}(\pre{\kappa}(t,\vartheta))=\tau^{-1} \sech t
        \mbox{ for } (t,\vartheta) \in \cyl_a^+, \\
      &\cate|_{\cate^{-1}\left(\pre{\Sigma}_{cat}\right)}^*\pre{\chi}=dt^2+d\vartheta^2,
        \mbox{ and} \\
      &\cate|_{\cate^{-1}\left(\pre{\Sigma}_{cat}\right)}^*\pre{A}=\tau(dt^2-d\vartheta^2),
    \end{aligned}
  \end{equation}
from (\ref{kappacatediff}) we secure items (\ref{precatchiest}) and (\ref{precatAest}).
Item (\ref{precatHest}) is clear since $\pre{\Sigma}_{cat}$ is a subset of a catenoid.

Items (\ref{preintplanar0est}) and (\ref{preintplanar1est}) follow from
(\ref{prepidef}) and the observation $\pre{\Sigma}$ is exactly planar
on the regions in question.
Since we already have estimates for $\pre{A}$ and $\pre{H}$
on $\pre{\Sigma}_i \cap \pre{\Sigma}_{cat}$ (for each $i \in \{0,1\}$),
to confirm the global estimates (\ref{preAest}) and (\ref{preHest})
it remains only to check them on the region
$\pre{\Sigma}_i \backslash \left(\pre{\Sigma}_{cat} \cup \{1-\pre{\pi}_i^*\rr>1/m\}\right)$.
To confirm the estimate (\ref{preplanarchiest})
it likewise remains to check it on this region,
but we must also translate (\ref{precatchiest})
into an estimate for $\pre{\pi}_i^*\pre{\chi}$
on $\pre{\Sigma}_i \cap \pre{\Sigma}_{cat}$.
 
With these purposes in mind, for each $i \in \{0,1\}$
we introduce the map $\pihat_i: \pre{\Sigma}_i \to \Bbar^2$
defined by
  \begin{equation}
  \label{pihat}
    \pihat_i(\xx,\yy,\zz)=(\xx,\yy),
  \end{equation}
so that $\pihat_i$ is a diffeomorphism onto its image,
a subset of $\Bbar^2$.
It is then evident upon comparing (\ref{prerho})
with (\ref{rhohat}) that
  \begin{equation}
    \pre{\rho}=\pihat_i^*\rhohat_i.
  \end{equation}
Using (\ref{cate}) we also find
  \begin{equation}
  \label{overlap}
    \cate^*\pre{\rho}^{-2}\pihat_i^*(d\xx^2+d\yy^2)=\tanh^2 t \, dt^2 + d\vartheta^2,
  \end{equation}
verifying (\ref{preplanarchiest}) on $\pre{\Sigma}_i \cap \pre{\Sigma}_{cat}$.
From the definition (\ref{preinitsdef}) of the pre-initial surface $\pre{\Sigma}$
we see that the remaining region
$\pre{\Sigma}_i \backslash \left(\pre{\Sigma}_{cat} \cup \{1-\pre{\pi}_i^*\rr>1/m\}\right)$
is contained in the union of $2m$ congruent copies of the union of the graphs of
the functions $\phi_0,\phi_1: \{\rr \geq \tau\} \to \R$
(defined in (\ref{phi01}))
over the annulus $\Lambda=\{\rr \in \left(\frac{1}{20m},\frac{1}{m}\right)\}$.

If we write $g_i$ for the metric induced by $\geuc$
on the graph of $\phi_i$ over $\Lambda$
and if we write $A_i$ and $H_i$ for the corresponding second fundamental form
and mean curvature, both relative to the upward unit normal,
then
  \begin{equation}
    \begin{aligned}
      &g_i=d\xx^2+d\yy^2+d\phi_i^2 \mbox{ and} \\
      &(A_i)_{ab}=\frac{(\phi_i)_{,ab}}{\sqrt{1+\abs{d\phi_i}^2}}, \mbox{ so} \\
      &\norm{D^k\left(H_i - \Delta \phi_i\right)}_{C^0(\Lambda,d\xx^2+d\yy^2)}
        \lesssim \sum_{k_1+k_2+k_3=k} \norm{D^{2+k_1}\phi_i}_{C^0}
          \norm{D^{1+k_2}\phi_i}_{C^0}
          \norm{D^{1+k_3}\phi_i}_{C^0}.
    \end{aligned}
  \end{equation}
Since (recalling (\ref{prerho}))
  \begin{equation}
  \label{rhoest}
    \norm{m\pre{\rho}^{-1} : C^k(\Lambda,m^2(d\xx^2+d\yy^2)}
      +\norm{m^{-1}\pre{\rho} : C^k(\Lambda,m^2(d\xx^2+d\yy^2))}
    \leq C(k),
  \end{equation}
the $C^k$ norms induced over $\Lambda$ by $\pre{\rho}^2(d\xx^2+d\yy^2)$ and $m^2(d\xx^2+d\yy^2)$
are equivalent, each bounded by the other times a constant independent of $m$.

We estimate
  \begin{equation}
  \label{phiests}
    \begin{aligned}
      &\norm{D^k\phi_1}_{C^0(\Lambda,d\xx^2+d\yy^2)} \leq C(k) m^k\tau
        \mbox{ for $k \geq 1$}, \\
      &\norm{D^k\phi_0}_{C^0(\Lambda,d\xx^2+d\yy^2)} \leq C(k) m^k\tau(\abs{\xi}+1)
        \mbox{ for $k \geq 1$}, \\
      &\norm{\Delta \phi_1}_{C^k(\Lambda,d\xx^2+d\yy^2)} \leq C(k) m^{2+k}\tau
        \mbox{ for $k \geq 0$}, \mbox{ and} \\
      &\norm{\Delta \phi_0 - \tau \xi 
          \Delta \left(\cutoff{\frac{1}{10m}}{\frac{1}{20m}} \circ \rr\right)}_{
            C^k(\Lambda,d\xx^2+d\yy^2)}
        \leq C(k) m^{2+k}\tau \mbox{ for $k \geq 0$}.
    \end{aligned}
  \end{equation}
Defining $\overline{v},\widetilde{\overline{v}}: \Bbar^2 \to \R$
to be the unique $\Grp[m]$-odd functions having restrictions to $W[m]$
  \begin{equation}
    \begin{aligned}
      &\overline{v}|_{W[m]}
      :=
      \cutoff{\frac{1}{10m}}{\frac{1}{20m}} \circ \sqrt{(1-\rr)^2+\theta^2}
        \mbox{ and} \\
      &\widetilde{\overline{v}}|_{W[m]}
      :=
      \cutoff{\frac{1}{10m}}{\frac{1}{20m}} \circ \sqrt{(1-\xx)^2+\yy^2},
    \end{aligned}
  \end{equation}
we also find
  \begin{equation}
  \label{vcomp}
    \norm{\overline{v}-\widetilde{\overline{v}}: C^k\left(\Bbar^2,m^2(d\xx^2+d\yy^2)\right)}
    \leq
    C(k)m^{-1}.
  \end{equation}

Applying (\ref{overlap}),(\ref{phiests}), and (\ref{vcomp})
in conjunction with (\ref{prewbar}) and the already established items
(\ref{precatAest}), (\ref{precatHest}), (\ref{preintplanar0est}),
and (\ref{preintplanar1est}),
we therefore obtain
  \begin{equation}
  \label{pihatests}
    \begin{aligned}
      &\norm{\pre{\chi}-\pre{\rho}^2\pihat_i^*(d\xx^2+d\yy^2)
      : C^k(\pre{\Sigma}_i, \pre{\chi})}
      \leq C(k) \left(\sech^2 \frac{m}{4}+(1+\xi^2)m^2\tau^2\right), \\
      &\norm{\pre{A}:
      C^k(\pre{\Sigma},\pre{\chi})}
      \leq C(k)(1+\abs{\xi})\tau, \mbox{ and} \\
      &\norm{\pre{H}-\pre{\rho}^2\tau \xi \pihat_0^*\pre{\pi}_0^{*-1}\pre{\wbar}:
      C^k(\pre{\Sigma},\pre{\chi})}
      \leq C(k)\left(1+m^2\tau^2\abs{\xi}^3\right) m^2\tau.
    \end{aligned}
  \end{equation}

To complete the proofs of items
(\ref{preplanarchiest}), (\ref{preAest}), and (\ref{preHest})
from the three estimates in (\ref{pihatests})
we need to compare the maps $\pre{\pi}_i$ and $\pihat_i$,
so we will need estimates for the height function $\zz|_{\pre{\Sigma}}$.
First we use (\ref{phi01}) to obtain the global $C^0$ estimate
  \begin{equation}
  \label{phiC0est}
  \abs{\zz}
  \leq
  \abs{\zz_0}+\tau \phi_{cat}\left(\frac{1}{10m}\right)
  \leq
  (C+\abs{\zeta}+\abs{\xi}+m)\tau.
  \end{equation}
Next we observe $\zz$ is constant on
$\pre{\Sigma}_i \cap \{\pihat_i^*\rr<1-1/m\}$
for each $i \in \{0,1\}$,
while on $\pre{\Sigma}_{cat}$ we have
  \begin{equation}
    \zz(\cate(t,\vartheta))=\zz_0+\tau t,
  \end{equation}
and elsewhere $\zz$ is controlled by the estimates for $\phi_0$ and $\phi_1$
in (\ref{phiests}).
Thus (using also item (\ref{precatchiest})
and the first estimate in (\ref{pihatests}))
we are able to verify (\ref{prezest}).
It follows immediately,
referring to the definitions (\ref{prepidef}) and (\ref{pihat}),
that for each $i \in \{0,1\}$
  \begin{equation}
    \norm{\pre{\pi}_i-\pihat_i:
      C^k(\pre{\Sigma}_i,\pihat_i^*\pre{\rho}^2(d\xx^2+d\yy^2))}
    \leq
    C(k)\left(\sech^2 \frac{m}{4}+[\abs{\zeta}+\abs{\xi}+m]\tau\right),
  \end{equation}
which in conjunction with (\ref{pihatests}),
as well as definition (\ref{prewbar}), completes the proof of
(\ref{preplanarchiest}), (\ref{preAest}), and (\ref{preHest}).
Item (\ref{prerhoest}) is now obvious from
(\ref{categeom}) and (\ref{rhoest})
together with items (\ref{precatchiest}) and (\ref{preplanarchiest}).

Now we turn to the boundary geometry.
Items (\ref{precatbdymet}) and (\ref{preplanarbdymet})
follow directly from
(\ref{precatchiest}) and (\ref{preplanarchiest})
respectively,
since $\pre{\kappa}(\{\abs{t}=\pi/2\})=\partial_{\Sph^2}\pre{\Sigma}_{cat}$
and $\pre{\pi}_i\partial_{\Sph^2}\pre{\Sigma}_i \subset \partial \Bbar^2$.
Items (\ref{precatconormal}) and (\ref{preplanarconormal})
then follow from these last two, in view of (\ref{prechiconormal}).
Next, from (\ref{precatchiest}) we see that
$\left(\pre{s}_{\pre{\chi}},\pre{\sigma}_{\pre{\chi}}\right)$
is a smooth coordinate system on
$\pre{\kappa}(\cyl_a^+ \cap \{\vartheta \neq 0\})$,
where moreover
  \begin{equation}
  \begin{aligned}
    &\norm{d\pre{\kappa}^*\pre{s}_{\pre{\chi}}-dt:
      C^k(\cyl_a^+ \cap \{\vartheta \neq 0\}, dt^2+d\vartheta^2)}
    \leq
    C(k)m^{-1} \mbox{ and} \\
    &\norm{\pre{\kappa}^*\pre{\sigma}_{\pre{\chi}}-(\pi/2-\abs{\vartheta}):
      C^k(\cyl_a^+ \cap \{\vartheta \neq 0\}, dt^2+d\vartheta^2)}
    \leq
    C(k)m^{-1},
  \end{aligned}
  \end{equation}
verifying (\ref{prebdy}) and (\ref{pressigmaest}) on this region.
Similarly, from (\ref{preplanarchiest}) we see that
$\left(\pre{s}_{\pre{\chi}},\pre{\sigma}_{\pre{\chi}}\right)$
is also a smooth coordinate system
on the region
$M:=\left\{d\left[\widehat{L}_0,\R^3,\geuc\right] \geq \frac{1}{m} \right\}$
(recalling (\ref{Wmdef}) and the last paragraph of Section \ref{nc})
on which $\rho=m$ and moreover for each $i \in \{0,1\}$
  \begin{equation}
  \begin{aligned}
    &\norm{d\pre{s}_{\pre{\chi}}-m^2\pre{\pi}_i^*d\theta^2:
      C^k(\pre{\Sigma}_i \cap \{\rho=m\}, \pre{\chi})}
    \leq
    C(k) \left(\sech^2 \frac{m}{4}+m\tau\right) \mbox{ and} \\
    &\norm{\pre{\sigma}_{\pre{\chi}}-m^2\pre{\pi}_i^*(1-\rr):
      C^k(\pre{\Sigma}_i \cap \{\rho=m\}, \pre{\chi})}
    \leq
    C(k) \left(\sech^2 \frac{m}{4}+m\tau \right),
  \end{aligned}
  \end{equation}
verifying (\ref{prebdy}) and (\ref{pressigmaest}) on this region as well.

To complete the proof of (\ref{prebdy}) and (\ref{pressigmaest})
we consider the region
$N_o:=\pre{\Sigma} \cap \left\{\frac{1}{40m}<d\left[\widehat{L}_0,\R^3,\geuc\right]
  <\frac{3}{m}\right\}$,
and its compact subset
$N:=\pre{\Sigma} \cap \left\{\frac{1}{30m}
\leq d\left[\widehat{L}_0,\R^3,\geuc\right] \leq \frac{2}{m}\right\}$
so that
$\pre{\Sigma}=\Grp[m] \left(N \cup M \cup \pre{\Sigma}_{cat}\right)$
and in fact, for $m$ sufficiently large,
every closed disc in $\pre{\Sigma}$
of $\pre{\chi}$ radius $\frac{1}{100}$ is completely contained
in $M$, $N$, or $\pre{\Sigma}_{cat}$.
Since $(N,\pre{\chi})$ has geometry uniformly bounded in $m$
(actually converging to a half cylinder of order-$1$ length
attached to a half disc of order-$1$ radius),
this ensures the existence of $\epsilon>0$ so that
(\ref{prebdy}) and (\ref{pressigmaest}) hold globally.

Finally we estimate the intersection angle between $\pre{\Sigma}$ and $\partial \Bbar^3$.
Working first on $\pre{\Sigma}_{cat}$, from (\ref{cate}) we have
  \begin{equation}
     \pre{\nu}(\cate(t,\vartheta))=\left(\sech t \cos \vartheta, -\sech t \sin \vartheta, \tanh t \right)
  \end{equation}
and obviously
  \begin{equation}
    \pre{X}(\cate(t,\vartheta))
    =(1-\tau \cosh t \cos \vartheta, \tau \cosh t \sin \vartheta, \tau t + \zz_0),
  \end{equation}
so
  \begin{equation}
      \left(\cate^*\left\langle \pre{X}, \pre{\nu} \right\rangle\right)(t,\vartheta)
      =
      \langle \pre{X}(\cate(t,\vartheta)), \pre{\nu}(\cate(t,\vartheta)) \rangle
      =
      \sech t \cos \vartheta - \tau + \tau t \tanh t + \zz_0 \tanh t.
  \end{equation}
Restricting to $\partial_{\Sph^2} \pre{\Sigma}_{cat}$ and using (\ref{catbdycurve}) we find
  \begin{equation}
    \cate^*\pre{\Theta}(t)
    =
    \frac{1}{2}\tau^{-1} (\zz_0 + \tau t)^2 \sech^2 t - \frac{1}{2}\tau + (\zz_0+\tau t) \tanh t,
  \end{equation}
which,
using also (\ref{precatbdymet}),
proves (\ref{preThetaest}) on $\pre{\Sigma}_{cat}$.

On the other hand,
by (\ref{pihat}) and (\ref{preinitsdef}),
on $P_0:=\pre{\Sigma}_0 \cap \left\{ 0 \leq \theta \leq \frac{\pi}{2m} \right\}$
and $P_1:=\pre{\Sigma}_i \cap \left\{ 0 \leq \theta \leq \frac{\pi}{m} \right\}$
we have
  \begin{equation}
    \pre{\nu}|_{P_i}(\pihat_i^{-1}(\xx,\yy))
    =
    (-1)^i\frac{\left(\phi_{i,\xx}(\xx-1,\yy), \phi_{i,\yy}(\xx-1,\yy), -1 \right)}
      {\sqrt{1+\abs{(\nabla \phi_i)(\xx-1,\yy)}^2}},
  \end{equation}
and
  \begin{equation}
    X|_{P_i}(\pihat_i^{-1}(\xx,\yy))=(\xx,\yy,\phi_i(\xx-1,\yy)),
  \end{equation}
so
  \begin{equation}
  \label{planarpreTheta}
    {{}\pihat_i^{-1}}^*\left\langle \pre{X}, \pre{\nu} \right\rangle |_{P_i}(\xx,\yy)
    =
    (-1)^i\frac{\xx \phi_{i,\xx}(\xx-1,\yy) + \yy \phi_{i,\yy}(\xx-1,\yy) - \phi_i(\xx-1,\yy)}
      {\sqrt{1+\abs{(\nabla \phi_i)(\xx-1,\yy)}^2}}.
  \end{equation}

By applying the estimates (\ref{phiC0est}), (\ref{phiests}),
and item (\ref{preplanarbdymet}) of the proposition
to (\ref{planarpreTheta})
we obtain estimates for $\pre{\Theta}$ on the remainder of $\partial_{\Sigma}$.
The derivative estimate for $\phi_0$ in (\ref{phiests}),
however, is not refined enough to ensure
that  (\ref{preThetaest}) holds on $\partial_{\Sph^2} \Sigma_0$.
In order secure an estimate independent of the parameter $\xi$
we use (\ref{vcomp}) along with the fact that
  \begin{equation}
    \partial_{\rr}\overline{v}|_{\partial \Bbar^2}=0;
  \end{equation}
it follows that
  \begin{equation}
    \norm{\xx \phi_{0,\xx}(\xx-1,\yy)+\yy \phi_{0,\yy}(\yy-1,\yy)
    : C^k(\partial W[m], m^2 d\theta^2)}
    \leq
    C(k)m\tau,
  \end{equation}
which we use in place of (\ref{phiests})
to finish the proof.
\end{proof}

\begin{remark}
\label{precty}
It is clear from their construction
that for each fixed $m$
the pre-initial surfaces $\pre{\Sigma}[m,\zeta,\xi]$
depend smoothly on the parameters $\zeta$ and $\xi$
in the sense that there exists a smooth map
$\pre{Y}_m: \R \times \R \times \pre{\Sigma}[m,0,0] \to \Bbar^3$
such that for each $\zeta,\xi \in \R$
the map $\pre{Y}_m(\zeta,\xi,\cdot)$ is an embedding
with image $\pre{\Sigma}[m,\zeta,\xi]$.
This allows us to identify functions on pre-initial surfaces
with the same $m$ value but different parameter values.
\end{remark}

\section{The initial surfaces}
\label{S:i}

We will now bend $\pre{\Sigma}$ near its boundary
to make it intersect the sphere $\partial \Bbar^3$ orthogonally.
We do this by picking a small function on $\pre{\Sigma}$ whose graph has
the desired property. Namely we define $\pre{u} \in C^\infty(\pre{\Sigma})$ by
  \begin{equation}
  \label{preu}
    \pre{u}\left(\pre{s}_{\pre{\chi}},\pre{\sigma}_{\pre{\chi}}\right)
    =
    -\pre{\rho}^{-1}\pre{\sigma}_{\pre{\chi}}
    \frac{\pre{\Theta}\left(\pre{s}_{\pre{\chi}}\right)}
      {\sqrt{1-\pre{\Theta}^2\left(\pre{s}_{\pre{\chi}}\right)}}
    \cutoff{\frac{\epsilon}{2}}{\frac{\epsilon}{4}}\left(\pre{\sigma}_{\pre{\chi}}\right),
  \end{equation}
with $\epsilon$ as in Proposition \ref{preinitest},
so that the coordinates $\left(\pre{s}_{\pre{\chi}},\pre{\sigma}_{\pre{\chi}}\right)$
are well-defined and smooth on the support of the cut-off function appearing in
the definition of $\pre{u}$
and we understand that $\pre{u}$ identically vanishes elsewhere.
Then we define the corresponding deformation
$\pre{X}_{\pre{u}}[m,\zeta,\xi]: \pre{\Sigma}[m,\zeta,\xi] \to \R^3$ of $\pre{X}$ by
  \begin{equation}
  \label{preXu}
    \pre{X}_{\pre{u}}[m,\zeta,\xi](p)=\pre{X}[m,\zeta,\xi](p)+\pre{u}(p)\pre{\nu}(p).
  \end{equation}
and the resulting initial surface $\Sigma[m,\zeta,\xi]$
as the image of $\pre{\Sigma}[m,\zeta,\xi]$
under $\pre{X}_u$:
  \begin{equation}
    \Sigma[m,\zeta,\xi]=\pre{X}_u[m,\zeta,\xi](\pre{\Sigma}[m,\zeta,\xi]).
  \end{equation}
In accordance with Remark \ref{suppress}
we will routinely abbreviate $\Sigma[m,\zeta,\xi]$ by $\Sigma$.

We write $X: \Sigma \to \R^3$ for the inclusion map of $\Sigma$ in $\R^3$,
$g=X^*\geuc$ for the induced metric,
$\nu: \Sigma \to \R^3$ for the unit normal which points downward at the origin $(0,0,0) \in \Sigma$,
$A_{ab}=-\nu_{i|a}X^i_{\;,b}$ for the corresponding second fundamental form,
and $H=-g^{ab}A_{ab}$
for the corresponding mean curvature.
In the proof of the following proposition we will establish that
$\pre{X}_{\pre{u}}$ defines a smooth diffeomorphism
from $\pre{\Sigma}$ to its image $\Sigma$.
Whenever convenient
we will permit ourselves to shrink the target of $\pre{X}_{\pre{u}}$
as originally defined
from $\R^3$ to $\Sigma$
to yield a diffeomorphism $\pre{X}_{\pre{u}}: \pre{\Sigma} \to \Sigma$.

With this interpretation in mind we now define
the functions $\rho, \wbar \in C^\infty(\Sigma)$,
the metric $\chi \in C^\infty\left(T^*\Sigma^{\otimes 2}\right)$,
the catenoidal region $\Sigma_{cat} \subset \Sigma$,
along with the spherical part of its boundary $\partial_{\Sph^2}\Sigma_{cat}$,
the diffeomorphism $\kappa: \cyl_a^+ \to \Sigma_{cat}$,
and, for each $i \in \{0,1\}$,
the planar region $\Sigma_i \subset \Sigma$,
along with the spherical part of its boundary $\partial_{\Sph^2}\Sigma_i$,
and the diffeomorphism $\pi_i: \Sigma_i \to \Bbar^2$
by
  \begin{equation}
  \label{uncheck}
    \begin{aligned}
      &\rho:={{}\pre{X}_u^{-1}}^*\pre{\rho},
        \quad \chi=\rho^2g:={{}\pre{X}_u^{-1}}^*\pre{\chi},
        \quad \wbar:={{}\pre{X}_u^{-1}}^*\pre{\wbar}, \\
      &\Sigma_{cat}:=\pre{X}_{\pre{u}}(\pre{\Sigma}_{cat}),
        \quad \partial_{\Sph^2} \Sigma_{cat}:=\partial \Sigma \cap \partial \Sigma_{cat},
        \quad \kappa:=\pre{X}_{\pre{u}} \circ \pre{\kappa}, \\
      &\Sigma_i:=\pre{X}_{\pre{u}}(\pre{\Sigma}_i),
        \quad \partial_{\Sph^2} \Sigma_i:=\partial \Sigma \cap \partial \Sigma_i, \mbox{ and }
        \pi_i:=\pre{\pi}_i \circ \pre{X}_{\pre{u}}^{-1}.
    \end{aligned}
  \end{equation}
\begin{defn}
We also introduce the coordinates $(s,\sigma)$
on a sufficiently small neighborhood in $\Sigma$ of $\partial \Sigma$,
so that for each $p$ in this neighborhood
$\sigma(p)$ is the $g$ distance in $\Sigma$ from $\partial \Sigma$ to $p$,
while $s(p)$ is the $g$ distance along $\partial \Sigma$
from an arbitrarily fixed reference point in $\partial \Sigma$
to the $g$-nearest-point projection of $p$ onto $\partial \Sigma$.
\end{defn}
\begin{remark}
\label{conormal}
In particular $\partial_\sigma\left|_{\partial \Sigma}\right.$
is the inward unit conormal vector field,
relative to the ambient Euclidean metric, along $\partial \Sigma$.
\end{remark}

\begin{prop}[Estimates on the initial surfaces]
\label{initest}
Given $c>0$, there exists $m_0>0$ such that
for every $m>m_0$ and for every $\zeta,\xi \in [-c,c]$
the set $\Sigma$ is a smooth
orientable
manifold with boundary,
smoothly embedded by $X$ in $\Bbar^3$;
$\Sigma$ is invariant under $\Grp[m]$;
$\Sigma$ has genus $m-1$; $\Sigma$ has connected boundary $\partial \Sigma$,
which is smoothly embedded in $\partial \Bbar^3$;
and $\Sigma$ intersects $\partial \Bbar^3$ orthogonally along $\partial \Sigma$.
Moreover, for each $k \in \Z^+$ there is a constant $C(k)>0$
such that given any $c>0$ there exists $m_0>0$ such that
for every integer $m>m_0$ and for every $\zeta,\xi \in [-c,c]$
the following estimates hold:
  \begin{enumerate}[(i)]
    \item \label{preuest}
    $\norm{\pre{u}: C^k\left(\pre{\Sigma},\pre{\chi},
    \tau\left[m\pre{\rho}^{-1} + m^2\tau\right]\right)}
    \leq
    C(k)$.

    \item \label{intplanar0est}
      on $\Sigma_0 \cap \left\{\pi_0^*\rr<1-\frac{1}{m}\right\}$ we have
      $\chi=m^2\pi_0^*(d\xx^2+d\yy^2)$,
      $A=0$, and $H=0$,

    \item \label{intplanar1est}
      on $\Sigma_1 \cap \left\{\pi_0^*\rr<1-\frac{1}{m}\right\}$ we have
      $\chi=\frac{m^2}{1-(\zz_0+a\tau)^2}\pi_1^*(d\xx^2+d\yy^2)$,
      $A=0$, and $H=0$,

    \item \label{catchiest} 
      $\norm{\kappa^*\chi-(dt^2+d\vartheta^2):
      C^k(\cyl_a^+, dt^2+d\vartheta^2)} \leq C(k)m^{-1}$,

    \item \label{planarchiest} 
      $\norm{\chi-\pi_i^*\rhohat_i^2(d\xx^2+d\yy^2)
      : C^k(\Sigma_i, \chi)}
      \leq C(k)\left(\sech^2 \frac{m}{4}+m^2\tau\right)$ for each $i \in \{0,1\}$,

    \item \label{Aest}
      $\norm{A: C^k(\Sigma,\chi)} \leq C(k)(1+\abs{\xi})\tau$
      and $\norm{\rho^{-2}\abs{A}^2|_{\Sigma_{cat}}-2\tau^2\rho^2:
        C^k\left(\Sigma_{cat},\chi\right)} \leq C(k)m^2\tau$,

    \item \label{Hest}
      $\norm{\rho^{-2}H-\tau\xi\wbar : C^k(\Sigma,\chi,m\rho^{-1}+m^2\tau)} \leq C(k)\tau$,

    \item \label{catconormalest}
      $\norm{\rho^{-1}\partial_{\sigma}
        +(\sgn \vartheta)\kappa_*\partial_{\vartheta}:
      C^k\left(\partial_{\Sph^2}\Sigma_{cat},
        \chi|_{\partial_{\Sph^2}\Sigma_{cat}}\right)}
      \leq C(k)m^{-1}$, and

    \item \label{planarconormalest}
      $\norm{\pi_{i*}\rho^{-1}\partial_{\sigma}+\rhohat_i^{-1}\partial_{\rr}:
      C^k\left(\partial_{\Sph^2}\Sigma_i, \rhohat_i^2 d\theta^2 \right)}
      \leq C(k) \left(\sech^2 \frac{m}{4}+m^2\tau\right)$ for each $i \in \{0,1\}$.
  \end{enumerate}
\end{prop}

\begin{proof}
From the definition (\ref{preu}) of $\pre{u}$
and from items (\ref{prerhoest}), (\ref{pressigmaest}), and (\ref{preThetaest})
of Proposition \ref{preinitest} we deduce the estimate (\ref{preuest}).
In particular
$\lim_{m \to \infty} \tau^{-1}\norm{\pre{u}}_{C^0}=0$,
uniformly in $\zeta,\xi \in [-c,c]$,
while for large $m$ each magnified pre-initial surface $\tau^{-1}\pre{\Sigma}$
resembles a union of widely separated discs
and truncated half catenoids with waist radii of unit size,
all intersecting the origin-centered sphere of radius $\tau^{-1}$ almost orthogonally.
The embeddedness and other topological assertions then follow
from Proposition \ref{preinittop},
as does the fact that $\Sigma$ doesn't leave the ball,
and the smoothness claims are obvious.
That $\Sigma$ is $\Grp[m]$-invariant follows from
the $\Grp[m]$ invariance of $\pre{\Sigma}$,
the definition of $\pre{X}_u$,
and the fact that $\pre{u}$ is $\Grp[m]$-odd.

Items (\ref{intplanar0est}) and (\ref{intplanar1est})
are obvious from the support of $\pre{u}$
and the corresponding items, 
(\ref{preintplanar0est}) and (\ref{preintplanar1est})
of Proposition \ref{preinitest}.
Items
(\ref{precatAest}), (\ref{preintplanar0est}), (\ref{preintplanar1est}),
and (\ref{preAest}) of Proposition \ref{preinitest} yield the estimate
  \begin{equation}
  \label{prerhoAest}
    \norm{\pre{\rho} \pre{A} : C^k(\pre{\Sigma}, \pre{\chi})} \leq C(k),
  \end{equation}
where, we emphasize, the constant $C(k)$ is independent of $m$, $\zeta$, $\xi$, and $c$,
while (\ref{preuest}) obviously ensures
  \begin{equation}
  \label{prerhouest}
    \norm{\pre{\rho}\pre{u} : C^k(\pre{\Sigma},\pre{\chi})} \leq C(k)m^2\tau,
  \end{equation}
so from
  \begin{equation}
    \pre{X}_{\pre{u}}^*g_{ab}
    =
    \pre{g}_{ab}-2\pre{u}\pre{A}_{ab}
      +\pre{u}^2\pre{A}_{ac}\pre{A}_{bd}\pre{g}^{cd}
      +\pre{u}_{,a}\pre{u}_{,b}
  \end{equation}
we thus obtain
  \begin{equation}
  \label{chiest}
    \norm{\pre{X}_{\pre{u}}^*\chi_{ab}-\pre{\chi}_{ab}:
      C^k(\pre{\Sigma},\pre{\chi})}
    \leq C(k)m^2\tau
  \end{equation}
(for $C(k)$ independent of $c$ and $m$ and for $m$ large enough compared to $c$),
which in conjunction with
items (\ref{precatchiest}) and (\ref{preplanarchiest})
of Proposition \ref{preinitest}
yields items (\ref{catchiest}) and (\ref{planarchiest}) of the present proposition.
Items (\ref{catconormalest}) and (\ref{planarconormalest})
then follow in turn,
since $\rho^{-1}\partial_{\sigma}$ is the $\chi$ inward unit conormal.

Next, since $\pre{u}|_{\partial \pre{\Sigma}}=0$,
we have
$\partial \Sigma=\partial \pre{\Sigma}$,
$\pre{X}_{\pre{u}}|_{\partial \pre{\Sigma}}=\pre{X}|_{\partial \pre{\Sigma}}$,
and $\pre{u}_{,\pre{s}}|_{\partial \pre{\Sigma}}=0$,
so along the boundary $\pre{X}_{\pre{u}}$ has unit normal
  \begin{equation}
    \pre{\nu}_{\pre{u}}|_{\partial\pre{\Sigma}}
    =
    \left.\frac{\pre{\nu}
      -\pre{u}_{,\pre{\sigma}}\pre{X}_{,\pre{\sigma}}}
      {\sqrt{1+\pre{u}_{,\pre{\sigma}}^2}}\right|_{\partial \pre{\Sigma}}.
  \end{equation}
Using (\ref{preconormal}), (\ref{prechiconormal}), and (\ref{preu}), it follows that
  \begin{equation}
    \langle X,\nu \rangle|_{\partial \Sigma}
    =
    \left.\left\langle \pre{X}_{\pre{u}}, \pre{\nu}_{\pre{u}} \right\rangle\right|_{\partial \pre{\Sigma}}
    =
    \left.\left\langle \pre{X}, \pre{\nu}_{\pre{u}} \right\rangle\right|_{\partial \pre{\Sigma}}
    =
    0,
  \end{equation}
proving that $\Sigma$ intersects $\partial \Bbar^3$ orthogonally.

More generally
  \begin{equation}
    \pre{\nu}_{\pre{u}}
    =
    \frac{\pre{\nu}-\widetilde{g}^{cd}\pre{u}_{,c}
      \left(\pre{X}_{,d}-\pre{u}\pre{A}_{de}\pre{X}_{,f}\pre{g}^{ef}\right)}
      {\sqrt{1+\abs{d\pre{u}}_{\widetilde{g}}^2}},
  \end{equation}
where $\widetilde{g}^{cd}$ is the inverse of the metric
  \begin{equation}
  \label{gtilde}
    \widetilde{g}_{cd}
    =
    \pre{g}_{cd}-2\pre{u}\pre{A}_{cd}+\pre{u}^2\pre{A}_{ce}\pre{A}_{df}\pre{g}^{ef}.
    =
    \pre{X}_{\pre{u}}^*g_{cd}-\pre{u}_{,c}\pre{u}_{,d}
  \end{equation}
Thus $\widetilde{g}|_p$ is the metric on $\pre{\Sigma}$, at point $p$,
induced by the immersion $\pre{X}_{u(p)}=\pre{X}+\pre{u}(p)\pre{\nu}$
and above
$\left.\left(\pre{X}_{,d}-\pre{u}\pre{A}_{de}\pre{X}_{,f}\pre{g}^{ef}\right)\right|_p
  =\left.\pre{X}_{u(p),d}\right|_p$.
We then find
  \begin{equation}
  \label{A}
  \begin{aligned}
    \pre{X}_{\pre{u}}^*A_{ab}
    =
    &\left(1+\abs{d\pre{u}}_{\widetilde{g}}^2\right)^{-1/2}
    \left(\pre{A}_{ab}+\pre{u}_{;ab}
      -\pre{u}\pre{A}_{ac}\pre{A}_{bd}\pre{g}^{cd}
      +\pre{u}\pre{u}_{;c}\widetilde{g}^{cd}\pre{A}_{ab;d} \right. \\
    &\left.+2\pre{u}_{;c}\pre{u}_{(;a}\pre{A}_{b)d}\widetilde{g}^{cd}
      -2\pre{u}\pre{u}_{;c}\pre{u}_{(;a}\pre{A}_{b)e}\pre{A}_{df}\pre{g}^{ef}\widetilde{g}^{cd}
      -\pre{u}^2\pre{u}_{;c}\pre{A}_{ab;e}\pre{A}_{df}\pre{g}^{ef}\widetilde{g}^{cd}\right),
  \end{aligned}
  \end{equation}
where semicolons indicate covariant differentiation relative to $\pre{g}$
and parentheses indicate symmetrization in the indices they enclose,
normalized so as to fix tensors already symmetric in the enclosed indices.

Since $\pre{g}^{-1}=\pre{\rho}^2\pre{\chi}^{-1}$, 
from item (\ref{prerhoest}) of Proposition \ref{preinitest} we get
  \begin{equation}
    \norm{\pre{g}^{-1}: C^k(\pre{\Sigma},\pre{\chi},\pre{\rho}^2)} \leq C(k),
  \end{equation}
which applied to (\ref{gtilde}) in conjunction with
(\ref{prerhoAest}) and (\ref{prerhouest})
implies
  \begin{equation}
  \label{gtildeest}
    \norm{\widetilde{g}^{-1}: C^k(\pre{\Sigma},\pre{\chi},\pre{\rho}^2)} \leq C(k)
  \end{equation}
as well.
Note also that
  \begin{equation}
    \pre{u}_{;ab}
    =\left(D_{\pre{g}}^2\pre{u}\right)_{ab}
    =\left(D_{\pre{\chi}}^2\pre{u}\right)_{ab}
      +\left(D_{\pre{g}}-D_{\pre{\chi}}\right)(du)_{ab}
    =\left(D_{\pre{\chi}}^2\pre{u}\right)_{ab}
      +2\pre{u}_{(,a}(\ln \rho)_{,b)}-\pre{u}_{,c}(\ln \rho)_{,d}\pre{\chi}^{cd}\pre{\chi}_{ab}.
  \end{equation}
Applying (\ref{preuest}), (\ref{prerhoAest}), and (\ref{gtildeest})
to (\ref{A})
we conclude
  \begin{equation}
    \norm{\pre{X}_{\pre{u}}^*A-\pre{A}:
      C^k\left(T^*\pre{\Sigma}^{\otimes 2},\pre{\chi},
      m\rho^{-1}+m^2\tau \right)}
    \leq C(k)\tau,
  \end{equation}
which along with (\ref{chiest})
and items (\ref{preAest}) and (\ref{preHest}) of Proposition \ref{preinitest}
proves items (\ref{Aest}) and (\ref{Hest})
of the present proposition.
\end{proof}

\begin{remark}
\label{cty}
For each fixed $m$
the function $\pre{u}$
clearly depends continuously,
in the sense of Remark \ref{precty},
on the parameters $\zeta$ and $\xi$,
and so in turn
there exists a smooth map
$Y_m: \R \times \R \times \Sigma[m,0,0] \to \Bbar^3$
such that for each $\zeta,\xi \in \R$
the map $Y_m(\zeta,\xi,\cdot)$ is an embedding
with image $\Sigma[m,\zeta,\xi]$.
This allows us to identify functions on initial surfaces
with the same $m$ value but different parameter values.
\end{remark}

\section{Graphs over the initial surfaces}
\label{graphs}

To deform the surface without leaving the ball
it will be useful to introduce on $\Bbar^3$ a metric $\auxm$,
called the \emph{auxiliary metric},
that makes the boundary sphere $\partial \Bbar^3$ totally geodesic
but preserves the intersection angle with $\Sigma$
and agrees with the Euclidean metric $\geuc$ away from the boundary.
In fact we might as well define $\auxm$ on all of $\R^3$.
A simple choice is the metric
  \begin{equation}
    \auxm:=\auxcf^2\geuc
  \end{equation}
with spherically symmetric conformal factor
  \begin{equation}
  \label{Omegadef}
    \auxcf:=\cutoff{\frac{2}{3}}{\frac{1}{3}} \circ \rr
      +\frac{1}{\rr}\cutoff{\frac{1}{3}}{\frac{2}{3}} \circ \rr,
  \end{equation}
so that under $\auxm$
the annular neighborhood $\left\{ \rr \in \left(\frac{2}{3},1\right)\right\}$ of $\partial \Bbar^3$
is isometric to the round cylinder
$\Sph^2 \times \left(\ln \frac{2}{3},0\right)$
and in fact the entire exterior of $\Bbar^3$ is isometric to a half cylinder.

We note that the $\auxm$ unit normal on $\Sigma$
pointing in the same direction as $\nu$
is $(\auxcf \circ X)^{-1}\nu$.
Given any function $\widetilde{u} \in C^2_{loc}(\Sigma)$,
we define the perturbation $X_{\widetilde{u}}: \Sigma \to \R^3$ by $\widetilde{u}$
of the inclusion map $X: \Sigma \to \R^3$ by
  \begin{equation}
  \label{Xu}
    X_{\widetilde{u}}(p):=\exp_{X(p)}^\auxm \widetilde{u}(p)\frac{\nu(p)}{(\auxcf \circ X)(p)},
  \end{equation}
where $\exp^\auxm: T\R^3 \to \R^3$ is the exponential map for $\auxm$ on $\R^3$.

For $\widetilde{u}$ sufficiently small $X_{\widetilde{u}}$ is an immersion with well-defined Euclidean
unit normal $\nu_{\widetilde{u}}$---taken to have positive inner product
with the velocity of the $\auxm$ geodesics generated by $\nu$---and
well-defined scalar mean curvature
relative to $\geuc$ and $\nu_{\widetilde{u}}$,
which we denote by $\Hcal[\widetilde{u}]: \Sigma \to \R$ (so $\Hcal[0]=H$).
Defining also $\Theta[\widetilde{u}]: \partial \Sigma \to \R$ by
  \begin{equation}
    \Theta[\widetilde{u}]:=(\geuc \circ X_{\widetilde{u}})\left(X_{\widetilde{u}}, \, \nu_{\widetilde{u}}\right)
  \end{equation}
(so $\Theta[0]=0$),
our task is to find $\widetilde{u} \in C^\infty(\Sigma)$ solving the system
  \begin{equation}
  \label{firstsys}
    \begin{aligned}
      &\Hcal[\widetilde{u}]=0 \\
      &\Theta[\widetilde{u}]=0
    \end{aligned}
  \end{equation}
and small enough that $X_{\widetilde{u}}$ is an embedding.
It will then follow from the maximum principle
(or directly from the estimates of the construction)
that the image of $X_{\widetilde{u}}$ is contained in $\Bbar^3$
and meets $\partial \Bbar^3$ only along $X_{\widetilde{u}}(\partial \Sigma)$.

A simple symmetry argument, as follows, 
reveals that the requirement $\Theta[\widetilde{u}]=0$
is equivalent to the Neumann condition
$\widetilde{u}_{,\sigma}\left|_{\partial \Sigma}\right.=0$ on $\widetilde{u}$,
but relative to $X_{\widetilde{u}}^*\auxm$ rather than $g=X_{\widetilde{u}}^*\geuc$;
here we recall (see Remark \ref{conormal})
that $\partial_\sigma$ is the inward unit conormal for $\partial \Sigma$.
We start by recalling that under $\auxm$
the neighborhood $\{\rr \geq 2/3\}$ of $\partial \Bbar^3$
in $\R^3$ is isometric to the Riemannian product
$\Sph^2 \times (\ln 2/3,\infty)$,
with $\partial \Bbar^3$ itself a cross section.
Consider now a second copy ${\Bbar^3}'$ of $\Bbar^3$ 
and denote the union of the two copies 
with the two boundaries identified 
by $\underline{\Bbar^3}$. 
Clearly 
$\underline{\Bbar^3} = \Bbar^3 \cup {\Bbar^3}'$ is a differential topological doubling of $\Bbar^3$ along $\partial \Bbar^3$ 
and there is a reflection 
$\Sbar : \underline{\Bbar^3} \to \underline{\Bbar^3}$ 
exchanging the two copies and keeping $\partial \Bbar^3$ pointwise fixed.  
$\Sbar$ is a smooth isometry with respect to the metric $\auxm$ 
(appropriately extended to $\underline{\Bbar^3}$). 
We define also $\underline{\Sigma}:=\Sigma\cup\Sigma'\subset \underline{\Bbar^3}$, 
where $\Sigma':=\Sbar(\Sigma) \subset{\Bbar^3}'$.  

$\underline{\Sigma}$ is clearly $C^0$ by construction and the orthogonality of the  intersection with $\partial \Bbar^3$ 
implies that it is $C^1$. 
Moreover since $\Sigma$ is smooth, 
the reflectional symmetry under $\Sbar$ implies that $\underline{\Sigma}$ is $C^2$. 
Now suppose that $\widetilde{u} \in C^1(\Sigma)$.
We define an extension of $\Sigma_{\widetilde{u}}$ by 
$\underline{\Sigma_{\widetilde{u}}} := \Sigma_{\widetilde{u}} \cup \Sbar(\Sigma_{\widetilde{u}}) $.  
We want to show that $\Sigma_{\widetilde{u}}$
intersects $\partial \Bbar^3$ orthogonally if and only
$\widetilde{u}_{,\sigma}\left|_{\partial \Sigma}\right.=0$.
Since $\Sigma_{\widetilde{u}}$ is clearly $C^1$, 
its extension $\underline{\Sigma_{\widetilde{u}}}$ is $C^1$ if and only if $\Sigma_{\widetilde{u}}$
intersects $\partial \Bbar^3$ orthogonally.
Let 
$\underline{\widetilde{u}}: \underline{\Sigma} \to \R$ 
be the function defined 
by requiring that 
$\underline{\Sigma_{\widetilde{u}}}$ is the graph of $\underline{\widetilde{u}}$ over $\underline{\Sigma}$;  
clearly then $\underline{\widetilde{u}}\left|_{\Sigma}\right.=\widetilde{u}$ 
and $\underline{\widetilde{u}} = \underline{\widetilde{u}} \circ \Sbar$.  
It follows that $\underline{\Sigma_{\widetilde{u}}}$
is $C^1$ precisely when $\underline{\widetilde{u}}$ is,
which in turn holds if and only if $\widetilde{u}_{,\sigma}\left|_{\partial \Sigma}\right.=0$,
establishing our claim.
The next lemma, giving an exact expression for $\Theta$,
provides slightly more information and
explains the above equivalence by a different argument.

\begin{lemma}[The boundary operator]
\label{Thetalemma}
For every $\widetilde{u} \in C^2_{loc}(\Sigma)$ sufficiently small
  \begin{equation}
    \Theta[\widetilde{u}]
    =
    \left.\left[\widetilde{g}\left(\partial_\sigma,\partial_\sigma\right)
       \left(1+\abs{d\widetilde{u}}_{\widetilde{g}}^2\right)\right]^{-1/2}
      \partial_\sigma \widetilde{u}\right|_{\partial \Sigma},
  \end{equation}
where 
$\partial_\sigma\left|_{\partial \Sigma}\right.$ is the inward unit conormal 
(recall Remark \ref{conormal}) 
and 
$\widetilde{g}$ is the metric on $\Sigma$ defined by
  \begin{equation}
    \widetilde{g}(p)
    :=
    X_{\widetilde{u}(p)}^*\auxm\left|_p\right.
    =
    X_{\widetilde{u}}^*\auxm\left|_p\right.-d\widetilde{u} \otimes d\widetilde{u}\left|_p\right..
  \end{equation}
\end{lemma}

\begin{proof}
First we observe that,
since $\auxcf|_{\partial \Bbar^3}=1$,
on $\partial \Sigma$
the vector field $X_u$ is $\auxm$ unit and $\auxm$ orthogonal to $\partial \Bbar^3$
and the Euclidean normals $\nu$ and $\nu_{\widetilde{u}}$ for $X$ and $X_{\widetilde{u}}$
agree with the corresponding $\auxm$ unit normals.
For any $p \in \partial \Sigma$
  \begin{equation}
    \nu_{\widetilde{u}}(p)
    =
    \frac{N(p,\widetilde{u}(p))-\widetilde{u}_cX_{\widetilde{u}(p),d}\widetilde{g}^{cd}(p)}
      {\sqrt{1+\abs{d\widetilde{u}(p)}_{\widetilde{g}(p)}^2}}
    =
    \frac{N(p,\widetilde{u}(p))-dX_{\widetilde{u}(p)}\left(\nabla_{\widetilde{g}}\widetilde{u}\right)(p)}
      {\sqrt{1+\abs{d\widetilde{u}(p)}_{\widetilde{g}(p)}^2}},
  \end{equation}
where $N(p,\widetilde{u}(p))$ is the tangent vector at
$X_{\widetilde{u}(p)}(p)=X_{\widetilde{u}}(p)$
to the $\auxm$ geodesic generated by $\nu(p)$.
Since $N$ is $\auxm$ parallel along these geodesics
and $\partial \Bbar^3$ is $\auxm$ totally geodesic,
we have
  \begin{equation}
    \geuc\left(X_{\widetilde{u}}(p),N(p,\widetilde{u}(p))\right)
    =
    \auxm\left(X_{\widetilde{u}}(p),N(p,\widetilde{u}(p))\right)
    =
    0
  \end{equation}
for every $p \in \partial \Sigma$.
Consequently
  \begin{equation}
    \Theta[\widetilde{u}]
    =
    -\frac{\left(\geuc \circ X_{\widetilde{u}}\right)
      \left(X_{\widetilde{u}}, dX_{\widetilde{u}(\cdot)}\nabla_{\widetilde{g}}\widetilde{u} \right)}
    {\sqrt{1+\abs{d\widetilde{u}}_{\widetilde{g}}^2}}
    =
    -\frac{\left(\auxm \circ X_{\widetilde{u}}\right)
      \left(X_{\widetilde{u}(\cdot)}, dX_{\widetilde{u}(\cdot)}\nabla_{\widetilde{g}}\widetilde{u} \right)}
    {\sqrt{1+\abs{d\widetilde{u}}_{\widetilde{g}}^2}},
  \end{equation}
but $X_{\widetilde{u}(p)}$ is the outward unit conormal for $X_{\widetilde{u}(p)}$
(under both $X_{\widetilde{u}(p)}^*\geuc$ and $X_{\widetilde{u}(p)}^*\auxm$),
so in fact
  \begin{equation}
    \Theta[\widetilde{u}]
    =
    \frac{\partial_{\widetilde{\sigma}}\widetilde{u}}
    {\sqrt{1+\abs{d\widetilde{u}}_{\widetilde{g}}^2}},
  \end{equation}
where $\partial_{\widetilde{\sigma}}$ is the $\widetilde{g}$ inward unit conormal for $\Sigma$.
To finish the proof we need only show that
the direction of the unit conormal is preserved under deformations of constant height:
$\partial_{\widetilde{\sigma}}=\sqrt{\widetilde{g}^{\sigma\sigma}}\partial_\sigma$,
working relative to the coordinates $(s,\sigma)$ on a neighborhood of $\partial \Sigma$.

To this end we note that
$\widetilde{g}(p)=h(1)$,
where $h(t):=X_{t\widetilde{u}(p)}^*\auxm\left|_p\right.$
is the unique solution to the initial-value problem
  \begin{equation}
    \begin{aligned}
      &\ddot{h}_{ab}(t)=\frac{1}{2}\dot{h}_{ac}\dot{h}_{bd}h^{cd}
        + 2\left(R^\auxm_{ijk\ell} \circ X_{t\widetilde{u}(p)}(p)\right)
        N^j(p,t\widetilde{u}(p))N^\ell(p,t\widetilde{u}(p))
        \left.X_{t\widetilde{u}(p),a}^i\right|_p
        \left.X_{t\widetilde{u}(p),b}^k\right|_p \\
      &h_{ab}(0)=\left.g_{ab}\right|_p \\
      &\dot{h}_{ab}(0)=-2\left.A^\auxm_{ab}\right|_p,
    \end{aligned}
  \end{equation}
where $R^\auxm$ is the Riemann curvature of $\auxm$
and $A^\auxm_{ab}$ is the second fundamental form of $X$ relative to $\auxm$.
Obviously $h_{s\sigma}(0)=\left.g_{s\sigma}\right|_p=0$.
Moreover, $\partial \Bbar^3$ is totally geodesic under $\auxm$,
so $\partial_\sigma=-X$ is $\auxm$ parallel along $\partial \Sigma$,
which means $\dot{h}_{s\sigma}(0)=-2A_{s\sigma}^\auxm=0$ too. 
Using again the fact that $\partial \Bbar^3$ is $\auxm$ totally geodesic,
as a consequence of the Codazzi equation we have
  \begin{equation}
    \auxm\left(\left.R^{\auxm}(u,v)w\right|_{X_{t\widetilde{u}(p)}(p)},
      X_{t\widetilde{u}(p)}(p)\right)=0 
  \end{equation}
for any $p \in \partial \Sigma$, $t \in [0,1]$, and $u,v,w \perp X_{t\widetilde{u}(p)}$.
Since
  \begin{equation}
    X_{t\widetilde{u}(p),\sigma}
    =
    \auxm\left(X_{t\widetilde{u}(p),\sigma},X_{t\widetilde{u}(p)}\right)X_{t\widetilde{u}(p)}
      + h_{ss}^{-1}h_{s\sigma}X_{t\widetilde{u}(p),s},
  \end{equation}
we find
  \begin{equation}
    \begin{aligned}
    \ddot{h}_{s\sigma}
    =
    &\frac{1}{2}\dot{h}_{ss}\dot{h}_{s\sigma}h^{ss}
      +\frac{1}{2}\dot{h}_{ss}\dot{h}_{\sigma\sigma}h^{s\sigma}
      +\frac{1}{2}\dot{h}_{s\sigma}^2h^{s\sigma}
      +\frac{1}{2}\dot{h}_{s\sigma}\dot{h}_{\sigma\sigma}h^{\sigma\sigma} \\
    &+2\left(R^\auxm_{ijk\ell} \circ X_{t\widetilde{u}(p)}\right)
        N^j(p,t\widetilde{u}(p))N^\ell(p,t\widetilde{u}(p))X_{t\widetilde{u}(p),s}^iX_{t\widetilde{u}(p),s}^k
        h_{ss}^{-1}h_{s\sigma}.
    \end{aligned}
  \end{equation}
Clearly $h_{s\sigma}(t) \equiv 0$ satisfies this equation
(whatever the values of the diagonal components $h_{ss}$ and $h_{\sigma\sigma}$)
as well as the trivial initial conditions established just above.
Thus $X_{\widetilde{u}(p),\sigma}(p)$ is orthogonal to $X_{\widetilde{u}(p),s}(p)$,
so parallel to $X_{\widetilde{u}(p)}(p)$, for every $p \in \partial \Sigma$,
as claimed.
\end{proof}

Accordingly we may replace the system (\ref{firstsys}) with
  \begin{equation}
    \begin{aligned}
      &\Hcal[\widetilde{u}]=0 \\
      &\left.\partial_\sigma\widetilde{u}\right|_{\partial \Sigma}=0,
    \end{aligned}
  \end{equation}
where the boundary condition is manifestly linear in $\widetilde{u}$.
The mean curvature operator is of course nonlinear;
we denote its linearization at $\widetilde{u}=0$ by $\widetilde{\Lcal}$,
whose relation to the familiar Jacobi operator for $\Sigma$ in $(\R^3,\geuc)$,
  \begin{equation}
    \Lcal=\Delta_g + \abs{A}_g^2,
  \end{equation}
is given by the next lemma,
which also gives yet another condition equivalent to
$\Theta[\widetilde{u}]=0$.

\begin{lemma}[The boundary condition and the linearized operator
relative to the Euclidean metric]
\label{utildeulemma}
For every $\widetilde{u} \in C^2_{loc}(\Sigma)$,
if we define $u \in C^2_{loc}(\Sigma)$ by
  \begin{equation}
  \label{utildeu}
    u:=(\auxcf \circ X)^{-1}\widetilde{u},
  \end{equation}
then
  \begin{enumerate}[(i)]
    \item $\widetilde{\Lcal} \widetilde{u}=\Lcal u$ and
    \item $\left.\partial_\sigma \widetilde{u}\right|_{\partial \Sigma}
      =\left.\left(\partial_\sigma + 1\right)u\right|_{\partial \Sigma}$.
  \end{enumerate}
In particular the boundary condition $\Theta[\widetilde{u}]=0$
is equivalent to the Robin condition $\partial_\sigma u=-u$
along $\partial \Sigma$.
(At this point we remind the reader
(see Remark \ref{conormal}) that $\partial_\sigma$
is the \textbf{inward} unit conormal.)
\end{lemma}

\begin{proof}
The first item is an immediate consequence
of the standard expression for the variation of mean curvature
under a normal deformation,
since the velocity vector field for the deformation (\ref{Xu})
is indeed everywhere orthogonal to $\Sigma$
under $\geuc$ as well as $\auxm$
and has Euclidean magnitude $(\auxcf \circ X)^{-1}\widetilde{u}$.
For the second item we have
  \begin{equation}
    \partial_\sigma \widetilde{u}
    =
    \partial_\sigma[(\auxcf \circ X)u]
    =
    (\auxcf \circ X)u_{,\sigma}-(\auxcf_{,\rr} \circ X)u,
  \end{equation}
since $X_{,\sigma}=-\partial_{\rr}$ along $\partial \Sigma$,
but, recalling (\ref{Omegadef}),
$\auxcf|_{\partial \Bbar^3} \equiv 1$
and $\auxcf_{,\rr}|_{\partial \Bbar^3} \equiv -1$.
\end{proof}

In the next section we study the linearized problem.

\section{The linearized problem}
\label{S:lp}

As noted at the end of the the previous section
we can replace
the nonlinear boundary condition $\Theta[\widetilde{u}]=0$
on $\widetilde{u}$
by the (linear) Robin condition
$\left.(\partial_\sigma+1)u\right|_{\partial \Sigma}=0$
on $u:=\widetilde{u}/(\Omega \circ X)$.
This equivalence followed from
Lemma \ref{Thetalemma} and item (ii) of Lemma \ref{utildeulemma}.
Thus, by item (i) of Lemma \ref{utildeulemma},
to solve the system (\ref{firstsys}) we
are led to study the linearized problem
  \begin{equation}
    \begin{aligned}
      &\Lcal u = E \\
      &(\partial_\sigma+1)u|_{\partial \Sigma}=0,
    \end{aligned}
  \end{equation}
where $E: \Sigma \to \R$ is a prescribed inhomogeneity,
which we are forced to accept
because $\Lcal u$ only approximates $\Hcal[\widetilde{u}]$,
but the boundary condition is homogeneous.

The analysis will be simplified by working with the $\chi$ metric,
so we multiply the first equation above by $\rho^{-2}$
and the second equation by $\rho^{-1}$ to obtain the equivalent system
  \begin{equation}
  \label{linsys}
    \begin{aligned}
      &\Lchi u = \rho^{-2}E \mbox{ on } \Sigma \\
      &\left(\rho^{-1}\partial_\sigma + \rho^{-1}\right)u
        =0 \mbox{ on } \partial \Sigma, 
    \end{aligned}
  \end{equation}
where
  \begin{equation}
    \Lchi:=\rho^{-2}\Lcal=\Delta_\chi+\rho^{-2}\abs{A}_g^2
  \end{equation}
and $\rho^{-1}\partial_\sigma$ is the $\chi$ inward unit conormal for $\Sigma$.

Significantly,
in the boundary operator,
the zeroth-order term $\rho^{-1}$ is small while the first-order term
$\rho^{-1}\partial_\sigma$ is, as measured by the $\chi$ metric, unit,
so we can treat the Robin condition imposed in (\ref{linsys})
as a small perturbation of a Neumann condition.
In fact we will solve (\ref{linsys}),
modulo certain obstructions,
separately on the catenoidal and two planar regions,
where $\Lchi$ also can be treated as a small perturbation of a simple operator,
and through an iterative procedure we will paste together these
``semilocal'' solutions to produce a global one.

\subsection*{Approximate solutions on the catenoidal region}
$\phantom{ab}$
\nopagebreak

To motivate the following proposition we mention now
the following three consequences of Proposition \ref{initest}:
(i) the catenoidal region $\Sigma_{cat}$ with the $\chi$ metric
is approximated by the long half cylinder $\cyl_a^+$ with the flat metric
  \begin{equation}
    \chiK=dt^2+d\vartheta^2,
  \end{equation}
(ii) the operator $\Lchi$ is there approximated by
  \begin{equation}
  \label{Lcat}
    \Lcat=\partial_t^2+\partial_\vartheta^2+2\sech^2 t,
  \end{equation}
and (iii) the Robin operator $\rho^{-1}\partial_\sigma + \rho^{-1}$ is there approximated by
$-(\sgn \vartheta)\partial_\vartheta$.
For each fixed $\zeta$ the length parameter $a$ tends to infinity as $m$ does.

We write $\cyl$ and $\cyl^+$ respectively
for the infinite-length cylinder and half cylinder
  \begin{equation}
    \cyl:=\R \times \Sph^1 \quad \mbox{and} \quad
    \cyl^+:=\cyl \cap \{\abs{\vartheta} \leq \pi/2\},
  \end{equation}
where, as for $\cyl_a$, we make use of the standard coordinates
$(t,\vartheta)$ on the universal cover of $\cyl$.
Additionally we let $\GrpK$ be the two-element subgroup of $\Grp[m]$
preserving $\Sigma_{cat}$ as a set,
which acts on $\cyl^+$ via $\kappa$ in the obvious way,
so that the nontrivial element takes $(t,\vartheta)$
to $(t,-\vartheta)$.
Note that $\GrpK$ then also preserves $\partial \cyl^+$
and that all elements of $\GrpK$
preserve not only the sides (the two unit normals) of $\Sigma_{cat}$ in $\R^3$
but also the sides (the two unit conormals)
of $\partial_{\Sph^2} \Sigma_{cat}$ in $\Sigma$.
Since the functions we are now considering on $\cyl^+$ or on $\partial \cyl^+$
represent either normal deformations or conormal derivatives of normal deformations,
the appropriate action of $\GrpK$ on any such function $f$ is given by
$\mathfrak{g}.f=f \circ \mathfrak{g}^{-1}$ for every $\mathfrak{g} \in \GrpK$.
Finally, to ensure that solutions on $\Sigma_{cat}$ decay
sufficiently rapidly toward its waist,
we include an exponential weight in the norms below.

Specifically, given a nonnegative integer $k$,
reals $\alpha,\gamma \in (0,1)$,
a submanifold $S$ of $\cyl^+$
(always either $\cyl^+$ itself or its boundary $\partial \cyl^+$)
and a function $u: \cyl^+ \to \R$, we define
  \begin{equation}
  \label{catnorms}
    \norm{u}_{k,\alpha,\gamma}
    :=
    \norm{u}_{k,\alpha,\gamma;S}
    :=
    \norm{u: C^{k,\alpha}\left(S,\chiK\left|_{_S}\right.,e^{\gamma \abs{t}}\right)}.
  \end{equation}
Whenever context permits,
we will omit, as indicated, from our notation for these norms
the domain S.
We also define the Banach spaces
  \begin{equation}
      C^{k,\alpha,\gamma}_{\GrpK}(S)
      :=
      \left\{ \left. u \in C^{k,\alpha}\left(S,\chiK\left|_{_S}\right.,e^{\gamma \abs{t}}\right)
        \; \right| \;
        \forall \mathfrak{g} \in \GrpK \;\; u \circ \mathfrak{g} = u \right\}
    \mbox{ for $S=\cyl^+$ or $S=\partial \cyl^+$}.
  \end{equation}

\begin{prop}[Solvability of the model problem on the half catenoids]
\label{catsol}
Let $\alpha,\gamma \in (0,1)$.
There exist a linear map
  \begin{equation}
    \Rcat: C_{\GrpK}^{0,\alpha,\gamma}\left(\cyl^+\right)
      \times C_{\GrpK}^{1,\alpha,\gamma}\left(\partial \cyl^+\right)
    \to C_{\GrpK}^{2,\alpha,\gamma}\left(\cyl^+\right)
  \end{equation}
and a constant $C(\alpha,\gamma)>0$ such that
for any $E \in C_{\GrpK}^{0,\alpha}\left(\cyl^+\right)$
and for any
$f \in C_{\GrpK}^{1,\alpha}\left(\partial \cyl^+\right)$
we have
  \begin{enumerate}[(i)]
    \item $\Lcat \Rcat(E,f) = E$
    \item $-(\sgn \vartheta) \partial_\vartheta \Rcat(E,f)|_{\partial \cyl^+}=f$, and
    \item $\norm{\Rcat(E,f)}_{2,\alpha,\gamma;\cyl^+}
      \leq C(\alpha,\gamma)\left(\norm{E}_{0,\alpha,\gamma;\cyl^+}
        +\norm{f}_{1,\alpha,\gamma;\partial \cyl^+}\right)$.
  \end{enumerate}
\end{prop}

\begin{proof}
We start with data $(E,f)$ as in the statement of the proposition.
Fixing a smooth compactly supported function
$\varphi \in C_c^\infty(\R)$
with integral $\int_{-\infty}^\infty \varphi(t) \, dt = 1$,
we define the functions
$u_{_B}, F: \cyl^+ \to \R$
by
  \begin{equation}
    \begin{aligned}
      &u_{_B}(t,\vartheta)
      :=
      \cutoff{\frac{\pi}{8}}{\frac{\pi}{4}}(\abs{\vartheta}) \cdot
      \left(\frac{\pi}{2}-\abs{\vartheta}\right)\int_{-\infty}^\infty
        \varphi(s)
          f\left(t-\left(\frac{\pi}{2}-\abs{\vartheta}\right)s,\frac{\pi}{2}\right) \, ds
          \mbox{ and} \\
      &F
      :=
      E-\Lcat u_{_B}.
    \end{aligned}
  \end{equation}
Then
  \begin{equation}
    \begin{aligned}
      &-(\sgn \vartheta)\partial_{\vartheta}u_{_B}\left|_{\partial \cyl^+}\right.=f, \\
      &u_{_B} \in C_{\GrpK}^{2,\alpha,\gamma}\left(\cyl^+\right),
        \quad
        \norm{u_{_B}}_{2,\alpha,\gamma} \leq C_1\norm{f}_{1,\alpha,\gamma}, \\
      &F \in C_{\GrpK}^{0,\alpha,\gamma}\left(\cyl^+\right),
        \mbox{ and }
        \norm{F}_{0,\alpha,\gamma} \leq 
          C_1\left(\norm{E}_{0,\alpha,\gamma}+\norm{f}_{1,\alpha,\gamma}\right)
    \end{aligned}
  \end{equation}
for some constant $C_1>0$ independent of the data $(E,f)$.

We will now find $u \in C^{2,\alpha,\gamma}_{\GrpK}\left(\cyl^+\right)$
satisfying 
$\Lcat u = F$ and
$\partial_{\vartheta}u\left|_{\partial \cyl^+}\right.=0$
and appropriately bounded by the data,
so that the proof can be completed by taking
$\Rcat(E,f)=u+u_{_B}$.
To begin we extend $F$, without relabelling, by
even reflection across $\partial \cyl^+$
to a function $F$
of the same name but
defined on the entire cylinder $\cyl$,
satisfying
$\norm{F: C^{0,\alpha}\left(\cyl, \chiK, e^{\gamma \abs{t}}\right)}
\leq C_1 \left(\norm{E}_{0,\alpha,\gamma}+\norm{f}_{1,\alpha,\gamma}\right)$,
and invariant under the two reflections
$(t,\vartheta) \mapsto (t,-\vartheta)$
and $(t,\vartheta) \mapsto (t,\pi-\vartheta)$.
Next we repeat the proof of \cite[Proposition 5.15]{wiygul:stacking}.

For each nonnegative integer $n$ we define the functions $F_n^{\pm}: \R \to \R$ by
  \begin{equation}
    F_n^+(t):=\int_0^{2\pi} F(t,\vartheta) \cos n\vartheta \, d\vartheta
     \quad \mbox{ and } \quad
    F_n^-(t):=\int_0^{2\pi} F(t,\vartheta) \sin n\vartheta \, d\vartheta,
  \end{equation}
but by the symmetries just mentioned $F_n^-(t) \equiv 0$ for every $n$
and $F_n^+(t) \equiv 0$ for every odd $n$,
so that
$F(t,\vartheta)=\frac{1}{2\pi}F_0^+(t)+\frac{1}{\pi}\sum_{n=1}^\infty F_{2n}^+(t) \cos 2n\vartheta$,
at least distributionally.
Now for each nonnegative integer $n \neq 1$ we define $u_n: \R \to \R$ by
  \begin{equation}
    u_n(t)
    :=
    \begin{cases}
      \int_0^t \left[(t-s) \tanh s \tanh t + (\tanh t - \tanh s)\right]F_n(s) \, ds
        \mbox{ for $n=0$} \\
      \frac{n+\tanh t}{2n(1-n^2)}e^{-nt} \int_{-\infty}^t (n-\tanh s) e^{ns} F_n(s) \, ds
        +\frac{n-\tanh t}{2n(1-n^2)}e^{nt} \int_t^\infty (n+\tanh s)e^{-ns} F_n(s) \, ds
        \mbox{ for $n \geq 0$},
    \end{cases}
  \end{equation}
so that $u_n$ solves
$\left(\partial_t^2+2\sech^2 t - n^2\right)u_n=F_n$
with $u_0(0)=\dot{u}_0(0)=0$ and $u_n$ bounded whenever $F_n$
is compactly supported and $n>1$.

Thus the distribution
  \begin{equation}
    u:=\frac{1}{2\pi}u_0+\frac{1}{\pi}\sum_{n=1}^\infty u_{2n}
      \mbox{ solves }
    \Lcat u = f,
  \end{equation}
at least in the distributional sense,
and is even (also as a distribution)
under the reflections $(t,\vartheta) \mapsto (t,-\vartheta)$
and $(t,\vartheta) \mapsto (t,\pi-\vartheta)$.
It is elementary to verify that
  \begin{equation}
    \begin{aligned}
      &\abs{u_n(t)}
      \leq
      \frac{C(\gamma)}{n^2+1}\norm{F}_{0,\alpha,\gamma}e^{\gamma \abs{t}}, \mbox{ so}
      &\norm{u}_{0,0,\gamma}
      \leq
      C(\gamma)\norm{F}_{0,\alpha,\gamma}
    \end{aligned}
  \end{equation}
for some constant $C(\gamma)$ independent of the data $(E,f)$.
Standard elliptic theory, using in particular the Schauder estimates,
then implies that in fact $u$ is a classical solution satisfying
  \begin{equation}
    \norm{u: C^{2,\alpha}\left(\cyl,\chiK,e^{\gamma \abs{t}}\right)}
    \leq
    C_2\left(\norm{E}_{0,\alpha,\gamma}+\norm{f}_{1,\alpha,\gamma}\right)
  \end{equation}
for some constant $C_2>0$ independent of the data $(E,f)$.
Since $u(t,\vartheta)=u(t,\pi-\vartheta)$ for all $\vartheta$,
$\partial_{\vartheta}u|_{\partial \cyl^+}=0$.
Setting $\Rcat(E,f):=u|_{\cyl^+}+u_{_B}$ concludes the proof.
\end{proof}

\subsection*{Approximate solutions on the planar regions}
$\phantom{ab}$
\nopagebreak

Proposition \ref{initest} further implies that for each $i \in \{0,1\}$
(i) each planar region $\Sigma_i$ under the $\chi$ metric
is approximated by an indented copy $\pi_i\Sigma_i$
of the unit disc $\Bbar^2$
under the conformally flat metric
  \begin{equation}
    \chiP{i}:=\rhohat_i^2\left(d\xx^2+d\yy^2\right),
  \end{equation}
for which definition we recall (\ref{rhohat}),
(ii) $\Lchi$ is there approximated by the corresponding Laplacian $\Delta_{\chiP{i}}$,
and (iii) the Robin operator $\rho^{-1}\partial_\sigma+\rho^{-1}$
is there approximated by $-\rhohat_i^{-1}\partial_{\rr}$.
Note that $\chiP{i}$ depends on $m$ but not on $\zeta$ or $\xi$
and that for large $m$ the region $\pi_i\Sigma_i$
tends to the full disc $\Bbar^2$.
Under $\chiP{i}$
the unindented disc $\Bbar^2 \backslash L_i$
(missing only the points in the configuration $L_i$ defined in (\ref{Wmdef}))
resembles a disc of radius $m$ with half cylinders attached
near each point of $L_i$.

We also observe that
$\Sigma_0$ is preserved as a set by every element of $\Grp[m]$,
but the subgroup of $\Grp[m]$ preserving $\Sigma_1$ is precisely
$\Grp^+[m]$, defined above (\ref{Grp+def}).
For convenience we will sometimes write $\GrpP{0}$ for $\Grp[m]$
and $\GrpP{1}$ for $\Grp^+[m]$,
as mentioned above (\ref{Grp+def}).
Every element of $\Grp^+[m]=\GrpP{1}$ preserves the sides of $\Sigma_1$,
but $\Grp[m]=\GrpP{0}$ includes elements that reverse the sides of $\Sigma_0$;
moreover every element of $\Grp^+[m]$ preserves the sides of
$\partial_{\Sph^2} \Sigma_1$ in $\Sigma_1$
and every element of $\Grp[m]$ preserves the sides of
$\partial_{\Sph^2} \Sigma_0$ in $\Sigma_0$.
Since the functions we are considering now on either
$\Bbar^2$ or $\partial \Bbar^2$
represent sections of the normal bundle of various subsets of $\Sigma$
(mean curvature, generators of normal deformations,
or conormal derivatives of generators of normal deformations),
the appropriate action of $\Grp^+[m]$
on such a function $f$ defined on either $\Sigma_1$ or $\partial \Sigma_1$
is given by $\mathfrak{g}.f=f \circ \mathfrak{g}^{-1}$
for every $\mathfrak{g} \in \Grp^+[m]$,
while the appropriate action of $\Grp[m]$
on such a function defined on either $\Sigma_0$ or $\partial \Sigma_0$
is given by
$\mathfrak{g}.f=f \circ \mathfrak{g}^{-1}$
for all $\mathfrak{g}$ preserving the sides of $\Sigma_0$
and by $\mathfrak{g}.f=-f \circ \mathfrak{g}^{-1}$
for all $\mathfrak{g}$ reversing the sides of $\Sigma_0$.
For each $i \in \{0,1\}$, each $\alpha \in (0,1)$,
each nonnegative integer $k$,
and each submanifold $S$ of $\Bbar^2 \backslash L_i$
(in practice always $\Bbar^2 \backslash L_i$ or $\partial \Bbar^2 \backslash L_i)$
we are therefore led to introduce the Banach spaces
  \begin{equation}
      C^{k,\alpha}_{\GrpP{i}}\left(S \backslash L_i,\chiP{i}\right)
      :=
      \left\{ \left. u \in C^{k,\alpha}\left(\Bbar^2 \backslash L_i,\chiP{i}\right) \; \right| \;
         \forall \mathfrak{g} \in \GrpP{i} \;\;
         u \circ \mathfrak{g}
         =\left\langle \mathfrak{g}_*\partial_{\zz}, \partial_{\zz} \right\rangle u
      \right\}
  \end{equation}
along with the abbreviated notation for their norms
  \begin{equation}
    \norm{u}_{k,\alpha}
    :=
    \norm{u}_{k,\alpha:S}
    :=
    \norm{u: C^{k,\alpha}\left(S,\chiP{i}\left|_{_S}\right.\right)}
  \end{equation}
where $S$ is a submanifold of $\Sigma$
(below always either $\Sigma$ or $\partial \Sigma$),
which, as indicated, we will frequently omit from the notation,
and where the choice of $i$ will always be inferred from context.

In order to obtain a bound for the solution independent of
$m$---in spite of the fact that, ignoring the attached cylinders,
$(\Bbar^2,\chiP{1})$ looks like a disc of radius $m$---it is necessary
to assume that away from the boundary the inhomogeneous term is small in terms of $m$.
Specifically, given $\alpha \in (0,1)$ and $E: \Bbar^2 \to \R$
we define the weighted H\"{o}lder norm
  \begin{equation}
  \label{primenorm}
    \norm{E}'_{0,\alpha}
    :=
    \norm{E}'_{0,\alpha,\Sigma}
    :=
    \norm{E:C^{0,\alpha}
      \left(\Bbar^2,\chiP{1},
        \left(\cutoff{\frac{2}{m}}{\frac{1}{m}}+m^{-2}\cutoff{\frac{1}{m}}{\frac{2}{m}}\right)
          \circ (1-\rr)\right)}.
  \end{equation}

Because Proposition \ref{catsol} allows us
to solve the approximate linearized problem with arbitrarily prescribed data
on the entirety of the catenoidal region $\Sigma_{cat}$,
we may assume that the inhomogeneous term and boundary data on each planar region
are supported mostly outside its intersection with $\Grp[m]\Sigma_{cat}$.
In the statement of the proposition below $\supp E$ and $\supp f$
denote the supports of the functions $E$ and $f$ respectively.

This restricted support is helpful because, away from this intersection,
the quantity $m^{-1}\rho$, which will figure in estimates for the solution,
has $C^1$ norm bounded by a constant independent of $m$.
Additionally, the solution, if it can be shown to exist,
will be $\chiP{i}$-harmonic on $\Grp[m]\Sigma_{cat}$,
so (because of the conformality of $\chiP{i}$ to the flat metric)
also harmonic in the classical sense.
This last property will be useful in arranging for the solution to decay
towards $\Grp\Sigma_{cat}$, which will be necessary to ensure
convergence of the iterative procedure by which we paste together a global solution
and also, in view of the exponential tapering of the catenoid,
to ensure that the solution is everywhere small enough to manage the nonlinear terms
and to maintain embeddedness under the final deformation.
To quantify the decay we will weight our H\"{o}lder norms
with the function $m\rhohat_i^{-1}$,
which is constantly $1$ away from the cylindrical regions of
$\left(\Bbar^2 \backslash L_i,\chiP{i}\right)$
and decays exponentially in the length parameter along the half cylinders.

On $\Bbar^2$ the Laplacian $\Delta_{\chiP{1}}$ acting
on $\Grp^+[m]$-equivariant functions with vanishing Neumann data
has one-dimensional kernel, spanned by the constant functions.
To obtain a solution we need to introduce one-dimensional \emph{substitute kernel},
used to adjust the inhomogeneous term to become orthogonal to the constants.
The substitute kernel is spanned by the function
$\what: \Bbar^2 \to \R$ defined by
  \begin{equation}
  \label{what}
    \what
    :=
    \left(
      \cutoff{\frac{1}{20m}}{\frac{1}{15m}}
      \cdot \cutoff{\frac{1}{5m}}{\frac{1}{10m}}
    \right) \circ (1-\rr),
  \end{equation}
so that $\what$ is everywhere nonnegative and supported in the annulus
$\left\{1-\frac{1}{5m}<\rr<1-\frac{1}{20m}\right\}$.
Finally, the fact that $\Delta_{\chiP{1}}1=0$ also plays a helpful role in that
we can adjust the solution by a constant in order to obtain the desired decay
near the indentations.

\begin{prop}[Solvability of the model problem on the top (and bottom) disc]
\label{planar1sol}
Let $\alpha \in (0,1)$ and $m \in \Z^+$.
There exist a linear map
  \begin{equation}
  \begin{aligned}
    \Rplanar{1}[m]:
    &\left\{ \left. E \in C^{0,\alpha}_{\Grp^+[m]}\left(\Bbar^2 \backslash L_1, \chiP{1} \right)
      \; \right| \; \supp E \subset
      \left\{ d\left[L_1,\Bbar^2,\geuc\right]>\frac{1}{30m}\right\} \right\} \\
    &\times
    \left\{ \left. f \in C^{1,\alpha}_{\Grp^+[m]}
      \left(\partial \Bbar^2 \backslash L_1,\rhohat_1^2 d\theta^2\right)
      \; \right| \; \supp f \subset 
      \left\{ d\left[L_1,\Bbar^2,\geuc\right]>\frac{1}{30m}\right\} \right\}
    \to
    C^{2,\alpha}_{\Grp^+[m]}\left(\Bbar^2 \backslash L_1,\chiP{1}\right) \times \R
  \end{aligned}
  \end{equation}
and a constant $C(\alpha)>0$ such that
for any $(E,f)$ in the domain of $\Rplanar{1}[m]$,
if $(u,\mu)=\Rplanar{1}[m](E,f)$, then
  \begin{enumerate}[(i)]
    \item $\Delta_{\chiP{1}}u=E+\mu \what$,
    \item $-\rhohat_1^{-1}\partial_{\rr}u|_{\partial \Bbar^2}=f$, and
    \item $\norm{u:C^{2,\alpha}\left(\Bbar^2 \backslash L_1,\chiP{1},m\rhohat_1^{-1}\right)}
      +\abs{\mu} 
      \leq C(\alpha)\left(\norm{E}'_{0,\alpha;\Bbar^2 \backslash L_1}
        +\norm{f}_{1,\alpha;\partial \Bbar^2 \backslash L_1}\right)$.
  \end{enumerate}
\end{prop}

We will routinely write $\Rplanar{1}$ in place of $\Rplanar{1}[m]$.

\begin{proof}
Suppose $E \in C_{\Grp^+[m]}^{0,\alpha}\left(\Bbar^2 \backslash L_1,\chiP{1}\right)$
and
$f \in C_{c,\Grp^+[m]}^{1,\alpha}\left(\partial \Bbar^2 \backslash L_1,\rhohat_1^2d\theta^2 \right)$
both have support contained in $\left\{ d\left[L_1\right]>\frac{1}{30m} \right\}$
Throughout the proof $C$ denotes a constant,
chosen large enough to validate the estimates asserted,
but independent of $m$, $E$, and $f$.
We fix a smooth, even function $\varphi \in C_c^\infty(\R)$
with support contained in $[-1,1]$
and with integral $\int_{-\infty}^\infty \varphi(s) \, ds =1$,
and we define the functions $u_{_B}, F:\Bbar^2 \backslash L_1 \to \R$ by
  \begin{equation}
    \begin{aligned}
      &u_{_B}(\rr,\theta)
      :=
      \rhohat_1(\rr,\theta) \cdot \cutoff{\frac{2}{m}}{\frac{1}{m}}(1-\rr) \cdot
        (1-\rr)\int_{-\infty}^\infty \varphi(s) f(\theta-(1-\rr)s) \, ds \mbox{ and} \\
      &F:=E-\Delta_{\chiP{1}}u_{_B}.
    \end{aligned}
  \end{equation}

Then
  \begin{equation}
    u_{_B} \in C^{2,\alpha}_{\Grp^+[m]}\left(\Bbar^2 \backslash L_1,\chiP{1}\right),
    \quad
    F \in C^{0,\alpha}_{\Grp^+[m]}\left(\Bbar^2 \backslash L_1,\chiP{1}\right),
    \mbox{ and}
    \quad
    -\rhohat_1^{-1}\partial_{\rr}u_0=f,
  \end{equation}
Moreover, since
  \begin{equation}
  \label{rhohat1est}
    \norm{m^{-1}\rhohat_1: C^1\left(\Bbar^2 
      \cap \left\{d\left[L_1\right]>\frac{1}{100m}\right\},\chiP{1}\right)} \leq C,
  \end{equation}
taking $m$ sufficiently large ensures
  \begin{equation}
  \label{Fprimenorm}
    \begin{aligned}
      &\supp u_{_B} \cup \supp F \subset \left\{ d\left[L_1\right] > \frac{1}{40m}\right\}, \\
      &\norm{u_{_B}: C^{2,\alpha}\left(\Bbar^2 \backslash L_1,\chiP{1},m\rhohat_1^{-1}\right)}
        \leq C\norm{f}_{1,\alpha}, \mbox{ and} \\
      &\norm{F}'_{0,\alpha} \leq C\left(\norm{E}'_{0,\alpha}+\norm{f}_{1,\alpha}\right).
    \end{aligned}
  \end{equation}
To prove the proposition it now suffices to find $\left(u_{_N},\mu\right)$ so that
items (i)-(iii) are satisfied with $u_{_N}$ in place of $u$, $F$ in place of $E$,
and $0$ in place of $f$.

To proceed we set
  \begin{equation}
  \label{mu}
    \mu
    :=
    -\frac{\int_{\Bbar^2} \rhohat_1^2 F \, d\xx \, d\yy}
      {\int_{\Bbar^2} \rhohat_1^2 \what \, d\xx \, d\yy}.
  \end{equation}
From (\ref{rhohat}) it is clear that $\rhohat_1 \geq m$ everywhere
and from (\ref{what}) we can see that
$\int_{\Bbar^2} \what \, d\xx \, d\yy \geq 1/(Cm)$,
so the denominator of (\ref{mu}) is at least $m/C$.
On the other hand, $\rhohat_1 \leq Cm$ on the support of $F$
and by (\ref{primenorm}) and (\ref{Fprimenorm})
$\int_{\Bbar^2} \abs{F} \, d\xx \, d\yy \leq \norm{F}'_{0,\alpha}(1/m^2+C/m)$,
so the numerator of (\ref{mu}) is no greater than $Cm\norm{F}'_{0,\alpha}$.
Thus
  \begin{equation}
  \label{muest}
    \abs{\mu}
      + \norm{F+\mu \what}'_{0,\alpha} 
    \leq C\left(\norm{E}'_{0,\alpha}+\norm{f}_{1,\alpha}\right).
  \end{equation}

Furthermore, it follows immediately from (\ref{mu}) that
$\int_{\Bbar^2} \rhohat_1^2(F+\mu \what) \, d\xx \, d\yy = 0$,
and consequently the classical Poisson equation
  \begin{equation}
  \label{poisson}
    \Delta_{\geuc}\widehat{u}=\rhohat_1^2\left(F+\mu \what\right)
  \end{equation}
(where $\geuc=d\xx^2+d\yy^2$ is now the Euclidean metric on $\Bbar^2$)
has a unique $C^{2,\alpha}$ solution $\widehat{u}$
satisfying $\partial_{\rr}\widehat{u}|_{\partial \Bbar^2}=0$
and $\widehat{u}(0,0)=0$.
Since $\geuc$ is preserved by $\Grp^+[m]$
and both $1$ and $F+\mu \what$ are $\Grp^+[m]$-even,
$\widehat{u}$ must be $\Grp^+[m]$-even as well.

Then, since $\chiP{1}=\rhohat_1^2\geuc$,
$\widehat{u}$ also solves
  \begin{equation}
    \Delta_{\chiP{1}}\widehat{u}=F+\mu\what.
  \end{equation}
By standard elliptic Schauder estimates
and the bounded geometry of $\left(\Bbar^2 \backslash L_1,\chiP{1}\right)$
there exists $\epsilon>0$ such that
for each $p \in \Bbar^2 \backslash L_1$ (including the possibility that $p \in \partial \Bbar^2$),
if $D[p,s]$ denotes the set of all points in $\Bbar^2$ of $\chi$ distance
from $p$ strictly less than $s$, then
  \begin{equation}
  \label{Schauder}
  \begin{aligned}
    \norm{\widehat{u}: C^{2,\alpha}\left(D[p,\epsilon],\chiP{1}\right)}
    &\leq
    C\left(\norm{\widehat{u}: C^0\left(D[p,2\epsilon]\right)}
      +\norm{F+\mu\what: C^{0,\alpha}\left(D[p,2\epsilon]\right)}\right) \\
    &\leq
    C\left(\norm{\widehat{u}: C^0\left(D[p,2\epsilon]\right)}
      +\norm{E}'_{0,\alpha}+\norm{f}_{1,\alpha}\right),
  \end{aligned}
  \end{equation}
where for the second inequality we have made use of
(\ref{primenorm}) and (\ref{muest}).

To estimate the $C^0$ norm of $\widehat{u}$
we study (\ref{poisson}) by separation of variables.
Note that because $\widehat{u}$, $\rhohat_1$, and $F+\mu \what$
are all $\Grp^+[m]$-even, for each $\rr \in [0,1]$
  \begin{equation}
    \begin{aligned}
      &\int_0^{2\pi} \widehat{u}(\rr \cos \theta, \rr \sin \theta) \sin n\theta \, d\theta
        =\int_0^{2\pi} \rhohat_1^2(F+\mu \what)(\rr \cos \theta, \rr \sin \theta)
          \sin n\theta \, d\theta
        =0 \mbox{ for every $n \in \Z^+$, and} \\
      &\int_0^{2\pi} \widehat{u}(\rr \cos \theta, \rr \sin \theta) \cos n\theta \, d\theta
        =\int_0^{2\pi} \rhohat_1^2(F+\mu \what)(\rr \cos \theta, \rr \sin \theta)
          \cos n\theta \, d\theta
        =0 \mbox{ unless $n \in m\Z$}.
    \end{aligned}
  \end{equation}
Defining for each nonnegative integer $n$ the functions
$\widehat{u}_n: [0,1] \to \R$ and $\widehat{F}_n: [0,1] \to \R$ by
  \begin{equation}
  \label{fouriercoeffs}
    \widehat{u}_n(\rr):=\int_0^{2\pi} \widehat{u}(\rr \cos \theta, \rr \sin \theta)
      \cos n\theta \, d\theta
    \quad \mbox{and} \quad
    \widehat{F}_n(\rr):=\int_0^{2\pi} \rhohat_1^2(F+\mu \what)(\rr \cos \theta, \rr \sin \theta)
      \cos n\theta \, d\theta,
  \end{equation}
we find from (\ref{poisson}) that
  \begin{equation}
    \left(\partial_{\rr}^2+\rr^{-1}\partial_{\rr}-\rr^{-2}n^2\right)\widehat{u}_n(r)
    =\widehat{F}_n(r)
  \end{equation}
for each $n$.
Imposing the boundary conditions
$\widehat{u}_n(0)=\partial_{\rr}\widehat{u}_n(1)=0$
yields the solutions
  \begin{equation}
  \label{intrep}
  \begin{aligned}
    &\widehat{u}_0(\rr)
    =
    \int_0^{\rr} \sss^{-1} \int_0^{\sss} \ttt \widehat{F}_0(\ttt) \, d\ttt \, d\sss
        \mbox{ and} \\
    &\widehat{u}_n(\rr)
    =
    -\frac{1}{2n}\left[\left(\rr^{n}+\rr^{-n}\right)
          \int_0^{\rr} \sss^{n+1} \widehat{F}_n(\sss) \, d\sss
        +\rr^n \int_{\rr}^1
          \left(\sss^{1+n}+\sss^{1-n}\right) \widehat{F}_n(\sss) \, d\sss\right]
        \mbox{ for $n \geq 1$}.
  \end{aligned}
  \end{equation}

It now follows easily from (\ref{intrep}),
using (\ref{rhohat1est}), (\ref{muest}), and (\ref{fouriercoeffs}),
that
  \begin{equation}
  \begin{aligned}
    \norm{\widehat{u}_0: C^0([0,1])}
    &\leq
    \int_0^1 \sss^{-1} \int_0^s \ttt \abs{\widehat{F}_0(t)} \, d\ttt \, d\sss \\
    &\leq
    \int_0^{1-\frac{2}{m}} \sss^{-1} \int_0^s \ttt m^2m^{-2}\norm{F}'_{0,\alpha} 
        \, d\ttt \, d\sss \\
      &\;\;\;\; + \int_{1-\frac{2}{m}}^1 \sss^{-1}
        \left(
          \int_0^{1-\frac{2}{m}} \ttt m^2 m^{-2}\norm{F}'_{0,\alpha} \, d\ttt
          +\int_{1-\frac{2}{m}}^1 \ttt m^2\norm{F}'_{0,\alpha} \, d\ttt
        \right) \, d\sss \\
    &\leq
    \frac{1}{4}\norm{F}'_{0,\alpha}
      + \frac{4}{m}\norm{F}'_{0,\alpha}\left(\frac{1}{2}+2m\right) \\
    &\leq
    C\norm{F}'_{0,\alpha}
  \end{aligned}
  \end{equation}
and for $n \geq 3$ and $\rr \in [0,1]$
  \begin{equation}
    \abs{\widehat{u}_n(\rr)}
    \leq
    \frac{1}{n^2-4}\left(\rr^2-\frac{2}{n}\rr^n\right)\norm{\widehat{F}_n: C^0([0,1])},
  \end{equation}
so
  \begin{equation}
    \norm{\widehat{u}_n: C^0([0,1])}
    \leq
    \frac{C}{n^2}m^2\norm{F}'_{0,\alpha},
  \end{equation}
where of course $C$ does not depend on $n$.
Therefore
  \begin{equation}
    \norm{\widehat{u}: C^0\left(\Bbar^2\right)}
    \leq
    \norm{\widehat{u}_0: C^0([0,1])}
      + \sum_{n=1}^\infty \norm{\widehat{u}_{mn}: C^0([0,1])}
    \leq
    C\norm{F}'_{0,\alpha}.
  \end{equation}
Applying this last estimate in conjunction with (\ref{Schauder}) we obtain
  \begin{equation}
    \norm{\widehat{u}: C^{2,\alpha}\left(\Bbar^2 \backslash L_1,\chiP{1}\right)}
    \leq
    C\left(\norm{E}'_{0,\alpha}+\norm{f}_{1,\alpha}\right).
  \end{equation}

It remains to arrange for the solution to decay toward the catenoidal waists.
For this note that because $F+\mu \what$ has support contained in
$\{\rho<40m\}$, the solution $\widehat{u}$ is harmonic
(in the classical sense that $\Delta_{\geuc}\widehat{u}=0$)
on $\{\rho \geq 40m\}$.
By the symmetries it suffices to focus on the component
of $\{\rho \geq 40m \}$ whose closure includes $(1,0)$,
namely the intersection of $\Bbar^2$
and the closed Euclidean disc $D$ with center $(1,0)$
and radius $\frac{1}{40m}$.
By even inversion through $\partial \Bbar^2$
we extend $\widehat{u}$ to
  \begin{equation}
    \overline{u}(\rr \cos \theta, \rr \sin \theta)
    :=
    \begin{cases}
      \widehat{u}(\rr \cos \theta, \rr \sin \theta) \mbox{ for } \rr \leq 1 \\
      \widehat{u}\left(\frac{1}{\rr} \cos \theta, \frac{1}{\rr} \sin \theta \right)
        \mbox{ for } \rr>1,
    \end{cases}
  \end{equation}
whose restriction to $D$ is harmonic.
If $D'$ is the Euclidean disc with center $(1,0)$
and radius $\frac{1}{80m}$,
from the classical theory of harmonic functions we have
  \begin{equation}
    \norm{\overline{u}: C^1(D',\geuc)}
    \leq
    C\norm{\overline{u}: C^0(\partial D)}
  \end{equation}
and therefore
  \begin{equation}
    \norm{\overline{u}-\overline{u}(1,0) : C^0(D',\geuc,\rr)}
    \leq
    C\left(\norm{E}'_{0,\alpha}+\norm{f}_{1,\alpha}\right).
  \end{equation}
The proof is now concluded by taking $u=u_{_B}+\widehat{u}-\widehat{u}(1,0)$.
\end{proof}

On $\pi_0\Sigma_0$ the larger symmetry group $\Grp[m]$ will permit us
to dispense with the $m^{-2}$ interior weight for the inhomogeneous term.
Furthermore, every $\Grp[m]$-odd function on $\Bbar^2$
is $L^2(\geuc)$ orthogonal to the constants,
so there are no obstructions to producing a solution.
On the other hand this also means that the constants are unavailable
to help arrange decay, so for this purpose
we still need to introduce an additional dimension
of \emph{extended substitute kernel}, 
spanned by the function
$\wbarhat: \Bbar^2 \to \R$
defined to be the unique $\Grp[m]$-odd function
having restriction to $W[m]$
  \begin{equation}
    \wbarhat|_{W[m]}
    :=
    \Delta_{\chiP{0}}
      \left(
        \cutoff{\frac{1}{10m}}{\frac{1}{20m}}
        \circ
        \sqrt{(1-\rr)^2+\theta^2}
      \right).
  \end{equation}
Comparing with (\ref{prewbar}) and using (\ref{uncheck}), we see that
  \begin{equation}
    \wbar=\pi_0^*\wbarhat.
  \end{equation}

\begin{prop}[Solvability of the model problem on the middle disc]
\label{planar0sol}
Let $\alpha \in (0,1)$ and $m \in \Z^+$.
There exist a linear map
  \begin{equation}
  \begin{aligned}
    \Rplanar{0}[m]:
    &\left\{ \left. E \in C^{0,\alpha}_{\Grp[m]}\left(\Bbar^2 \backslash L_0, \chiP{0} \right)
      \; \right| \; \supp E \subset
      \left\{ d\left[L_0,\Bbar^2,\geuc\right]>\frac{1}{30m}\right\} \right\} \\
    &\times
    \left\{ \left. f \in C^{1,\alpha}_{\Grp[m]}
      \left(\partial \Bbar^2 \backslash L_0,\rhohat_0^2 d\theta^2\right)
      \; \right| \; \supp f \subset
      \left\{ d\left[L_0,\Bbar^2,\geuc\right]>\frac{1}{30m}\right\} \right\}
    \to
    C^{2,\alpha}_{\Grp[m]}\left(\Bbar^2 \backslash L_0,\chiP{0}\right) \times \R
  \end{aligned}
  \end{equation}
and a constant $C(\alpha)>0$ such that
for any $(E,f)$ in the domain of $\Rplanar{0}[m]$,
if $(u,\mubar)=\Rplanar{0}[m](E,f)$, then
  \begin{enumerate}[(i)]
    \item $\Delta_{\chiP{0}}u=E+\mubar \, \wbarhat$,
    \item $-\rhohat_0^{-1}\partial_{\rr}u|_{\partial \Bbar^2}=f$, and
    \item $\norm{u:C^{2,\alpha}\left(\Bbar^2 \backslash L_0,\chiP{0},m\rhohat_0^{-1}\right)}
      +\abs{\mubar} 
      \leq C(\alpha)\left(\norm{E}_{0,\alpha;\Bbar^2 \backslash L_0}
        +\norm{f}_{1,\alpha;\partial \Bbar^2 \backslash L_0}\right)$.
  \end{enumerate}
\end{prop}

We will routinely shorten $\Rplanar{0}[m]$ to $\Rplanar{0}$.

\begin{proof}
Most of the proof is almost identical to that of Proposition \ref{planar1sol}.
In fact we can follow the proof exactly, making only the obvious modifications,
up until the last step,
where we added an appropriately chosen constant
to arrange for the rapid decay of the solution toward the catenoidal waists.
In particular the obvious analog of $\mu$ in that proof ((\ref{mu}))
will necessarily vanish here,
because of the reflections through lines included in $\Grp^+[m]$.
Likewise the analog of the constant mode $\widehat{u}_0$ must vanish identically,
justifying the use of the unweighted norm of $E$ in item (iii).

In this way we find $\widetilde{u} \in C^{2,\alpha}\left(\Bbar^2 \backslash L_0,\chiP{0}\right)$
satisfying
  \begin{equation}
    \begin{aligned}
      &\Delta_{\chiP{0}}\widetilde{u}=E, \\
      &-\rhohat_0^{-1}\partial_{\rr}\widetilde{u}|_{\Bbar^2}=f, \mbox{ and} \\
      &\norm{\widetilde{u}: C^{2,\alpha}\left(\Bbar^2 \backslash L_0,\chiP{0}\right)}
        \leq
        C\left(\norm{E}_{0,\alpha}+\norm{f}_{1,\alpha}\right).
    \end{aligned}
  \end{equation}
For the last step of arranging the decay
we again write $D$ and $D'$
for the Euclidean discs with common center $(1,0)$
and radii $\frac{1}{40m}$ and $\frac{1}{80m}$ respectively.
As before $\widetilde{u}|_{D \cap \Bbar^2}$ is harmonic
and can be extended to a harmonic function $\overline{u}: D \to \R$.
Now define $\vbarhat: \Bbar^2 \to \R$
to be the unique $\Grp[m]$-odd function having restriction to $W[m]$
  \begin{equation}
    \vbarhat|_{W[m]}:=\cutoff{\frac{1}{10m}}{\frac{1}{20m}} \circ \sqrt{(1-\rr)^2+\theta^2}.
  \end{equation}
Then $\norm{\vbarhat: C^{2,\alpha}\left(\Bbar^2 \backslash L_0,\chiP{0}\right)} \leq C$.
Moreover, since $\vbarhat|_{D}=1$,
$\overline{u}-\left(\overline{u}(1,0)\right)\vbarhat|_D$
is harmonic on $D$ and vanishes at $(1,0)$, so
  \begin{equation}
    \norm{\overline{u}-\left(\overline{u}(1,0)\right)\vbarhat: C^0(D',\geuc,\rr)}
    \leq
    C\left(\norm{E}_{0,\alpha}+\norm{f}_{1,\alpha}\right),
  \end{equation}
and the proof is concluded by setting
$\mubar=-\widetilde{u}(1,0)$
and $u=\widetilde{u}+\mubar \, \vbarhat$.
\end{proof}

\subsection*{Exact global solutions}
$\phantom{ab}$
\nopagebreak

We will now construct and estimate global solutions,
modulo extended substitute kernel,
to the linearized problem on the initial surfaces.
To secure adequate estimates on the solutions
we introduce,
for each $\alpha \in [0,1)$, $\gamma>0$, integer $k \geq 0$,
and submanifold $S$ of $\Sigma$
(always either $\Sigma$ or $\partial \Sigma$)
the weighted H\"{o}lder norms
  \begin{equation}
  \label{globalnorms}
    \begin{aligned}
      &\norm{u}_{k,\alpha,\gamma}
      :=
      \norm{u}_{k,\alpha,\gamma;S}
      :=
      \norm{u: C^{k,\alpha}(S,\chi\left|_{_S}\right., m^\gamma \rho^{-\gamma})} \mbox{ and} \\
      &\norm{E}'_{0,\alpha,\gamma}
      :=
      \norm{E}'_{0,\alpha,\gamma;\Sigma}
      :=
      \norm{E}_{0,\alpha,\gamma;\Sigma}
        +\norm{E|_{\Sigma_1}:C^{0,\alpha}\left(\Sigma_1,\chi,
            \left(\cutoff{\frac{2}{m}}{\frac{1}{m}}+m^{-2}\cutoff{\frac{1}{m}}{\frac{2}{m}}\right)
          \circ (1-\pi_1^*\rr)\right)}
    \end{aligned}
  \end{equation}
for any $u: S \to \R$ and $E: \Sigma \to \R$.
Whenever context permits we will omit from our notation the domain $S$, as indicated. 

We remark that there exists a constant $C>0$
such that for any $k \in \Z \cap [0,2]$, $\alpha,\gamma \in (0,1)$, $\zeta,\xi \in \R$,
function $u: \Bbar^2 \to \R$ supported in $\pi_1(\Sigma_1)$,
function $v: \Sigma \to \R$ supported in $\Sigma_{cat}$,
and function $f: \partial \Sigma \to \R$ supported in $\partial_{\Sph^2} \Sigma_{cat}$,
provided $m$ is chosen sufficiently large in terms of $\zeta$ and $\xi$, we have
  \begin{equation}
    \begin{aligned}
      &C^{-1}\norm{\pi_1^*u}'_{0,\alpha,\gamma;\Sigma}
        \leq \norm{u}'_{0,\alpha;\Bbar^2}
        \leq C\norm{\pi_1^*u}'_{0,\alpha,\gamma;\Sigma}, \\
      &C^{-1}\norm{\kappa^*v}_{k,\alpha,\gamma;\cyl^+}
        \leq e^{\gamma a}\norm{v}_{k,\alpha,\gamma;\Sigma}
        \leq C\norm{\kappa^*v}_{k,\alpha,\gamma;\cyl^+}, \mbox{ and} \\
      &C^{-1}\norm{\kappa^*f}_{k,\alpha,\gamma;\partial \cyl^+}
        \leq e^{\gamma a}\norm{f}_{k,\alpha,\gamma;\partial \Sigma}
        \leq C\norm{\kappa^*f}_{k,\alpha,\gamma;\partial \cyl^+},
    \end{aligned}
  \end{equation}
where the pullbacks are extended to vanish identically where they are otherwise undefined
and we recall (\ref{phi01}), (\ref{catnorms}), and (\ref{primenorm}).

We also define the following H\"{o}lder spaces of $\Grp[m]$-odd functions:
  \begin{equation}
  \label{globalHoelder}
      C^{k,\alpha}_{\Grp[m]}(S)
      :=
      \left\{u \in C^{k,\alpha}(S,\chi\left|_{_S}\right.) \; | \;
        \forall \mathfrak{g} \in \Grp[m] \;\;
        u \circ \mathfrak{g}=\left\langle \mathfrak{g}_*\partial_{\zz},\partial_{\zz}\right\rangle u
      \right\},
  \end{equation}
where $S$ is either $\Sigma$ or $\partial \Sigma$.
For each such $S$ we denote by 
$C^{k,\alpha,\gamma}_{\Grp[m]}(S)$
the vector space $C^{k,\alpha}_{\Grp[m]}(S)$
equipped with the norm $\norm{\cdot}_{k,\alpha,\gamma;S}$,
and we denote by $C'^{,k,\alpha,\gamma}_{\Grp[m]}(\Sigma)$
the vector space $C^{k,\alpha}_{\Grp[m]}(\Sigma)$
equipped with the norm $\norm{\cdot}'_{k,\alpha,\gamma;\Sigma}$.

Finally we recall the definition of $\wbar: \Sigma \to \R$
through (\ref{prewbar}) and (\ref{uncheck}),
we recall the definition of $\what: \Bbar^2 \to \R$ in (\ref{what}),
and we define $w: \Sigma \to \R$
by
  \begin{equation}
    w:=\pi_1^*\what,
  \end{equation}
which we understand to be extended to vanish identically
on $\Sigma \backslash \Grp[m] \Sigma_1$.
When we apply Proposition \ref{globalsol} below
in the final section,
we will always take $f=0$,
in which case item (ii) of the Proposition
simply expresses that the Robin condition
in (\ref{linsys}) is enforced,
which, we recall from the discussion in Section \ref{graphs},
guarantees (without incurring any nonlinear error)
orthogonal intersection of the graph generated by the solution.

\begin{prop}[Solvability of the linearized problem on the initial surface]
\label{globalsol}
Given $\alpha, \gamma \in (0,1)$ and $c>0$
there exists $m_0>0$
such that for every integer $m>m_0$
and for every $\zeta,\xi \in [-c,c]$
there is a linear map
  \begin{equation}
    \Rcal[m,\zeta,\xi]: C^{,0,\alpha}_{\Grp[m]}\left(\Sigma[m,\zeta,\xi]\right)
      \times C_{\Grp[m]}^{1,\alpha}\left(\partial\Sigma[m,\zeta,\xi]\right)
    \to
    C^{2,\alpha}_{\Grp[m]}\left(\Sigma[m,\zeta,\xi]\right) \times \R \times \R
  \end{equation}
and there is a constant $C>0$---depending on $\alpha$ and $\gamma$
but not on $c$, $\zeta$, $\xi$, or $m$---such that
if $(E,f)$ belongs to the domain of $\Rcal[m,\zeta,\xi]$
and $(u,\mu,\mubar)=\Rcal[m,\zeta,\xi](E,f)$, we have
  \begin{enumerate}[(i)]
    \item $\Lchi u = E + \mu w + \mubar \, \wbar$,
    \item $-\rho^{-1}\left(\partial_{\sigma}+1\right)u|_{\partial \Sigma}=f$, and
    \item $\norm{u}_{2,\alpha,\gamma}+\abs{\mu}+\abs{\mubar}
            \leq C\left(\norm{E}'_{0,\alpha,\gamma}+\norm{f}_{1,\alpha,\gamma}\right)$.
    \item Furthermore, for each fixed $m>m_0$,
      the map $(\zeta,\xi) \mapsto \Rcal[m,\zeta,\xi]$
      is continuous. (See Remark \ref{cty}.)
  \end{enumerate}
\end{prop}

\begin{proof}
Let $\alpha,\gamma \in (0,1)$ and $c>0$.
The continuity assertion in item (iv)
will be clear from the construction of $\Rcal[m,\zeta,\xi]$
to follow in light of
the smooth dependence of the geometry of the initial surfaces on the parameters
as expressed in Proposition \ref{initest}.
Suppose now $m \in \Z^+$ and $\zeta,\xi \in [-c,c]$.
To make the proof easier to read
we will henceforth suppress from our notation 
the dependence on $m$, $\zeta$, and $\xi$
of the operator constructed.
We will begin by constructing an approximate solution operator
by pasting together solutions obtained from
$\Rcat$, $\Rplanar{0}$, and $\Rplanar{1}$
and imported to $\Sigma$ via $\kappa$, $\pi_0$, and $\pi_1$
as follows.

First we define the two linear operators (of the same name)
  \begin{equation}
  \label{Psicat}
    \begin{aligned}
      &\Psi_{cat}: C\left(\Sigma_{cat}\right) \to C(\Sigma) \quad \mbox{and} \quad
        \Psi_{cat}: C\left(\partial_{\Sph^2} \Sigma_{cat}\right) 
        \to C(\partial \Sigma) \mbox{ by} \\
      &\Psi_{cat}E:=\kappa^{*-1}\left(\left(\cutoff{a-1}{a-2}\circ \abs{t}\right) 
        \cdot \kappa^*E\right)
    \end{aligned}
  \end{equation}
where we understand the right-hand side to be extended
to the unique $\Grp[m]$-odd function on $\Sigma$ (or $\partial \Sigma$)
which vanishes off $\Grp[m]\Sigma_{cat}$
(or the latter's intersection with $\partial \Sigma$).

Now we define
  \begin{equation}
  \label{Rcaltcat}
    \begin{aligned}
      &\Rcalt_{cat}: C_{\Grp[m]}^{0,\alpha}(\Sigma) \times C_{\Grp[m]}^{1,\alpha}(\partial \Sigma)
        \to C_{\Grp[m]}^{2,\alpha}(\Sigma) \mbox{ by} \\
      &\Rcalt_{cat}(E,f)
      :=
      \Psi_{cat}\kappa^{*-1}v_{cat}, \mbox{ where} \\
      &v_{cat}
      :=
      \Rcat
        \left(\kappa^*\Psi_{cat}\left(E|_{\Sigma_{cat}}\right),
          \kappa^*\Psi_{cat}\left(f|_{\partial_{\Sph^2} \Sigma_{cat}}\right)\right).
    \end{aligned}
  \end{equation}
If $u_{cat}=\Rcalt_{cat}(E,f)$, then, recalling (\ref{Lcat}) and Proposition \ref{catsol},
  \begin{equation}
  \label{LRcat}
    \begin{aligned}
      &\Lchi u_{cat}
      =
      \Psi_{cat}^2\left(E|_{\Sigma_{cat}}\right)
        + \left[\Lchi, \Psi_{cat}\right]\kappa^{*-1}v_{cat}
        +\Psi_{cat}\kappa^{*-1}\left(\kappa^*\Lchi\kappa^{*-1}-\Lcat\right)v_{cat} \mbox{ and} \\
      &\partial_\sigma u_{cat}|_{\partial_\Sigma}
      =
      \Psi_{cat}^2\left(f|_{\partial_{\Sph^2} \Sigma_{cat}}\right)
        +\left.\left(\left[\partial_\sigma,\Psi_{cat}\right]\kappa^{*-1}v_{cat}
        +\Psi_{cat}\kappa^{*-1}\left(\kappa^*\partial_\sigma\kappa^{*-1}
          +(\sgn \vartheta)\partial_\vartheta\right)v_{cat}\right)\right|_{\partial \Sigma}.
    \end{aligned}
  \end{equation}

Next for each $i \in \{0,1\}$ we define the two operators (of the same name)
  \begin{equation}
  \label{Psii}
    \begin{aligned}
      &\Psi_i: C\left(\Sigma_i\right) \to C(\Sigma) \quad \mbox{and} \quad
        \Psi_i: C\left(\partial_{\Sph^2}\Sigma_i\right) \to C(\partial \Sigma) \mbox{ by} \\
      &\Psi_iE
      :=
      \pi_i^*\left(\cutoff{2m^2}{m^2}\circ \rhohat_i\right)\pi_i^{*-1}E,
    \end{aligned}
  \end{equation}
where we understand the right-hand side to be extended to be
the unique $\Grp[m]$-odd function on $\Sigma$ (or on $\partial \Sigma$)
that vanishes identically on
$\Sigma \backslash \Grp[m] \Sigma_i$
(or on $\partial \Sigma \backslash \Grp[m] \partial_{\Sph^2}\Sigma_i$).

For each $E \in C_{\Grp[m]}^{0,\alpha}(\Sigma)$
and $f \in C_{\Grp[m]}^{1,\alpha}(\partial \Sigma)$
we set
  \begin{equation}
  \label{EfP}
    \begin{aligned}
      E_{_P}
      :=
      E-\Psi_{cat}^2\left(E_{\Sigma_{cat}}\right)
        -\left[\Lchi, \Psi_{cat}\right]\kappa^{*-1}v_{cat} \mbox{ and} \\
      f_{_P}
      :=f-\Psi_{cat}^2\left(f|_{\partial_{\Sph^2} \Sigma_{cat}}\right)
        -\left.\left(\left[\partial_\sigma,\Psi_{cat}\right]\kappa^{*-1}v_{cat}\right)
          \right|_{\partial \Sigma},
    \end{aligned}
  \end{equation}
where $v_{cat}$ is as defined in (\ref{Rcaltcat}). 
Then we define
  \begin{equation}
    \begin{aligned}
      &\Rcalt_0: C_{\Grp[m]}^{0,\alpha}(\Sigma) \times C_{\Grp[m]}^{1,\alpha}(\partial \Sigma)
        \to C_{\Grp[m]}^{2,\alpha}(\Sigma) \times \R \mbox{ by} \\
      &\Rcalt_0(E,f)
      :=
      \left(\Psi_0 \pi_0^*v_0,\mubar\right), \mbox{ where} \\
      &(v_0,\mubar):=\Rplanar{0}\left(\pi_0^{*-1}E_{_P}, \pi_0^{*-1}f_{_P}\right),
    \end{aligned}
  \end{equation}
and
  \begin{equation}
    \begin{aligned}
      &\Rcalt_1: C_{\Grp[m]}^{0,\alpha}(\Sigma) \times C_{\Grp[m]}^{1,\alpha}(\partial \Sigma)
        \to C_{\Grp[m]}^{2,\alpha}(\Sigma) \times \R \mbox{ by} \\
      &\Rcalt_1(E,f)
      :=
      \left(\Psi_1 \pi_1^*v_1,\mu\right), \mbox{ where} \\
      &(v_1,\mu):=\Rplanar{1}\left(\pi_0^{*-1}E_{_P}, \pi_0^{*-1}f_{_P}\right).
    \end{aligned}
  \end{equation}

If $(u_0,\mubar)=\Rcalt_0(E,f)$
and $(u_1,\mu)=\Rcalt_1(E,f)$, then, recalling Proposition \ref{planar0sol}
and using definition (\ref{EfP}),
  \begin{equation}
  \label{LR0}
    \begin{aligned}
      &\Lchi u_0
      =
      \Psi_0\left(E_{_P}|_{\Sigma_0}\right)
        +\mubar \, \wbar
        +\left[\Lchi,\Psi_0\right]\pi_0^*v_0
        +\Psi_0 \pi_0^*\left(\pi_0^{*-1}\Lchi \pi_0^*-\Delta_{\chiP{0}}\right)v_0 \mbox{ and} \\
      &\partial_\sigma u_0|_{\partial \Sigma}
      =
      \Psi_0\left(f_{_P}|_{\partial_{\Sph^2}\Sigma_0}\right)
        +
        \left.\left(
          \left[\partial_\sigma,\Psi_0\right]\pi_0^*v_0
          +\Psi_0\pi_0^*\left(\pi_0^{*-1}\partial_\sigma \pi_0^*
            +\rhohat_0^{-1}\partial_{\rr}\right)v_0
        \right)\right|_{\partial \Sigma}
    \end{aligned}
  \end{equation}
and
  \begin{equation}
  \label{LR1}
    \begin{aligned}
      &\Lchi u_1
      =
      \Psi_1\left(E_{_P}|_{\Sigma_1}\right)
        +\mu \, w
        +\left[\Lchi,\Psi_1\right]\pi_1^*v_1
        +\Psi_1 \pi_1^*\left(\pi_1^{*-1}\Lchi \pi_1^*-\Delta_{\chiP{1}}\right)v_1 \mbox{ and} \\
      &\partial_\sigma u_1|_{\partial \Sigma}
      =
      \Psi_1\left(f_{_P}|_{\partial_{\Sph^2}\Sigma_1}\right)
        +
        \left.\left(
          \left[\partial_\sigma,\Psi_1\right]\pi_1^*v_1
          +\Psi_1\pi_1^*\left(\pi_1^{*-1}\partial_\sigma \pi_1^*
            +\rhohat_1^{-1}\partial_{\rr}\right)v_1
        \right)\right|_{\partial \Sigma}.
    \end{aligned}
  \end{equation}

Now we define
  \begin{equation}
    \begin{aligned}
      &L:
      C^{2,\alpha}_{\Grp[m]}(\Sigma) \times \R \times \R
      \to
      C^{0,\alpha}_{\Grp[m]}(\Sigma) \times C^{1,\alpha}_{\Grp[m]}(\partial \Sigma) \mbox{ by} \\
      &L:
      (u,\mu,\mubar)
      \mapsto
      \left(\Lchi u - \mu w - \mubar \, \wbar, 
        \left.\left(\rho^{-1}\partial_\sigma+\rho^{-1}\right)u\right|_{\partial \Sigma}\right)
    \end{aligned}
  \end{equation}
and the approximate solution operator
  \begin{equation}
    \begin{aligned}
      &\Rcalt:
      C_{\Grp[m]}^{0,\alpha}(\Sigma) \times C_{\Grp[m]}^{1,\alpha}(\partial \Sigma)
      \to
      C_{\Grp[m]}^{2,\alpha}(\Sigma) \times \R \times \R \mbox{ by} \\
      &\Rcalt:=\Rcalt_{cat}+\Rcalt_{0}+\Rcalt_{1}.
    \end{aligned}
  \end{equation}
Then, using (\ref{Psicat}), (\ref{LRcat}), (\ref{Psii}), (\ref{EfP}), (\ref{LR0}), and (\ref{LR1}),
  \begin{equation}
      L\Rcalt(E,f)=\left(\widetilde{E},\widetilde{f}\right),
  \end{equation}
where
  \begin{equation}
    \begin{aligned}
      \widetilde{E}
      =
      &E+\Psi_{cat}\kappa^{*-1}\left(\kappa^*\Lchi\kappa^{*-1}-\Lcat\right)v_{cat}
        +\Psi_0 \pi_0^*\left(\pi_0^{*-1}\Lchi \pi_0^*-\Delta_{\chiP{0}}\right)v_0 \\
      &+\Psi_1 \pi_1^*\left(\pi_1^{*-1}\Lchi \pi_1^*-\Delta_{\chiP{1}}\right)v_1
        +\left[\Lchi,\Psi_0\right]\pi_0^*v_0
        +\left[\Lchi,\Psi_1\right]\pi_1^*v_1
    \end{aligned}
  \end{equation}
and
  \begin{equation}
    \begin{aligned}
      \widetilde{f}
      =
      &f+\left(\Psi_{cat}\kappa^{*-1}\left(\kappa^*\partial_\sigma\kappa^{*-1}
          +(\sgn \vartheta)\partial_\vartheta\right)v_{cat}
        +\Psi_0\pi_0^*\left(\pi_0^{*-1}\partial_\sigma \pi_0^*
            +\rhohat_0^{-1}\partial_{\rr}\right)v_0\right. \\
      &\left.\left.+\Psi_1\pi_1^*\left(\pi_1^{*-1}\partial_\sigma \pi_1^*
            +\rhohat_1^{-1}\partial_{\rr}\right)v_1
        +\left[\partial_\sigma,\Psi_0\right]\pi_0^*v_0
        +\left[\partial_\sigma,\Psi_1\right]\pi_1^*v_1
      \right)\right|_{\partial \Sigma}
      +\left.\rho^{-1}\Rcalt(E,f)\right|_{\partial \Sigma}.
    \end{aligned}
  \end{equation}

Using Propositions \ref{initest}, \ref{catsol}, \ref{planar1sol}, and \ref{planar0sol}
along with the fact that each commutator term is supported
in $\{m^2/2 < \rho < 3m^2\}$,
we therefore obtain
  \begin{equation}
  \label{LRest}
    \norm{\widetilde{E}-E}'_{0,\alpha,\gamma}+\norm{\widetilde{f}-f}_{1,\alpha,\gamma}
    \leq
    C\left(m^{-1}+\sech^2 \frac{m}{4} + m^3\tau + m^{\gamma-1} \right)
      \left(\norm{E}'_{0,\alpha,\gamma}+\norm{f}_{1,\alpha,\gamma}\right);
  \end{equation}
we emphasize in particular that to estimate $\norm{\widetilde{E}-E}'_{0,\alpha,\gamma}$
rather than merely $\norm{\widetilde{E}-E}_{0,\alpha,\gamma}$
we have needed item (ii) of Proposition \ref{initest}.
Thus
$L\Rcalt$ is a small perturbation of the identity operator
on $C'^{,0,\alpha,\gamma}_{\Grp[m]}(\Sigma) \times C^{1,\alpha,\gamma}_{\Grp[m]}(\partial \Sigma)$
(recalling the Banach spaces defined immediately below (\ref{globalHoelder}))
and is consequently invertible (continuously in $\zeta$ and $\xi$).
Setting
  \begin{equation}
    \Rcal:=\Rcalt\left(L\Rcalt\right)^{-1}
  \end{equation}
concludes the proof of items (i) and (ii),
and item (iii) is now obvious from the foregoing construction of $\Rcal$
and the estimates afforded by Propositions
(\ref{catsol}), (\ref{planar1sol}), and (\ref{planar0sol}),
ending the proof.
\end{proof}

\section{The main theorem}
\label{end}

\subsection*{The first-order solution}
$\phantom{ab}$
\nopagebreak

Now we can solve our problem (\ref{firstsys}) to first order,
modulo extended substitute kernel.
First we summarize our estimates from Proposition (\ref{initest})
for the initial mean curvature,
making use of the $\norm{\cdot}'_{k,\alpha,\gamma;\Sigma}$ norm
defined in (\ref{globalnorms}).

\begin{lemma}[Estimate of the initial mean curvature]
\label{mcest}
Let $\alpha,\gamma \in (0,1)$.
There is a constant $C>0$
such for any $c>0$
there exists $m_0>0$
such that for every $m>m_0$
and $\zeta,\xi \in [-c,c]$
the mean curvature $H$ of the initial surface $\Sigma[m,\zeta,\xi]$
satisfies the estimate
  \begin{equation}
    \norm{\rho^{-2}H-\xi \tau \wbar}'_{0,\alpha,\gamma;\Sigma}
    \leq
    C\tau.
  \end{equation}
\end{lemma}

\begin{proof}
The estimate follows immediately
from definition (\ref{globalnorms})
and items (\ref{intplanar1est}) and (\ref{Hest})
of Proposition \ref{initest},
using the obvious inequality
$m\rho^{-1} \leq m\rho^{-\gamma}$
as well as the inequality
$m^2\tau \leq m\rho^{-\gamma}$,
which is clear, when $m$ is large,
from (\ref{phi01})
since $\rho^{-\gamma} \geq \tau^\gamma$.
\end{proof}

Now we can apply Proposition \ref{globalsol}.

\begin{lemma}[The solution to first order]
\label{firstordersol}
Let $\alpha,\gamma \in (0,1)$.
There is a constant $C>0$
such that
for any $c>0$ there exists $m_0>0$
such that whenever $m>m_0$ and $\zeta,\xi \in [-c,c]$,
  \begin{equation}
    \norm{u_1}_{2,\alpha,\gamma}+\abs{\mu_1}+\abs{\mubar_1} \leq C\tau
  \end{equation}
where
  \begin{equation}
    (u_1,\mu_1,\mubar_1):=-\Rcal\left(\rho^{-2}H-\tau \xi \wbar,0\right),
  \end{equation}
recalling that $H$ is the mean curvature of the initial surface $\Sigma[m,\zeta,\xi]$
and $\wbar$ is defined through (\ref{prewbar}) and (\ref{uncheck}).
Furthermore, for each fixed $m$
the map $(\zeta,\xi) \mapsto (u_1,\mu_1,\mubar_1)$ 
is continuous
(if we identify $\Sigma[m,\zeta,\xi]$ and $\Sigma[m,0,0]$
as in Remark \ref{cty}).
\end{lemma}

\begin{proof}
Take $m_0$ to be the maximum of the identically named quantities
featuring in the statements of Propositions \ref{initest} and \ref{globalsol}.
The estimate follows from the boundedness of $\Rcal$
stated in Proposition \ref{globalsol}
and from the estimate of $H$ in (\ref{Hest}) of Proposition \ref{initest}.
The continuity claim is clear from the smooth dependence
of the initial surfaces on the parameters
(as in Remark \ref{cty})
and from item (iv) of Proposition \ref{globalsol}.
\end{proof}

\subsection*{The nonlinear terms and the vertical force}
$\phantom{ab}$
\nopagebreak

We recall (\ref{Xu}),
defining the deformation of the inclusion map $X: \Sigma \to \Bbar^3$
by a given function $\widetilde{u}: \Sigma \to \R$,
and we also recall that $\Hcal[\widetilde{u}]$
denotes the mean curvature of this map,
relative to the Euclidean metric $\geuc$ and the unit normal $\nu_{\widetilde{u}}$
obtained by deforming
the unit normal $\nu$ for $\Sigma$.
Of course the initial surface, the inclusion map,
and the deformed inclusion all depend on the data $m$, $\zeta$, and $\xi$,
and so we may write $\Hcal[m,\zeta,\xi,\widetilde{u}]$ to emphasize the corresponding
dependence of the resulting mean curvature.
Now, given $u \in C^2_{loc}(\Sigma)$
we define $\widetilde{u} \in C^2_{loc}(\Sigma)$
by equation (\ref{utildeu})
and we define the nonlinear map
  \begin{equation}
  \label{Qcaldef}
    \begin{aligned}
      &\Qcal[m,\zeta,\xi,\cdot]:
      C^2_{loc}(\Sigma[m,\zeta,\xi])
      \to
      C^0_{loc}(\Sigma[m,\zeta,\xi]) \mbox{ by} \\
      &\Qcal[m,\zeta,\xi,u]
      :=
      \Hcal[m,\zeta,\xi,\widetilde{u}]
        -\Hcal[m,\zeta,\xi,0]
        -\left.\frac{d}{dt}\right|_{t=0}\Hcal[m,\zeta,\xi,t\widetilde{u}]
      =
      \Hcal[\widetilde{u}]-H-\Lcal u.
    \end{aligned}
  \end{equation}
In order to state estimates for $\Qcal$
independent of the parameters $\zeta$ and $\xi$
it is useful to define the quantity
  \begin{equation}
    \taubar:=\taubar[m]:=\tau[m,0,0]=\frac{1}{20m}e^{-m/2},
  \end{equation}
recalling (\ref{phi01}).

\begin{lemma}[Estimate of the nonlinear terms]
\label{quadest}
Let $C,c>0$ and $\alpha,\gamma \in (0,1)$.
There exists $m_0>0$ such that 
  \begin{equation}
    \norm{\rho^{-2}\Qcal[m,\zeta,\xi,u]}'_{0,\alpha,\gamma}
    \leq
    \taubar^{1+\gamma/2}
  \end{equation}
whenever $m>m_0$, $\zeta,\xi \in [-c,c]$,
and $u \in C^{2,\alpha}(\Sigma,\chi)$ satisfies
$\norm{u}_{2,\alpha,\gamma} \leq C\tau$.
Furthermore, for each fixed $m>m_0$, the map
  \begin{equation}
    \begin{aligned}
      &C^{2,\alpha}(\Sigma[m,0,0]) \times \R \times \R
      \to
      C^{0,\alpha}(\Sigma[m,0,0]) \\
      &(u,\zeta,\xi)
      \mapsto
      \Qcal[m,\zeta,\xi,u]
    \end{aligned}
  \end{equation}
is continuous
(where for all $\zeta,\xi \in [-c,c]$
we identify $\Sigma[m,\zeta,\xi]$ with $\Sigma[m,0,0]$
as in Remark \ref{cty}).
\end{lemma}

\begin{proof}
The continuity is obvious from the smooth dependence of the initial surfaces
on the parameters.
Let $u: \Sigma \to \R$ satisfy $\norm{u}_{2,\alpha,\gamma} \leq C\tau$.
We recall (\ref{Xu}) and consider the map $X_{\widetilde{u}}$,
where $\widetilde{u}$ is related to $u$ as in (\ref{utildeu}).
Now take any $p \in \Sigma$ and write $B$ for the $\chi$ geodesic ball in $\Sigma$
with center $p$ and radius $1$.
Clearly for each nonnegative integer $k$
there is a constant $C(k)$, depending on just $k$,
such that the first $k$ $\rho^2(p)\geuc$ covariant derivatives of $\rho^2(p)\auxm$
are, as measured by $\rho^2(p)\geuc$, bounded by $C(k)$.
Moreover, by Proposition \ref{initest},
we can choose $C(k)$ so that it also bounds the first $k$ covariant derivatives
of the second fundamental form in $(\R^3,\geuc)$ of $\rho(p)X|_B$,
the inclusion map of $B$ blown up by a factor of $\rho(p)$.
By scaling it follows that
  \begin{equation}
    \rho^{-1}(p)\norm{\Qcal[u]: C^{0,\alpha}(B,\chi)}
    \leq
    C_1\norm{\rho(p)u: C^{2,\alpha}(B,\chi)}^2
  \end{equation}
for some constant $C_1>0$ independent of $p$, $\zeta$, $\xi$, and $m$,
so
  \begin{equation}
    \rho^{-2}(p)\norm{\Qcal[u]: C^{0,\alpha}(B,\chi)}
    \leq
    C_1\rho(p)\norm{u:C^{2,\alpha}(B,\chi)}^2,
  \end{equation}
and therefore globally
  \begin{equation}
    \begin{aligned}
      \norm{\rho^{-2}\Qcal[u]}'_{0,\alpha,\gamma}
      \leq
      C_1C^2\left(m^3 + m^{\gamma}\tau^{\gamma-1}\right)\tau^2
      \leq
      m\tau^{1+\gamma} \leq me^{4c}\taubar^{1+\gamma}
      \leq
      \taubar^{1+\gamma/2}
    \end{aligned}
  \end{equation}
when $m$ is sufficiently large compared to $c$.
\end{proof}

The Killing field $\partial_{\zz}$
generating vertical translations
induces on the initial surface $\Sigma$ the Jacobi field
$\left(\geuc \circ X\right)\left(\nu,\partial_{\zz} \circ X \right)$,
which can be identified as the geometric origin of the kernel
we confronted when solving the linearized equation on $\Sigma_1$.
To solve that equation we were obliged to introduce substitute kernel,
spanned by the function $w$.
We will manage this kernel by adjusting the parameters $\zeta$ and $\xi$
to control the vertical force
  \begin{equation}
  \label{forcedef}
    \Fcal[m,\zeta,\xi,u]
    :=
    \int_{\Sigma_1 \cup \Grp[m]\kappa(\{t \geq 0\})}
      \left(\geuc \circ X_{\widetilde{u}}\right)
        \left(\nu_{\widetilde{u}},\partial_{\zz} \circ X_{\widetilde{u}}\right)
        \Hcal[\widetilde{u}] \, \sqrt{\abs{g_{\widetilde{u}}}}
  \end{equation}
through the portion of the surface $X_{\widetilde{u}}(\Sigma)$
arising from the deformation under $\widetilde{u}$, as in (\ref{Xu}),
of the deformation under $\pre{u}$, as in (\ref{preXu}),
of the portion of the pre-initial surface $\pre{\Sigma}$
above (and including) the catenoidal waists at height $\zz_0$,
recalling the latter's definition in (\ref{phi01}).
In definition (\ref{forcedef})
the function $\widetilde{u}$ is related to the given function $u$
through (\ref{utildeu}),
the vector field $\nu_{\widetilde{u}}$
is the unit normal for $X_{\widetilde{u}}$ obtained continuously from $\nu$
through the family $\left\{X_{\widetilde{tu}}\right\}_{t=0}^1$,
$g_{\widetilde{u}}$ is the metric on $\Sigma$ induced by $\geuc$
and $X_{\widetilde{u}}$,
and $\sqrt{\abs{g_{\widetilde{u}}}}$ is the corresponding area form.
Of course $\pre{\Sigma}$, $\Sigma$, $\widetilde{u}$, $\nu_{\widetilde{u}}$, and $g_{\widetilde{u}}$
all depend on the data $m$, $\zeta$, and $\xi$, as usual.

\begin{lemma}[Estimate of the vertical force]
\label{forcelemma}
Let $C>0$ and $\alpha,\gamma \in (0,1)$.
There exist $m_0, \overline{c}>0$ such that
for any $c>0$
  \begin{equation}
  \label{forcest}
    \abs{(2\zeta-\xi)-\frac{1}{2\pi\tau}\Fcal[m,\zeta,\xi,u]} \leq \overline{c}
  \end{equation}
whenever $m>m_0$, $\zeta,\xi \in [-c,c]$, and $u \in C^{2,\alpha,\gamma}_{\Grp[m]}(\Sigma)$
satisfies 
$\partial_\sigma \tilde{u}\left|_{\partial \Sigma}\right.=0$ and
$\norm{u}_{2,\alpha,\gamma} \leq C\tau$.
Furthermore, for each fixed $m>m_0$, the map
  \begin{equation}
    \begin{aligned}
      &C^{2,\alpha}_{\Grp[m]}(\Sigma[m,0,0]) \times \R \times \R
      \to
      \R \\
      &(u,\zeta,\xi)
      \mapsto
      \Fcal[m,\zeta,\xi,u]
    \end{aligned}
  \end{equation}
is continuous
(where again we are identifying each $\Sigma[m,\zeta,\xi]$ with $\Sigma[m,0,0]$
as in Remark \ref{cty}).
\end{lemma}

\begin{proof}
The continuity is obvious from the smooth dependence, for each fixed $m$,
of the initial surfaces on the parameters.
To estimate the force we begin with the observation that
the assumption $\widetilde{u}_{,\sigma}\left|_{\partial \Sigma}\right.=0$
means that the surface $X_{\widetilde{u}}$
intersects $\partial \Bbar^3$ orthogonally;
because $\partial_{\zz}$ is Killing,
the formula for the first variation of area implies
  \begin{equation}
    \Fcal[m,\zeta,\xi,u]
    =
    \int_{\Grp[m]\kappa(\{t=0\})}
      \left(\geuc \circ X_{\widetilde{u}}\right)
        \left(\eta_{\widetilde{u}},\partial_{\zz} \circ X_{\widetilde{u}}\right)
        \, ds_{\widetilde{u}}
    +\int_{\partial \Sigma \cap \{\zz \geq \zz_0\}}
      \zz \circ X_{\widetilde{u}} \,\, ds_{\widetilde{u}},
  \end{equation}
where (i) at each point of its domain
$\eta_{\widetilde{u}}: \Grp[m]\kappa(\{t=0\}) \to T\R^3$
is the downward $\geuc$ unit vector which is simultaneously $\geuc$ orthogonal
to $\nu_{\widetilde{u}}$ and $X_{\widetilde{u}}\left(\Grp[m]\kappa(\{t=0\})\right)$
(that is the downward unit conormal at the waists),
(ii) $ds_{\widetilde{u}}$ is the arc length form
induced on the designated curves by $X_{\widetilde{u}}$ and $\geuc$,
and (iii) we have also taken advantage of the equality
$\partial \pre{\Sigma}=\partial \Sigma$
to recognize that
$\partial \Bbar^3 \cap \partial\left(\Sigma_1 \cup \Grp[m]\kappa(\{t \geq 0\})\right)
 =\partial \Sigma \cap \{\zz \geq \zz_0\}$.

We now use (\ref{phi01}), Propositions \ref{preinitest} and \ref{initest},
and the bound assumed on $u$
to make the following estimates
for some constant $C_1>0$ independent of $m$, $\zeta$, $\xi$, and $u$.
On $\{t=0\} \subset \cyl_a^+$ we have the estimates
  \begin{equation}
    \begin{aligned}
      &\abs{\kappa^*ds_{\widetilde{u}}-\tau \, d\vartheta}_{\tau^2 \, d\vartheta^2}
        \leq
        C_1m^{-1}+C_1\tau^{-1}\abs{u \circ \kappa}+C_1\abs{du \circ \kappa}_g^2
        \leq
        C_1m^{-1} \mbox{ and} \\
      &\abs{\eta_{\widetilde{u}} \circ \kappa+\partial_{\zz} \circ \kappa}_{\geuc}
        \leq
        C_1\left(\abs{d\pre{u} \circ \pre{\kappa}}_{\pre{g}}
          +\abs{d\widetilde{u} \circ \kappa}_g\right)
        \leq
        CC_1m^\gamma \tau^\gamma.
    \end{aligned}
  \end{equation}

On
$\partial \Sigma \cap \{\zz \geq \zz_0\}=\partial \pre{\Sigma} \cap \{\zz \geq \zz_0\}$
we have
  \begin{equation}
    \zz=\zz \circ \pre{X}_{\pre{u}}^{-1},
    \qquad
    ds=\pre{X}_{\pre{u}}^*d\pre{s},
    \qquad
    \abs{\zz \circ X_{\widetilde{u}}-\zz} \leq \abs{\widetilde{u}}
      \leq \abs{u} \leq C\tau m^\gamma\rho^{-\gamma},
  \end{equation}
and
  \begin{equation}
    \begin{aligned}
      \abs{ds_{\widetilde{u}}-ds}_{ds^2}
      &\leq
      C_1\abs{A_\auxm}_\auxm \abs{\widetilde{u}}+\abs{d\widetilde{u}}_g^2
        +C_1\left(\abs{R_\auxm}_\auxm+\abs{A_\auxm}^2_\auxm\right)\widetilde{u}^2 \\
      &\leq
      C_1\left(1+\abs{A}_g\right)\abs{u}+\abs{du}_g^2 \\
      &\leq
      C_1C\tau m^\gamma \rho^{-\gamma}
        +C_1\rho^2\tau C\tau m^\gamma \rho^{-\gamma}
        +C_1m^2c\tau C\tau
        + \rho^2C^2\tau^2m^{2\gamma}\rho^{-2\gamma} \\
      &\leq
      CC_1m^\gamma\tau^\gamma,
    \end{aligned}
  \end{equation}
where we recall (\ref{preXu}),
where $R_\auxm$ is the Riemann curvature tensor of $\auxm$
and $A_\auxm$ is the second fundamental form of $\Sigma$
relative to $\auxm$,
and where $ds$ and $d\pre{s}$ are, as earlier, the arc length forms
on $\partial \Sigma$ and $\partial \pre{\Sigma}$ respectively.

On $\{\abs{\vartheta}=\pi/2\} \subset \cyl_a^+$ we have the estimates
  \begin{equation}
    \abs{\pre{\kappa}^*d\pre{s}-\tau \cosh t \, dt}_{\pre{\kappa^*}d\pre{s}^2}
    \leq
    C_1m^{-1}
    \quad \mbox{and} \quad
    \pre{\kappa}^*\zz=\tau t + \zz_0,
  \end{equation}
and finally on $\partial_{\Sph^2} \Sigma_1$ we have the estimates
  \begin{equation}
    \abs{d\pre{s}-\pre{\pi}_1^*d\theta}_{d\pre{s}^2}
    \leq
    C_1m^2\tau
    \quad \mbox{and} \quad
    \abs{\zz-(\tau a + \zz_0)} \leq C_1\tau.
  \end{equation}

From all of the above estimates as well as the elementary integrals
  \begin{equation}
    \begin{aligned}
      &\int_0^a (\tau t + \zz_0) \tau \cosh t \, dt
      =
      \frac{\tanh a}{20m}(\tau a+\zz_0)-\frac{\tau}{20m} + \tau^2
        \mbox{ (using (\ref{phi01}))},  \\
      &\int_0^\pi \tau \, d\vartheta = \tau \pi, \quad \mbox{ and } \quad
      \int_{\arcsin \frac{1}{20m}\sqrt{1-(20m)^{-2}/4}}^{\frac{\pi}{m}}
        (\tau a + \zz_0) \, d\vartheta
      =
      \left(\frac{\pi}{m}-\frac{1}{20m}+O\left(m^{-3}\right)\right)(\tau a + \zz_0)
    \end{aligned}
  \end{equation}
for, respectively, the approximate catenaries, the waists,
and the remaining approximate circle of latitude,
we conclude
  \begin{equation}
    \begin{aligned}
      &\abs{
        \int_{\kappa(\{t=0\})}
          \left(\geuc \circ X_{\widetilde{u}}\left|_{\Grp[m]\kappa(\{t=0\})}\right.\right)
            \left(\eta_{\widetilde{u}},
            \partial_{\zz} \circ X_{\widetilde{u}}
            \left|_{\Grp[m]\kappa(\{t=0\})}\right. \right) \, ds_{\widetilde{u}}
        +\pi\tau
      }
      \leq
      CC_2m^{-1}\tau, \\
      &\abs{
        \int_{\partial_{\Sph^2} \Sigma_{cat} \cap \{\zz \geq \zz_0\}}
          \zz \circ X_{\widetilde{u}} \,\, ds_{\widetilde{u}}
        -\frac{\tau a + \zz_0}{10m}
      }
      \leq
      CC_2m^{-1}\tau, \mbox{ and} \\
      &\abs{
        \int_{\partial_{\Sph^2}\Sigma_1 \backslash \Grp[m]\partial_{\Sph^2}\Sigma_{cat}}
          \zz \circ X_{\widetilde{u}} \,\, ds_{\widetilde{u}}
        -\left(2\pi-\frac{1}{10}\right)\left(\tau a+\zz_0\right)
      }
      \leq
      CC_2\tau
    \end{aligned}
  \end{equation}
for some constant $C_2>0$ independent of $m$, $\zeta$, $\xi$, and $u$,
and therefore, since $\zz_0=\tau \xi + \tau a$, as defined in (\ref{phi01})
  \begin{equation}
    \abs{\Fcal[m,\zeta,\xi,u]-2\pi\tau(2a+\xi-m)}
    \leq
    3CC_2\tau.
  \end{equation}
The result (\ref{forcest}) now follows from the estimate of $a$ in (\ref{phi01}).
\end{proof}

\subsection*{Explicitly defined diffeomorphisms between the initial surfaces}
$\phantom{ab}$
\nopagebreak

As we have already observed,
for each given $m$, the initial surfaces
depend smoothly on the parameters
$\zeta$ and $\xi$.
Consequently for each $\zeta,\xi \in \R$
there exists a diffeomorphism
  \begin{equation}
    P[m,\zeta,\xi]: \Sigma[m,\zeta,\xi] \to \Sigma[m,0,0]
  \end{equation}
depending smoothly on $\zeta,\xi$ in the sense that,
if $X[m,\zeta,\xi]: \Sigma[m,\zeta,\xi] \to \Bbar^3$
is the inclusion map of $\Sigma[m,\zeta,\xi]$ in $\Bbar^3$,
then the composite embedding
$X[m,\zeta,\xi] \circ P[m,\zeta,\xi]^{-1}: \Sigma[m,0,0] \to \Bbar^3$,
(which parametrizes $\Sigma[m,\zeta,\xi]$ over $\Sigma[m,0,0]$)
is smooth in $\zeta$ and $\xi$.

So far we have made use (implicitly) of the identifications such diffeomorphisms afford
only to assert the continuity in $(\zeta,\xi)$ of certain functions
defined on the initial surfaces.
These assertions do not at all depend on the details of the diffeomorphisms used, 
but in the proof of the main theorem below
we will need (modest) control over the effect of these diffeomorphisms
on the $\norm{\cdot}_{2,\alpha,\gamma;\Sigma}$ norms of functions they identify.
To achieve this we now make an explicit choice of $P[m,\zeta,\xi]$.

Actually we first define
  \begin{equation}
    \pre{P}[m,\zeta,\xi]: \pre{\Sigma}[m,\zeta,\xi] \to \pre{\Sigma}[m,0,0],
  \end{equation}
identifying pre-initial surfaces.
As a matter of notation, given $m \in \Z^+$ and $\zeta, \xi \in \R$
we will append $[m,\zeta,\xi]$ to the names for various regions
and maps defined on the pre-initial surfaces
in order to distinguish them, when desirable,
from the corresponding objects arising from a different choice of data.
For each $m \in \Z^+$ we set
  \begin{equation}
    \abar:=a[m,0,0],
  \end{equation}
the value of the quantity $a$, defined in (\ref{phi01}),
when $\zeta=\xi=0$,
and
we will continue to abbreviate $a[m,\zeta,\xi]$ by $a$
whenever the intended values of $\zeta$ and $\xi$
are clear from context.
(Of course $a$ does not depend at all on the $\xi$ parameter.)
The map
  \begin{equation}
    \begin{aligned}
      &T[m,\zeta,\xi]: \cyl_a^+ \to \cyl_{\abar}^+ \mbox{ given by} \\
      &T[m,\zeta,\xi]: (t,\vartheta) \mapsto \left(t \abar/a,\vartheta\right)
    \end{aligned}
  \end{equation}
is then obviously a diffeomorphism.

On the other hand,
for each $i \in \{0,1\}$,
if we set
  \begin{equation}
    \begin{aligned}
      &A_i[m,\zeta,\xi]
      :=
      \pre{\Sigma}_i[m,\zeta,\xi] \backslash
        \left(\Grp[m] \pre{\kappa}[m,\zeta,\xi]\left(\cyl_{a-2}^+\right)\right)
        \subset \pre{\Sigma}_i[m,\zeta,\xi], \\
      &B_i[m,\zeta,\xi]
      :=
      \pre{\Sigma}_i[m,\zeta,\xi] \backslash
        \left(\Grp[m] \pre{\kappa}[m,\zeta,\xi]\left(\cyl_{a-1}^+\right)\right)
        \subset A_i[m,\zeta,\xi], \mbox{ and} \\
      &C_i[m,\zeta,\xi]
      :=
      \pre{\Sigma}_{cat}[m,\zeta,\xi] \cap \left(A_i[m,\zeta,\xi]
        \backslash B_i[m,\zeta,\xi]\right)
        \subset \pre{\Sigma}_{cat}[m,\zeta,\xi] \cap \pre{\Sigma}_i[m,\zeta,\xi],
    \end{aligned}
  \end{equation}
then, provided $m$ is large enough compared to $\zeta$ and $\xi$,
  \begin{equation}
    \begin{aligned}
      &\pre{\pi}_i[m,\zeta,\xi] \left(A_i\right)
        \subset \pre{\pi}_i[m,0,0]\left(\pre{\Sigma}_i[m,0,0]\right)
        \mbox{ and} \\
      &\pre{\pi}_i[m,\zeta,\xi] \left(C_i\right)
        \subset \pre{\pi}_i[m,0,0]\left(\pre{\Sigma}_{cat}[m,0,0]
        \cap \pre{\Sigma}_i[m,0,0]\right),
    \end{aligned}
  \end{equation}
so that the maps
  \begin{equation}
    \begin{aligned}
      &\pre{P}_i[m,\zeta,\xi]: A_i \to \pre{\Sigma}_i[m,0,0] \mbox{ given by} \\
      &\pre{P}_i[m,\zeta,\xi]:= \pre{\pi}_i[m,0,0]^{-1} \circ \pre{\pi}_i[m,\zeta,\xi]
    \end{aligned}
  \end{equation}
and
  \begin{equation}
    \begin{aligned}
      &T_i[m,\zeta,\xi]: \left\{ (t,\vartheta) \in \cyl_a^+ \; | \;
        (-1)^{i+1}t \in [a-2,a-1] \right\} 
      \to 
      \cyl_{\abar}^+ \mbox{ given by} \\ 
      &T_i[m,\zeta,\xi]:=\pre{\kappa}[m,0,0]^{-1}
        \circ \pre{P}_i[m,\zeta,\xi] \circ \pre{\kappa}[m,\zeta,\xi]
    \end{aligned}
  \end{equation}
are well-defined and diffeomorphisms onto their images.

In fact, for large $m$,
$T_i[m,\zeta,\xi]$ is
a small perturbation of the map
$(t,\vartheta) \mapsto \left(t+(-1)^{i+1}\zeta,\vartheta\right)$
and in particular nearly agrees with $T[m,\zeta,\xi]$
for $\abs{t}$ close to $a$.
We will now glue these maps together into a single smooth diffeomorphism
by taking convex combinations, relative to the flat metric,
of their image points on $\cyl_{\abar}^+$.
Namely we define
  \begin{equation}
    \begin{aligned}
      &\widetilde{T}[m,\zeta,\xi]: \cyl_{a-1}^+ \to \cyl_{\abar}^+ \mbox{ by} \\
      &\widetilde{T}[m,\zeta,\xi]
      :=
      \left(\cutoff{-a+2}{-a+1} \circ t\right)T_0[m,\zeta,\xi]
         + \left(\cutoff{a-2}{a-1} \circ t\right)T_1[m,\zeta,\xi] \\
      &\;\;\;\;\;\;\;\;\;\;\;\;\;\;\;\;\;\;\;  + \left(\cutoff{-a+1}{-a+2} \circ t\right)
            \left(\cutoff{a-1}{a-2} \circ t\right)T[m,\zeta,\xi],
    \end{aligned}
  \end{equation}
where scalar multiplication and addition on $\cyl_{\abar}^+$
are defined via the $(t,\vartheta)$ coordinates componentwise
and we remark that each $T_i[m,\zeta,\xi]$ appearing in the definition
is defined on the support of the cutoff function multiplying it.

Finally we define
$\pre{P}[m,\zeta,\xi]: \pre{\Sigma}[m,\zeta,\xi] \to \pre{\Sigma}[m,0,0]$
to be the unique $\Grp[m]$-equivariant map having the restrictions
  \begin{equation}
    \begin{aligned}
      &\pre{P}[m,\zeta,\xi]\left|_{B_i[m,\zeta,\xi]}\right.
      :=
      \pre{P}_i[m,\zeta,\xi] \mbox{ for each $i \in \{0,1\}$ and} \\
      &\pre{P}[m,\zeta,\xi]\left|_{\pre{\kappa}[m,\zeta,\xi]\left(\cyl_{a-1}\right)}\right.
      :=
      \pre{\kappa}[m,0,0] \circ \widetilde{T}[m,\zeta,\xi] \circ \pre{\kappa}[m,\zeta,\xi]^{-1}.
    \end{aligned}
  \end{equation}

This concludes the definition of $\pre{P}[m,\zeta,\xi]$.
We then set
  \begin{equation}
    P[m,\zeta,\xi]:=\pre{X}_{\pre{u}[m,0,0]}[m,0,0] \circ
        \pre{P}[m,\zeta,\xi] \circ \pre{X}_{\pre{u}[m,\zeta,\xi]}[m,\zeta,\xi]^{-1},
  \end{equation}
where we recall from (\ref{preXu}) the map
$\pre{X}_{\pre{u}}: \pre{\Sigma} \to \Sigma$
defining the initial surface and identifying it with its pre-initial precursor.
Last we define the isomorphism
  \begin{equation}
  \label{Pcaldef}
    \Pcal_{m,\zeta,\xi}:=P[m,\zeta,\xi]^*,
  \end{equation}
taking functions on $\Sigma[m,0,0]$ (or $\partial \Sigma[m,0,0]$)
to functions on $\Sigma[m,\zeta,\xi]$ (or $\partial \Sigma[m,0,0]$).

For use in the proof of the main theorem
we observe,
referring to (\ref{uncheck})
and especially the definition of $\pre{P}[m,\zeta,\xi]$ above,  
that there is a constant $C>0$
such that for
each real $c>0$
there exists $m_0>0$ such that
for all $m>m_0$,
$u \in C^{2,\alpha}_{\Grp[m]}(\Sigma[m,\zeta,\xi])$
and $v \in C^{2,\alpha}_{\Grp[m]}(\Sigma[m,0,0])$
  \begin{equation}
  \label{Pest}
    \begin{aligned}
      &\norm{\Pcal_{m,\zeta,\xi}^{-1}u}_{C^{2,\alpha,\gamma}(\Sigma[m,0,0])}
        \leq 
        Ce^{c}\norm{u}_{C^{2,\alpha,\gamma}(\Sigma[m,\zeta,\xi])} \mbox{ and} \\
      &\norm{\Pcal_{m,\zeta,\xi}v}_{C^{2,\alpha,\gamma}(\Sigma[m,\zeta,\xi])}
        \leq
        Ce^{c}\norm{v}_{C^{2,\alpha,\gamma}(\Sigma[m,0,0])}.
    \end{aligned}
  \end{equation}
We mention that more refined estimates are available,
but the above will suffice for our purposes.

\subsection*{The main theorem}
$\phantom{ab}$
\nopagebreak

\begin{theorem}[Main Theorem]
\label{mainthm}
Fix $\alpha,\gamma \in (0,1)$.
There exist $m_0,\cbar,C>0$ such that
for every integer $m>m_0$
there are parameters $\zeta,\xi \in [-\cbar,\cbar]$
and a function $u \in C^\infty(\Sigma[m,\zeta,\xi])$
such that $\norm{u}_{2,\alpha,\gamma} \leq C\tau$
and $X_{\widetilde{u}}: \Sigma[m,\zeta,\xi] \to \Bbar^3$ (as defined in (\ref{Xu}))
is a smooth embedding whose image
$X_{\widetilde{u}}(\Sigma[m,\zeta,\xi])$ is a
two-sided
free boundary minimal surface
in $\Bbar^3$,
invariant under $\Grp[m]$ and
having connected boundary and genus $m-1$.
The minimal surfaces obtained tend in the appropriate sense as $m\to\infty$ 
to the equatorial disc with multiplicity three
and the length of their boundary tends to $6\pi$. 
\end{theorem}

\begin{proof}
For all $\zeta,\xi \in \R$ and $m$ sufficiently large
we define $(u_1,\mu_1,\mubar_1)$ as in Lemma \ref{firstordersol}.
For each $m$ we set
  \begin{equation}
    B_m:=\left\{ u \in C^{2,\frac{\alpha}{2}}_{\Grp[m]}(\Sigma[m,0,0]) 
               \; | \; \norm{u}_{2,\alpha,\gamma} \leq \taubar^{1+\gamma/3} \right\},
  \end{equation}
and for each $v \in B_m$ and $m$ sufficiently large let
  \begin{equation}
    (v',\mu_v,\mubar_v)
    :=
    -\Pcal_{m,\zeta,\xi}^{-1}
      \Rcal\left(\rho^{2}\Qcal[m,\zeta,\xi,u_1+\Pcal_{m,\zeta,\xi}v],0\right),
  \end{equation}
where $\Rcal$ is the solution (modulo extended substitute kernel)
operator introduced in Proposition \ref{globalsol},
$\Qcal[\cdots]$
represents the part, defined in (\ref{Qcaldef}), of the mean curvature which is nonlinear in 
the perturbing function,
and
$\Pcal$ is the map defined in (\ref{Qcaldef}) used to identify functions
defined on initial surfaces with common value of $m$ but possibly different values
of $\zeta$ and $\xi$.
According to Proposition \ref{globalsol}, Lemma \ref{firstordersol},
and Lemma \ref{quadest}, using also (\ref{Pest}),
we can choose $C>0$ large enough so that for any $c>0$
there exists $m_0$
sufficiently large so that for all $m>m_0$ and $\zeta,\xi \in [-c,c]$
  \begin{equation}
    \begin{aligned}
      &\norm{u_1}_{2,\beta,\gamma}+\abs{\mu_1}+\abs{\mubar_1} \leq C\tau \mbox{ and} \\
      &\norm{v'}_{2,\beta,\gamma}+\abs{\mu_v}+\abs{\mubar_v} \leq \taubar^{1+\gamma/3}.
    \end{aligned}
  \end{equation}
(Here we have used the fact, apparent from (\ref{phi01}),
that by taking $m_0$ large enough
any power of $\tau$ with negative exponent
can be made to dominate any given constant,
including the factor of $e^c$ appearing in  (\ref{Pest}).
We remark that one could alternatively establish
a uniform bound for
$\Pcal^{-1}_{m,\zeta,\xi}\Rcal[m,\zeta,\xi]
\left(\Pcal_{m,\zeta,\xi} \times \Pcal_{m,\zeta,\xi}\right)$;
this is possible but not asserted in Proposition \ref{globalsol}.)

In particular, for $m$ sufficiently large,
  \begin{equation}
    \abs{\frac{\mubar_1+\mubar_v}{\tau}} \leq 2C,
  \end{equation}
so by Lemma \ref{forcelemma}
  \begin{equation}
    \abs{\zeta-\frac{1}{2}\xi+\frac{\mubar_1+\mubar_v}{2\tau}
    -\frac{1}{4\pi\tau}\Fcal[m,\zeta,\xi,u_1+\Pcal_{m,\zeta,\xi}v]}
    \leq
    \frac{1}{2}\overline{c}+C
  \end{equation}

Now set
  \begin{equation}
    \cbar:=\max \left\{2C,\frac{1}{2}\overline{c}+C\right\}
  \end{equation}
and take $m_0$ to be the maximum of the quantities of the same name
in Proposition \ref{preinittop}, Proposition \ref{globalsol},
Lemma \ref{firstordersol}, Lemma \ref{quadest}, and Lemma \ref{forcelemma}
when the quantity $c$ in each is taken to be $\cbar$.
Suppose $m>m_0$ and define
  \begin{equation}
    \begin{aligned}
      &\Jcal:
      B_m \times [-\cbar,\cbar] \times [-\cbar,\cbar]
      \to
      B_m \times [-\cbar,\cbar] \times [-\cbar,\cbar] \mbox{ by} \\
      &\Jcal: \begin{bmatrix}v \\ \zeta \\ \xi\end{bmatrix}
      \mapsto
      \begin{bmatrix}
        -v' \\
        \zeta-\xi/2+\left(\mubar_1+\mubar_v\right)/(2\tau)
          -\Fcal[m,\zeta,\xi,u_1+\Pcal_{m,\zeta,\xi}v]/(4\pi\tau) \\
        -\left(\mubar_1+\mubar_v\right)/\tau
      \end{bmatrix},
    \end{aligned}
  \end{equation}
where $\Fcal[\cdots]$ is the vertical force defined in (\ref{forcedef});
that $\Jcal$ really maps into the stated target follows from the preceding inequalities
and the choice of $\cbar$.
Moreover, by the continuity assertions in Proposition \ref{globalsol},
Lemma \ref{firstordersol}, Lemma \ref{quadest}, and Lemma \ref{forcelemma},
$\Jcal$ is continuous, using the $C^{2,\frac{\alpha}{2}}$ norm on the first factor.
Since this first factor is compact in $C^{2,\frac{\alpha}{2}}(\Sigma)$,
$B_m$ is compact, and it is obviously convex.

The Schauder fixed-point theorem therefore applies to ensure
the existence of a fixed point $(v,\zeta,\xi)$ to the map $\Jcal$.
Now, setting $u_2:=\Pcal_{m,\zeta,\xi}v$
and defining $\widetilde{u_1}$ and $\widetilde{u_2}$
in terms of $u_1$ and $u_2$ as in (\ref{utildeu}),
  \begin{equation}
    \begin{aligned}
      \rho^{-2}\Hcal[\widetilde{u_1}+\widetilde{u_2}]
      &=
      \rho^{-2}H+\Lchi u_1+\Lchi u_2 + \rho^{-2}\Qcal[u_1+u_2] \\
      &=
      \left(\mu_1 + \mu_v\right)w,
    \end{aligned}
  \end{equation}
where to get the second line
we have used the definition of $u_1$ through Lemma \ref{firstordersol}
and the fact that $\Jcal(v,\zeta,\xi)=(v,\zeta,\xi)$.
By this same fact though we also conclude that $\Fcal=0$.
On the other hand, $w$ has a sign and the unit normal to $X_{\widetilde{u_1}+\widetilde{u_2}}$
has a positive vertical component on the support of $w$,
so from the above expression for the mean curvature of $X_{\widetilde{u_1}+\widetilde{u_2}}$
and the definition of $\Fcal$
we see that in fact $\mu_1+\mu_v=0$,
so $X_{\widetilde{u_1}+\widetilde{u_2}}$ is exactly minimal.
The remaining claims follow from the estimate for $u:=u_1+u_2$,
Proposition \ref{preinittop},
and the $\Grp[m]$-equivariance of $\Hcal$ and $\Rcal$.
\end{proof}

\bibliographystyle{amsplain}
\bibliography{paper}
\end{document}